\documentclass[10pt,a4paper]{article}

\usepackage{fullpage}
\usepackage{amsmath, amsfonts, amssymb, amsbsy}
\usepackage{latexsym} \usepackage{stmaryrd} 

\usepackage[utf8]{inputenc}

\usepackage{pifont} 

\usepackage{hyperref}
\hypersetup{
    colorlinks,%
    citecolor=black,%
    filecolor=black,%
    linkcolor=black,%
    urlcolor=black
}

\begin{document}


\newtheorem{theorem}{Theorem}[section]
\newtheorem{proposition}[theorem]{Proposition}
\newtheorem{corollary}[theorem]{Corollary}
\newtheorem{lemma}[theorem]{Lemma}
\newtheorem{definition}[theorem]{Definition}
\newtheorem{remark}[theorem]{Remark}
\newtheorem{example}[theorem]{Example}
\newtheorem{conjecture}[theorem]{Conjecture}

\newcommand{\proof}{\noindent \textbf{Proof. }}
\newcommand{\qed}{ \hfill {\vrule width 6pt height 6pt depth 0pt} \medskip }

\newcommand{\separe}{\medskip \centerline{\tt -------------- } \bigskip}
\newcommand{\topic}[1]{\medskip \centerline{\tt --------------- #1 ---------------} \bigskip}
\newcommand{\note}[1]{\medskip \noindent \fbox{\parbox{\textwidth}{\small \tt nota:\,#1}} \smallskip}

\renewcommand{\separe}{}
\renewcommand{\topic}[1]{}

\newcommand{\todo}[1]{\hfill {\small \tt [todo:\,#1]}}


\newcommand{\eps}{\varepsilon} 
\renewcommand{\det}{\mathrm{det}} 
\renewcommand{\div}{\mathrm{div}} 
\newcommand{\argmin}{ \mathrm{argmin} \,}
\newcommand{\Om}{\Omega} 
\newcommand{\weakto}{ \rightharpoonup}  
\newcommand{\weakstarto}{\stackrel{*}{\rightharpoonup}}
\newcommand{\R}{\mathbb{R}}

\newcommand{\F}{\mathcal{F}}
\newcommand{\E}{\mathcal{E}}

\newcommand{\stress}{\boldsymbol{\sigma}} 
\newcommand{\strain}{\boldsymbol{\epsilon}} 
\newcommand{\jump}[1]{\llbracket #1 \rrbracket}
\def\Xint#1{\mathchoice
 {\XXint\displaystyle\textstyle{#1}}%
 {\XXint\textstyle\scriptstyle{#1}}%
 {\XXint\scriptstyle\scriptscriptstyle{#1}}%
 {\XXint\scriptscriptstyle\scriptscriptstyle{#1}}%
 \!\int}
\def\XXint#1#2#3{{\setbox0=\hbox{$#1{#2#3}{\int}$}
 \vcenter{\hbox{$#2#3$}}\kern-.5\wd0}}
 \def\ddashint{\Xint=} 
 \def\dashint{\Xint-}


\newcommand{\flatto}{{\hspace{1pt}\flat}}
\newcommand{\sharpo}{{\hspace{0.5pt}\sharp}}
\newcommand{\0}{{\text{\tiny 0}}}  \newcommand{\1}{{\text{\tiny 1}}} \newcommand{\2}{{\text{\tiny 2}}}

\renewcommand{\1}{{\hspace{-3pt}\text{\fontsize{6}{6} \rm 1}\hspace{0.2pt}}}
\renewcommand{\2}{{\hspace{-3pt}\text{\fontsize{6}{6} \rm  2}\hspace{0.2pt}}}

\newcommand{\U}{\mathcal{U}}\newcommand{\V}{\mathcal{V}}\newcommand{\D}{\mathcal{D}}

\renewcommand{\;}{\hspace{0.5pt}}
\renewcommand{\div}{\text{\sl div}}
\newcommand{\C}{\, \boldsymbol{C}}

\newcommand{\Chat}{\,\widehat{\boldsymbol C}}

\newcommand{\Gc}{G^{\text{\sl \scriptsize c}}}
\renewcommand{\Gc}{G^{\hspace{-2.5pt}\text{\sl \fontsize{6}{6}  c}}}
\newcommand{\Ghat}{\hspace{3pt}\widehat{\hspace{-3pt}\G}} 

\newcommand{\A}{{\hspace{-4pt}\text{\fontsize{6}{6} $A$}}}	\newcommand{\B}{{\hspace{-3.5pt}\text{\fontsize{6}{6}  $B$}}}

\newcommand{\GcA}{\Gc_\A}	\newcommand{\GcB}{\Gc_\B}
\newcommand{\CA}{\boldsymbol{C}_\A}	\newcommand{\CB}{\boldsymbol{C}_\B}

\renewcommand{\hom}{{\hspace{-3pt}\text{\fontsize{7}{7} \sl hom}\hspace{0.2pt}}}
\newcommand{\smahom}{{\hspace{-3pt}\text{\fontsize{4}{4} \sl hom}\hspace{0.2pt}}}
\newcommand{\eff}{{\hspace{-3.5pt}\text{\fontsize{7}{7} \sl ef\hspace{0.15pt}f}\hspace{1pt}}}
\newcommand{\loc}{{\hspace{-3.5pt}\text{\fontsize{7}{7} \sl loc}\hspace{1pt}}}
\newcommand{\Ghom}{G_\hom}		\newcommand{\Gceff}{\Gc_\eff}
\newcommand{\Gchom}{\Gc_\hom}

\newcommand{\W}{\mathcal{W}}

\newcommand{\G}{\hspace{0.5pt}\mathcal{G}}
\newcommand{\hatU}{\hat\U}%
\newcommand{\hatC}{\,\hat{\boldsymbol C}}
\newcommand{\h}{{\hom}}

\renewcommand{\l}{{\hspace{0.5pt}l}}
\newcommand{\evhom}{\text{\sl ev-hom}}



\vspace{0.1cm}
{\Large
\noindent {\bf Homogenization of Griffith's Criterion for brittle Laminates} }
\vspace{6pt} 

\vspace{32pt}

\begin{small}
\noindent{\bf M.~Negri}

\noindent{Department of Mathematics -  University of Pavia} 

\noindent{Via A.~Ferrata 5 - 27100 Pavia - Italy}

\noindent{matteo.negri@unipv.it}

\vspace{32pt}
\noindent {\bf Abstract.} We consider a periodic, linear elastic laminate with a brittle crack, evolving along a prescribed path according to Griffith's criterion. We study the homogenized limit of this evolution, as the size of the layers vanishes. The limit evolution is governed again by Griffith's criterion, in terms of the energy release (of the homogenized elastic energy) and an effective toughness,
which in general differs from the weak$^*$ limit   of the periodic toughness. We provide a variational characterization of the effective toughness and, by the energy identity, we link the toughening effect (in the limit) to the micro-instabilities of the evolution (in the periodic laminate). Finally, we provide a couple of explicit calculations of the effective toughness in the anti-plane setting, showing in particular an example of toughening by  elastic contrast. 


\end{small}

\section*{Introduction}

From the mathematical point of view the homogenization of evolution equations is not a simple development of the (well established) homogenization theory for elliptic equations. 
A remarkable example was presented by Tartar in \cite{Tartar_ARMA90}. Consider a Cauchy problem of the form 
$$
	\begin{cases} \dot{u}_\eps ( x , t ) + a_\eps (x) u _\eps ( x, t) = f( x, t )  \\ 
	u_\eps ( x, 0 ) = u^0 (x)  .
	\end{cases}	
$$
If $a_\eps \weakstarto a_0$  (which is the case in homogenization) then, up to subsequences, $u_\eps \weakto u_0$ where $u_0$ solves an evolution equation of the form 
$$
	\dot{u}_0 ( x , t ) + a_0 (x) u_0 ( x, t) + \int_0^t K ( x, t- \tau) u (\tau) \, d\tau = f( x, t )  . 
$$
Therefore, a memory effect (in time) appears in the limit, as a consequence of the heterogeneity (in space) of the coefficients. This  example is studied also in \cite[\S  3.5.2]{Mielke_LNAMM16} from a different perspective: setting the problem in spaces of measures Mielke has shown that different limit evolution may occur, as the system is endowed with different energy-dissipation structures. For instance, if $a_\eps (x) = a ( x / \eps)$ and $f=0$, then suitable choices of energy and dissipation lead to the following limit equations:
$$
	\dot{u}_0 ( x , t ) + ( \min a ) \, u_0 ( x, t) = 0 , 
		\qquad
	\dot{u}_0 ( x , t ) +  ( \max a) \, u_0 ( x, t) = 0 
$$
(clearly, different subsequences are extracted from $u_\eps$). In \cite{Mielke_LNAMM16} this examples is presented in the context of evolutionary $\Gamma$-convergence, which gives a notion of convergence for a class of doubly non-linear parabolic evolutions of the form 
$$
	\begin{cases} 
	\partial_{v} \mathcal{R}_\eps ( u_\eps  (x,t) , \dot{u}_\eps (x,t) ) = - \nabla \E_\eps (  u_\eps (x,t) )  ,  \\
	u_\eps ( 0, x) = u^\eps (x) ,
	\end{cases}
$$
usually when $\E_\eps$ $\Gamma$-converge to $\E_0$ (for the general theory of $\Gamma$-convergence we refer to the book of Dal~Maso \cite{DalMaso93}). When the dissipation is simply of the form $\mathcal{R}_\eps ( u  , v ) = \tfrac12\|  v \|^2$ the abstract results of Sandier and Serfaty \cite{SandSerf04} provide sufficient conditions which guarantee the (weak) convergence of $u_\eps$ to the solution $u_0$ of the parabolic problem
$$
	\begin{cases} 
	\dot{u}_0 (x,t)  = - \nabla \E_0 (  u_0 (x,t) )  ,  \\
	u_0 ( 0, x) = u^0 (x) .
	\end{cases}
$$
For instance, under suitable regularity assumptions, convergence holds if the initial data are well-prepared, i.e., if $u^\eps \weakto u^0$ and $\E_\eps (u^\eps) \to \E_0 ( u^0)$, and if 
$$  \| \nabla \E_0 \| ( u )  \le \Gamma\text{-}\liminf_{\eps \to 0} \| \nabla \E_\eps \| ( u) .   $$
In the rate-independent setting, similar conditions are sufficient also for the convergence of balanced (or vanishing) viscosity evolutions, see e.g.~\cite{Negri_COCV14}. The conditions in \cite{SandSerf04} are actually far more general, and include the case of $\eps$-depending norms:
\begin{equation}  \label{e.SSlsc}    
		\| \nabla \E_0 \|_0 ( u )  \le \Gamma\text{-}\liminf_{\eps \to 0} \| \nabla \E_\eps \|_\eps ( u)  , 
\end{equation}
which will be relevant also in our context. Finally, consider a family of minimizing movement (or implicit Euler) schemes of the form 
$$
         u_\eps ( (k+1) \tau_\eps ) \in \mathrm{argmin} \left\{  \E_\eps (u) + \tfrac{1}{2}  \tau_\eps^{-1} \| u - u_\eps ( k \tau_\eps) \|^2   \right\},
$$
where $\tau_\eps \to 0$ as $\eps \to 0$. When the energy $\E_\eps$ is independent of $\eps$ it is well known that (the interpolation of) the above time-discrete evolutions converge, as $\eps \to 0$, to the gradient flow of $\E_\eps$. However, by virtue of the examples provided in Braides \cite{Braides_LNM13}, in general the limit evolution will depend on the relationship between $\eps$ (which parametrizes the family of energies) and the time step $\tau_\eps$.

\medskip
From the mechanical point of view our work is motivated by experimental and numerical results on the effective toughness of brittle composites (see for instance 
\cite{XiaPonsonRavichandranBhattacharya_PRL12, HossainHsuehBourdinBhattacharya_JMPS14} and the references therein). 
In general, when dealing with a periodically arranged composite, with a very fine microstructure, it is natural to introduce also an homogenized material reflecting ``in average'' the mechanical behavior of the composite. For instance, the equilibrium configuration (in terms of displacement, stress and energy) of a linear elastic composite can be approximated by the equilibrium configuration of a homogeneous linear elastic material, obtained with ``suitable averages'' of the stiffness tensors of the underlying constituents. 

For quasi-static fracture propagation, the mechanical behavior is instead described by the crack set, at time $t$, together with the equilibrium configuration of the body, given the crack at time $t$; further, the propagation is governed  by a rate-independent evolution law for the crack (Griffith's criterion)  together with momentum balance (an elliptic  {\sc pde}). For the latter (static) homogenization is suitable, and several mathematical approaches are available in the literature, see e.g.~\cite{MuratTartar_97, BraidesDefranceschi_98, MR503330} and the references therein. 
For the former, Griffith's criterion, it is needed instead some kind of ``evolutionary homogenization'' or, better, an ``evolutionary convergence'' for rate-independent systems (specified hereafter). In practice, see for instance \cite{XiaPonsonRavichandranBhattacharya_PRL12, HossainHsuehBourdinBhattacharya_JMPS14}, the goal would be to find an effective ``averaged'' toughness, in such a way that the quasi-static evolution of a crack in the brittle composite is approximated by the quasi-static evolution of a crack in a brittle homogeneous material characterized by: the (static) homogenization of the stiffness tensors together with the effective toughness. In \S\ref{explicit} we provide an example of a laminate with explicit computations of effective toughness and homogenized stiffness. This example highlights some peculiar features: the effective toughness depends on the volume fraction of the layers, on their toughness, and also on the coefficients of their stiffness matrices, further, it can be larger than the toughness of the single brittle layers, and larger than the weak$^*$ limit of the periodic toughness. These properties are qualitatively consistent with the numerical and experimental results of Hossain et al. \cite{HossainHsuehBourdinBhattacharya_JMPS14}, where toughening (i.e., the macroscopic increase in toughness) occurs solely as a result of the  microscopic elastic heterogeneities. In mechanics, most of the interest is indeed driven by the toughening and the strengthening effects of composites \cite{NM_14}. In general, toughening is also due to several microscopic ``geometrical features'' of the crack, such as micro-cracking, branching, debonding and tortuosity; a mathematical result on toughening by the homogenization of highly oscillating crack paths is given in \cite{Barchiesi_ARMA18}, where Barchiesi employs a modified $\Gamma$-convergence framework in order to take into account the irreversibility constraint.

\medskip
Let us turn to our results. We consider a brittle crack which propagates horizontally in a linear elastic laminate composed by $n$ periodic layers; each layer is made of a layer of material A and B with volume ratio $\lambda$ and $1-\lambda$, respectively. We distinguish between vertical \cite{HossainHsuehBourdinBhattacharya_JMPS14} and horizontal layers: in the former setting the crack crosses the layers, while in the latter it lays in the interface between two layers. We do not consider ``surfing boundary condition'' \cite{HossainHsuehBourdinBhattacharya_JMPS14} but a Dirichlet boundary condition $u = f(t) g$ for the displacement $u$ on a portion $\partial_D \Omega$ of the boundary. 

We begin with the more delicate case, that of vertical layers. We denote the toughness by $\Gc_n (l)$ and the energy release by $G_n (t , l )$, where $t $ is time and $l $ indicates the position of the crack tip. Note that $\Gc_n$ depends on $l$ and it is periodic. On the contrary $G_n (t, \cdot)$ is not periodic and, to our best knowledge, it is not known if it is defined at the interfaces between the two materials; for this reason, we actually employ the right lower semi-continuous envelope which is indeed the natural choice for Griffith's criterion (see Remark \ref{r.extend}). By linearity the dependence on time is much easier, since we can write $G_n (t, l) = f^2(t) \G_n (l)$, where $\G_n$ is computed with boundary condition $u = g$, i.e., for $f(t) =1$. Denoting by $\ell_n(t)$ the position of the tip at time $t$, the quasi-static propagation is governed by Griffith's criterion: 
\begin{itemize}\addtolength{\itemsep}{-0.2\baselineskip}
\item[i)] $G_n (t, \ell_n(t) ) \le \Gc_n (\ell_n(t))$, 
\item[ii)]  if $G_n ( t, \ell_n (t) ) < \Gc_n (\ell_n (t)) $ then $\dot{\ell}^+_n (t) =0$ (i.e., the right time derivative vanishes).
\end{itemize}
It is well known, see e.g.~\cite{NegriOrtner_M3AS08, Negri_ACV10}, that in general there are different solutions satisfying Griffith's criterion, and that these solutions may be discontinuous (in time). Discontinuities are basically due to the fact that Griffith's criterion is rate-independent (i.e., invariant under time reparametrization).  From the mathematical point of view, the behavior in the discontinuities characterizes the different notions of evolution, see \cite{Negri_ACV10}. Here, we denote by $J(\ell_n)$ the set times where $\ell_n$ jumps and following \cite{NegriOrtner_M3AS08} we require, besides i) and ii), that 
\begin{itemize}
\item[iii)]  if $t \in J ( \ell_n)$ then $G_n (t , l ) \ge \Gc_n (l)$ for every $l \in [\ell^-_n(t) , \ell^+_n(t) )$. 
\end{itemize}
This condition means that the transition between the equilibria $\ell^-_n(t)$ and $\ell^+_n(t)$ is unstable (or catastrophic). Notice that evolutions satisfying i), ii), and iii) can be obtained also by a vanishing viscosity approach \cite{KneesMielkZan08} taking the limit of continuous solutions of rate-dependent models. As observed in the numerical simulations \cite{HossainHsuehBourdinBhattacharya_JMPS14}, discontinuities often occur at those interfaces where the crack passes from the material with higher to that with lower toughness. Therefore, in our evolution $\ell_n$ we should expect many of these ``microscopic jumps'', in the order $n$ of the number of layers; these discontinuities, which ``disappear'' in the homogenized limit, will produce the toughening effect. 

At this point, let us describe the homogenized limit. Denote by $G_\hom (t, l)$ the energy release of the homogenized elastic energy. By Helly's theorem we know that $\ell_n \to \ell$ pointwise (upon extracting a subsequence). The goal is to define an effective toughness $\Gc_\eff$ in such a way that the limit $\ell$ satisfies Griffith's criterion: 
\begin{itemize}\addtolength{\itemsep}{-0.2\baselineskip}
\item[i)] $G_\hom (t, \ell (t) ) \le \Gc_\eff (\ell(t))$, 
\item[ii)]  if $G_\hom ( t, \ell (t) ) < \Gc_\eff (\ell (t)) $ then $\dot{\ell}^+ (t) =0$ (i.e., the right time derivative vanishes),
\item[iii)]  if $t \in J ( \ell )$ then $G_\hom (t , l ) \ge \Gc_\eff (l)$ for every $l \in [\ell^- (t) , \ell^+ (t) )$. 
\end{itemize}
Notice that in general $G_\hom (t, \cdot )$ is neither the pointwise nor the $\Gamma$-limit of $G_n (t,\cdot)$ (an explicit example is given in \S \ref{explicit}). However, a variational convergence in the spirit of \eqref{e.SSlsc} holds. Indeed, writing again $G_\hom (t , l ) = f^2(t) \G_\hom (l)$ the effective toughness is defined in such a way that 
\begin{equation} \label{e-i.eff-tough}
\displaystyle \frac{\G_\hom (l) }{\Gc_\eff (l) } = \Gamma\text{-}\liminf_{n \to \infty}  \frac{\G_n (l) }{ \Gc_n (l) } \,. 
\end{equation} 
First of all, note that in general $\Gc_\eff$ may depend on $l$ and, through the energy release, on several other parameters of the problem, in particular it does not depend only on the toughness of material A and B; in mechanical terms, it is an R-curve, see e.g.~\cite{And95}. 
At the current stage, it is hard to guess if $\Gc_\eff$ should always be a constant. In our examples it is so, on the contrary, in the experimental measures of \cite{NM_14} the effective toughness is surely not constant. 
What is clear is the fact that the effective toughness depends on the elastic contrast and it can be higher than the toughness of the underlying materials. Indeed, in the explicit calculations of \S \ref{explicit} we show that (under suitable assumptions)
$$
	\displaystyle \Gceff   =  \lambda  \max \left\{  \GcA  \, , \,  \GcB \, \frac{\mu_{\B,\1}}{\mu_{\A,\1}} \right\}   
		+ ( 1 - \lambda) \max \left\{  \GcA \frac{\mu_{\A,\1}}{\mu_{\B,\1}}  \, , \,  \GcB \right\} , 
$$
where $\Gc_\A$ and $\Gc_\B$ are the toughness while $\mu_{\A,\1}$ and $\mu_{\B,\1}$ are the components of the stiffness matrix, in the anti-plane setting. Moreover, independently of this example, if $\G_n$ converge uniformly to $\G_\hom$ then $\Gceff = \max\{ \GcA, \GcB \}$. 

Finally, let us compare \eqref{e-i.eff-tough} with \eqref{e.SSlsc}. In the context of fracture, the energy release $\G_n$ plays the role of the slope $\| \nabla \E_n \|$; 
the ratio $\G_n /  \Gc_n $ is instead some sort of energy release with respect to the dissipation metric induced by the toughness $\Gc_n$, and plays the role of $\| \nabla \E_n \|_n$; in these terms, the effective toughness plays the role of the norm $\| \cdot \|_0$ in \eqref{e.SSlsc} or, more precisely, of the smallest norm which makes \eqref{e.SSlsc} true, and for this reason in the definition of $\Gc_\eff$ an equality is needed.

In the case of horizontal layers with an interface crack the picture is much simpler. Indeed,  the toughness $\Gc$ is constant and independent of $n$, moreover the energy release $\G_n$ is well defined and converge locally uniformly to $\G_\hom$; this is enough to prove that the evolutions $\ell_n$ converge to an evolution $\ell$ which satisfies Griffith's criterion with energy release $G_\hom$ and toughness $\Gc_\eff= \Gc$.

\tableofcontents

\section{Vertical layers with horizontal crack} \label{s.e}

Let $\Omega = (0,L) \times (-H,H)$ be the (uncracked) reference configuration. For $ l \in (0,L]$ let $K_l = [0, l] \times \{0\}$ be the crack; in the presence of a crack $K_l$ the reference configuration is then $\Omega \setminus K_l$. 
We denote by $\partial_D \Omega$, independently of $l$, the union of the sets $\{ 0 \} \times ( -H,0) \cup (0,H)$ and $\{ L \} \times (-H,H)$.
\separe
Let $g \in H^{1/2} ( \partial_D \Omega)$. 
To fix the ideas, 
$$   g (x) = \begin{cases} -1  & x \in \{ 0 \} \times (-H,0),  \\ 1  & x \in \{ 0 \} \times (0,H) ,  \\ 0  & x \in \{ L \} \times ,  \end{cases} $$
is an admissible function.  By abuse of notation we still denote by $g$ a lifting of the boundary datum. 
Given $l \in (0,L]$ the spaces of admissible displacements and admissible variations are respectively
$$
\U_\l = \{ u \in H^1(\Omega \setminus K_l) : u = g \text{ in $\partial_D \Omega$} \} , 
\qquad
\V_\l = \{ v \in H^1(\Omega \setminus K_l) : v = 0 \text{ in $\partial_D \Omega$} \} .
$$
Note that $\U_l \neq \emptyset$ since $g \in H^{1/2} (  \{ L \} \times  (-H,H) )$.  Clearly if $l_1 < l_2$ then $\U_{\l_1} \subset \U_{\l_2}$ and thus  $\U_{\l} \subset \U_L \subset H^1(\Omega \setminus K_L)$ for every $l \in (0,L]$. All the spaces $\V_\l$ and the sets $\U_\l$ are endowed with the $H^1 ( \Omega \setminus K_L)$-norm. 
Splitting $\Omega \setminus K_L$ into the disjoint union of  $\Omega^+ = (0,L) \times (0,H)$ and $\Omega^- = (0,L) \times (-H,0)$, a uniform Poincar\'e inequality holds:~there exists $C_P>0$ such that for every $ l \in (0,L]$ and $v \in \V_l$ 
$$ 
      \int_{ \Omega \setminus K_l} | v |^2 \, dx \le C_P \int_{\Omega\setminus K_l} | \nabla v |^2 \, dx .
$$ 


The conditions given above are independent of the periodic structure and hold throughout the paper, for both abstract results and explicit computations. Actually, some of the results may hold under more general conditions, depending on the specific setting. For instance, we consider the antiplane setting since it allows more easily for explicit computations, however, the abstract results hold also in the plane-strain setting, discussed in \S \ref{s.pss}. 


\separe

\subsection{Energy and energy release} \label{s.2-1}

In this section we assume that $\Omega$ is a composite, made of $n \in \mathbb{N} \setminus \{ 0\}$ periodic vertical layers of thickness $l_n = L/n$; each layer itself is composed of a vertical layer of material $A$ with thickness $\lambda l_n$, for $\lambda \in (0,1)$, and a vertical layer of material $B$ with thickness $(1-\lambda) l_n$. 
More precisely, let $l_{n,k} = k\, l_n$ for $k=0,...,n$ and $l_{n,k+\lambda} =  \l_{n,k} + \lambda l_n = (k + \lambda) l_n$ for $k=0,...,n-1$; the layers of materials $A$ are of the form $(l_{n,k}, l_{n,k+\lambda}) \times (-H,H)$ while those of material $B$ are of the form $( l_{n,k+\lambda} , l_{n,k+1}) \times (-H,H)$.
For later convenience, we also introduce the notation $\Lambda_{n} = \{ l \in [0,L] : l = l_{n,k} \text{ or } l = l_{n,k+\lambda} \}$ so that the interfaces between material $A$ and $B$ (in the sound material) are of the form $\{ l \} \times (-H, H)$ for $l \in \Lambda_n \setminus \{ 0, L \}$. We will also employ the notation $\Omega_{n,k} = ( l_{n,k} , l_{n,k+1}) \times (-H,H) $ and 
$$   	\Omega_{n,\hspace{1pt}\A} = \bigcup_{k=0}^{n-1} \ ( l_{n,k} , l_{n,k+\lambda}) \times (-H,H) ,  \qquad 	\Omega_{n,\hspace{1pt}\B} = \bigcup_{k=0}^{n-1} \ ( l_{n,k +\lambda} , l_{n,k+1}) \times (-H,H) .	$$
We define in a similar way the crack sets $K_{n,k} =  ( l_{n,k} , l_{n,k+1} ) \times \{ 0 \} $, and then  $K_{n,\hspace{1pt}\A}$ and $K_{n,\hspace{1pt}\B}$.

We assume that $A$ and $B$ are elastic brittle materials. We denote by 
$\CA$ the stiffness matrix of material $A$, of the form 
$$
	\CA = \left( \begin{matrix} \mu_{\A,\1} & 0 \\ 0 & \mu_{\A,\2}  \end{matrix} \right) , 
	\quad \text{$\mu_{\A,i} >0$ for $i=1,2$, } 
$$
and we denote by $\GcA>0$ the toughness of $A$; we employ a similar notation for material $B$. 
Let $\C_n : \Omega \setminus ( \Lambda_n \times (-H,H)) \to \{ \CA , \CB \}$ be the periodic stiffness matrix and let $\Gc_n : [0,L] \setminus \Lambda_n \to \{ \GcA , \GcB \}$ be the periodic toughness. 

\separe

The linear elastic energy $\W_{l,n}: \U_l \to \R$ is given by 
$$
     \W_{l,n} (u) = \tfrac12 \int_{\Omega \setminus K_l} \nabla u \C_n \nabla u^T \, dx .
$$ 
Denote by $u_{l,n} \in \U_\l$ the (unique) minimizer of the energy $\W_{l,n}$; clearly $u_{l,n}$ is also the (unique) solution of the variational problem  
$$   \int_{\Omega \setminus K_l} \nabla u \C_n \nabla v^T \, dx = 0   \quad \text{for every $v \in \V_l$. }  $$
Note that the above bi-linear form is coercive and continuous, uniformly with respect to $l \in (0,L]$ and $n \in \mathbb{N}$.
In the sequel we will often employ the condensed (or reduced) elastic energy $\E_n : (0,L] \to \R$ given by 
$$
     \E_n ( l ) = \W_{l,n} ( u_{l,n} ) = \min \{   \W_{l,n} ( u) : u \in \U_\l \} .
$$
The following lemma contains the fundamental properties of this energy (for $n$ fixed).
 
\begin{lemma} \label{l.En} The energy $\E_n$ is decreasing, continuous, and of class $C^1 ( (0,L) \setminus \Lambda_n )$. In particular, if $l \in (l_{n,k} , l_{n,k+\lambda})$ the derivative takes the form 
\begin{equation} \label{e.E'=DE}
      \E'_n ( l) = \int_{\Omega \setminus K_l} 
      \nabla u_{l,n} \boldsymbol{C}_n \boldsymbol{E} \nabla u_{l,n}^T \, \phi' dx 
	 \quad \text{where} \quad 
    \boldsymbol{E} =   \left( \begin{matrix} -1 &  0 \\ 0 & 1 \end{matrix} \right) 
\end{equation}
while $\phi \in W^{1,\infty} (0,L)$, $\mathrm{supp} ( \phi ) \subset ( l_{n,k} , l_{n,k+\lambda} ) $, $\phi (l) = 1$ and, by abuse of notation, $\phi' = \partial_{l} \phi (x_1) $. A representation like \eqref{e.E'=DE} holds also for $l \in (l_{n,k+\lambda}, l_{n,k+1})$. 
\end{lemma}

\begin{remark} Note that $\E_n$ is not periodic, unless it is constant. Note also that the support of the auxiliary function $\phi$ shrinks when $n \to \infty$. 
\end{remark}

\proof We provide a short proof, following \cite{N_AMO17}. For $l_1 < l_2$ we have $\U_{\l_1} \subset \U_{\l_2}$ and thus  $\E_n ( l_1) \ge \E_{n} ( l_2)$. 
If $l_m \to l$ then $\E_n (l_m) \to \E_n ( l)$ and $u_{l_m,n} \to u_{l,n}$ strongly in $H^1 (\Omega \setminus K_L)$.

\separe

Let us compute the derivative of the energy $\E_n (l)$ for $l \not\in \Lambda_n$.  
For $h \ll 1$ let $\Psi_h ( x) = (x_1 + h \phi  (x_1), x_2)$ be a diffeomorphism of $\Omega \setminus K_{l}$ onto $\Omega \setminus K_{l+h}$. Note that $\Psi_h (x) = x$ out of the layer $(l_{n,k} , l_{n,k+\lambda} ) \times  (-H,H)
$. Using the change of variable $w = u \circ \Psi_h$ the energy $\W_{l+h,n} : \U_{\,l+h} \to \R$ reads
$$
	\W_{l+h,n} ( u ) = \overline{\W}_{l,h,n} (w) = \tfrac12 \int_{\Omega \setminus K_l} \nabla w  \C_n ( \boldsymbol{I} + h \boldsymbol{M}_{\!h} \phi' ) \nabla w^T  \, dx 
	\quad \text{where} \quad 
	\boldsymbol{M}_{\!h} = \left(  \begin{array}{cc}  - ( 1 + h \phi')^{-1} & 0 \\ 0 & 1 \end{array} \right). 
$$
If $u_{l+h,n} \in \argmin \{ \W_{l+h,n} (u) : u \in \U_{\,l+h} \}$ then $w_{l,h,n} = (u_{l+h,n} \circ \Psi_h) \in \argmin \{ \overline{\W}_{l,h,n} (w) : w \in \U_{\,l} \}$ and $\E_n (l+h) = \overline{\W}_{l,h,n} ( w_{l,h,n})$. Hence, 
$$
	 \E'_n (l) = \lim_{h \to 0} \frac{\E_n (l+h) - \E_n (l)}{h} = \lim_{h \to 0} \frac{\overline{\W}_{l,h,n} (w_{l,h,n}) - \W_{l,n} (u_{l,n})}{h} .
$$
Writing the variational formulations for $w_{l,h,n}$ and $u_{l,n}$ we have 
$$
	\int_{\Omega \setminus K_l} \nabla w_{l,h,n} \C_n ( \boldsymbol{I} + h \boldsymbol{M}_{\!h} \phi' )  \nabla v^T  \, dx  = 
	\int_{\Omega \setminus K_l} \nabla u_{l,n}   \C_n  \nabla v^T  \, dx = 0 
	\quad \text{for every $v \in \V_{\l}$.} 
$$
Then, being $w_{l,h,n} - u_{l,n} \in \V_l$, we can re-write the energies as 
\begin{align*}
	\overline{\W}_{l,h,n} (w_{l,h,n}) & = \tfrac12 \int_{\Omega \setminus K_l} \nabla w_{l,h,n} \C_n ( \boldsymbol{I} + h \boldsymbol{M}_{\!h} \phi' )  \nabla u_{l,n}^T  \, dx , 
	\\
	 \W_{l,n} (u_{l,n}) & = \tfrac12 \int_{\Omega \setminus K_l} \nabla w_{l,h,n}   \C_n   \nabla u_{l,n}^T  \, dx .
\end{align*}
Hence, 
$$
	\frac{ \overline{\W}_{l,h,n} (w_{l,h,n}) - \W_{l,n} (u_{l,n})}{h}=  \tfrac12 \int_{\Omega \setminus K_l} \nabla w_{l,h,n}  \C_n \boldsymbol{M}_{\!h} \nabla u_{l,n}^T  \, \phi'  \, dx .
$$
As $h \to 0$, $\boldsymbol{M}_{\!h} \to 2 \boldsymbol{E}$ in $L^\infty(\Omega ; \R^{2 \times 2} )$ and  $w_{l,h,n} \weakto u_{l,n}$ in $H^1( \Omega \setminus K_l)$, therefore we obtain \eqref{e.E'=DE}. 

If $l_m \to l$ and  $l_{n,k} <  l  < l_{n,k+\lambda}$ we can choose $\phi$ to be independent of $m$ (sufficiently large). Moreover, $\nabla u_{l_m,n} \to \nabla u_{l,n}$ in $L^2(\Omega \setminus K_L , \R^2)$. Therefore, by the representation \eqref{e.E'=DE} we get $\E'_n (l_m) \to \E'_n (l)$. \qed 



By virtue of the above Proposition for every $l \not\in \Lambda_n$ we can define the energy release 
\begin{equation}   \label{e.enrel}  \mathcal{G}_n (l) = - \partial^+_l \E_n (l)  = - \partial_l \E_n (l)   .  \end{equation} 
In the sequel it will be convenient to extend the toughness and the energy release to $\Lambda_n$ (even if in general the energy is not differentiable in $\Lambda_n$). 
We define 
\begin{equation} \label{e.extGn}
       \G_n ( l ) = \liminf_{s \to l^+} \G_n ( s) \text{      if $l \neq L$ } \quad \text{ and } 
       \quad 
       \G_n ( L) = 0 .
\end{equation}
With this definition the function $\G_n : [0,L] \to [0,+\infty] $ is right lower semi-continuous. In Remark \ref{r.extend} we will see that for $l \neq L$ this extension is actually the only one compatible with Griffith's criterion; for $l=L$ it is instead suggested by technical convenience (strictly speaking, the right derivative of the reduced energies $\E_n$ do not make sense for $l=L$). Note that $\G_n$ may possibly take value $+\infty$ in $\Lambda_n$. 
In a similar manner, we define 
\begin{equation} \label{e.extGc}
  \Gc_n (l) = \lim_{s \to l^+}  \Gc_n (s)\text{      if $l \neq L$ } \quad \text{ and } 
       \quad 
       \Gc_n (L) = \GcA. \end{equation}
Clearly, $G_n^c : [0,L] \to \{ \GcA, \GcB \}$ is right continuous.




\subsection{Quasi-static evolution by Griffith's criterion} \label{s.2-2}


In this subsection we provide a notion of quasi-static propagation, following \cite{NegriOrtner_M3AS08, KneesMielkZan08, Negri_ACV10}. We consider a time depending boundary datum of the form $f (t) g $ where $f \in C^1([0,T])$ with $f(0)=0$ and $f' > 0$ (non-monotone boundary conditions are considered in \S \ref{s.non-mono}). Let $L_0 \in (0,L)$ be the initial crack. 
For $l \in [L_0,L]$ the space of admissible displacement is
$$ U_{t,l} = \{ u \in H^1(\Omega \setminus K_l) : u = f(t) g \text{ in $\partial_D \Omega$} \}  = f(t) \, \U_\l .  $$
We denote by $E_n$ and $G_n$ the corresponding condensed elastic energy and energy release. 
For $t=0$ we simply have $E_n (0, \cdot) = G_n (0, \cdot) = 0$ since $u=0$ on $\partial_D \Omega$. If $t>0$ by linearity we can write 
$$   E_n ( t, l ) =  f^2 (t) \,\E_n (l) ,   \qquad  G_n ( t, l) = f^2(t) \, \G_n (l) .  $$
The latter inequality makes sense in generalized form also in the case $\G_n (l) = +\infty$. In particular, $E_n ( \cdot ,l)$ is of class $C^1 ( [0,T])$ while, by Lemma \ref{l.En}, $E_n ( t, \cdot)$ is decreasing, continuous and of class $C^1( (0,L) \setminus \Lambda_n)$. The energy release rate $G_n (\cdot , l)$ is instead of class $C^1 ( [0,T] )$ if $\G_n (l ) < + \infty$ (e.g. if $l \not\in \Lambda_n$), otherwise  $G_n (\cdot , l) = + \infty$ for  $t>0$ (in this case we will say that $G_n ( \cdot , l)$ is continuous, with extended values). 
Finally, $G_n ( t , \cdot)$ is right lower semi-continuous in $[0,L)$. 
%
%

\separe

The following proposition provides our formulation of  Griffith's criterion for the quasi-static propagation of  a crack in a brittle material with periodic vertical layers.

\begin{proposition} \label{p.existsn} There exists a unique non-decreasing, right continuous function $\ell_n : [0,T) \to [L_0, L]$ which satisfies the initial condition $\ell_n(0) = L_0$ and Griffith's criterion in the following form: 
\begin{itemize}
\item[i)] $G_n (t, \ell_n(t) ) \le \Gc_n (\ell_n(t))$ for every time $t \in [0,T)$;
\item[ii)]  if $G_n ( t, \ell_n (t) ) < \Gc_n (\ell_n (t)) $ then $\ell_n$ is right differentiable in $t$ and $\dot{\ell}^+_n (t) =0$; 
\item[iii)]  if $t \in J ( \ell_n)$ then $G_n (t , l ) \ge \Gc_n (l)$ for every $l \in [\ell^-_n(t) , \ell_n(t) )$;
\end{itemize}
where $J (\ell_n) $ denotes the set of discontinuity points of $\ell_n$.
\end{proposition} 

Before proving Proposition \ref{p.existsn} let us make a few comments. Uniqueness is true thanks to the monotonicity of the boundary condition together with the right continuity of the evolution.  
The discontinuity points $t \in J(\ell_n)$ correspond to instantaneous non-equilibrium transitions of the system between the equilibria $\ell_n^- (t)$ and $\ell_n(t) = \ell_n^+(t)$. As shown in the numerical simulations of \cite{HossainHsuehBourdinBhattacharya_JMPS14} discontinuities are often tailored to the layers: roughly speaking, the crack jumps through the layer with smaller toughness, reaching the closest interface, where it ``re-nucleates'' before advancing in the layer with higher toughness. In this case the non-equilibrium condition $G_n (t, l) \ge \Gc_n (l)$ may not hold for $l = \ell_n^+(t) \in \Lambda_n$. For instance, it may happen that 
$G_n (t, l) > \Gc_\A $ for $l \in (l_{n,k} , l_{n,k+\lambda} ) $ while $G_n (t, l_{n,k+\lambda} ) < \Gc_\B $. 
For this reason condition ii) holds with the right derivative $\dot\ell_n^+(t)$; however, out of interfaces the full time derivative $\dot\ell_n (t)$ can be used (see Corollary \ref{c.dotell}).

Strictly speaking, right continuity is not necessary for Griffith's criterion, however it is also not a restrictive since $\ell_n$ turns out to be the (only) right continuous representative of any evolution $\lambda_n$ which satisfies Griffith's criterion. More precisely, we have the following result, whose proof is postponed after Corollary  \ref{c.dotell}. 

\begin{corollary} \label{c.unica} If $\lambda_n : [0,T] \to [L_0,L]$ is non-decreasing, satisfies $\lambda_n (0) = L_0$ and Griffith's criterion, i.e., 
\begin{itemize}
\item[i)] $G_n (t,  \lambda_n(t) ) \le \Gc_n (  \lambda_n(t))$ for every time $t \in [0,T]$;
\item[ii)]  if $G_n ( t,  \lambda_n (t) ) < \Gc_n (  \lambda_n (t)) $ and $t<T$ then $\lambda_n$ is right differentiable in $t$ and $\dot\lambda^+_n (t) =0$; 
\item[iii)] if $t \in J (  \lambda_n)$ then $G_n (t , l ) \ge \Gc_n (l)$ for every $l \in [\lambda^-_n(t) , \lambda^+_n(t) )$;
\end{itemize}
then the right continuous representative of $\lambda_n$ coincides with $\ell_n$ in $[0,T)$. 
\end{corollary}

\medskip
\noindent {\bf Proof of Proposition \ref{p.existsn}.} 
The proof is based on the following explicit representation of the evolution: 
\begin{equation} \label{e.ell+}
 \ell_n (t) = \inf \{ l \in [L_0,L] : G_n (t , l) < \Gc_n (l) \} . 
\end{equation}
Note that the set $\{ l \in [L_0,L] : G_n (t , l) < \Gc_n (l) \}$ is not empty because, by definition, $G_n (t, L) =0$. Thus $\ell_n$ is well defined and takes values in $[L_0,L]$. We will check that the function $\ell_n$ given by \eqref{e.ell+} is indeed the unique solution.

\separe

Clearly $\ell_n$ satisfies the initial condition because $G_n ( 0 , L_0) = 0$.

\separe

Let us check monotonicity. Since the function $f$ is increasing and $G_n$ is non-negative, for every $0 \le t_1 < t_2 < T$ it holds $G_n ( t_\1, l )  = f^2(t_\1) \G_n ( l) \le f^2(t_\2) \G_n (l) = G_n ( t_\2, l )$, hence 
$$
	\{ l \in [L_0,L] :   G_n ( t_\1, l) < \Gc_n (l) \}  \supset 
	\{ l \in [L_0,L] :   G_n ( t_\2, l) < \Gc_n (l) \} \,. 
$$
Taking the infimum, yields $\ell_n (t_\1) \le \ell_n (t_\2)$.  

\separe

Next, we show that $\ell_n$ is right continuous. Since in general $G_n ( \cdot , l)$ is continuous (with extended values) only in $(0,T)$, we distinguish between $t=0$ and $t>0$. 
In the former case, let $L_0 < l$ with $l \not\in \Lambda_n$. Then as $\tau \to 0$ we have $G_n ( \tau , l) \to 0 < \Gc_n (l)$. Hence, \eqref{e.ell+} implies $\ell_n (\tau) \le l$ for $\tau > 0$ sufficiently small. By monotonicity $\ell_n^+(0)  \le l$ and thus $\ell_n^+ (0) \le L_0$ by the arbitrariness of $l$.  Let $t > 0$.  By monotonicity $\ell^+_n (t) \le \lim_{\tau \to t^+} \ell_n (\tau)$. Let us prove the opposite inequality. If $\lim_{\tau \to t^+} \ell_n (\tau) = L_0$ there is nothing to prove. 
Otherwise, let $L_0 \le l^* < \lim_{\tau \to t^+} \ell_n(\tau) $ with $l^* \not\in \Lambda_n$ (remember that $\Lambda_n$ is a discrete set). By monotonicity $l^* < \ell_n (\tau) $ for every $\tau > t$ and thus applying the definition \eqref{e.ell+} to $\ell_n(\tau)$ we know that $G_n (\tau, l^*) \ge  \Gc_n (l^*)$ for every  $\tau > t$. Passing to the limit as $\tau \to t^+$ by the continuity of $G_n ( \cdot, l^*)$ it follows that $G_n (t, l^*) \ge  \Gc_n(l^*)$ for every $L_0 \le l^* < \lim_{\tau \to t^+} \ell(\tau)$ with $l^* \not\in \Lambda_n$. Hence $\ell_n(t) = \inf \{ l \in [L_0,L] : G_n ( t, l) <  \Gc_n(l) \}  \ge l^*$. Taking the supremum in the right hand side for $L_0 \le l^* < \lim_{\tau \to t^+} \ell(\tau)$ with $l^* \not\in \Lambda_n$, we get $\ell_n (t) \ge \lim_{\tau \to t^+} \ell(\tau) $.

\separe

Let us check i). If $\ell_n(t) < L$ then by \eqref{e.ell+} there exists a sequence $l_m \searrow \ell_n(t)$ such that $G_n (t, l_m) <  \Gc_n (l_m)$ and $l_m \not\in \Lambda_n$. By the right lower semi-continuity of $G_n (t, \cdot)$ and the right continuity of $ \Gc_n$ it follows that 
\begin{equation} \label{e.KT1}
G_n ( t , \ell_n (t) ) \le \liminf _{m \to \infty} G_n (t, l_m ) \le \lim_{m \to \infty}  \Gc_n (l_m) =  \Gc_n ( \ell_n (t)).
\end{equation}
 If $\ell_n(t) =L$ there is nothing to prove since, by definition, $ 0 = G_n (t , L ) < \Gc_n (L)$.
 
\separe

Next we check ii). Again, we consider separately $t=0$ and $t>0$. If $t=0$ we have $G_n (0, L_0) < \Gc_n (L_0)$. If $\ell_n (\tau) = L_0$ for some $\tau >0$ then $\ell_n$ is constant up to $\tau$ and thus $\dot\ell_n^+ (0)=0$. On the contrary, assume that $\ell_n (\tau) > L_0$ for every $\tau >0$. By \eqref{e.ell+} we know that $G_n (\tau, l) \ge \Gc_n (l)$ for every $\tau>0$ and every $L_0 \le l < \ell_n(\tau)$. Hence $f^2(\tau) \G_n ( l ) \ge \min \{ \GcA  , \GcB \}$ and thus 
$$
f^2 (\tau) \big(  \E_n (L_0) - \E_n ( \ell_n(\tau) \big) =  f^2(\tau)  \int_{L_0}^{\ell_n(\tau)} \G_n ( l)  \, dl \ge ( \ell_n ( \tau) - L_0 ) \min \{ \GcA  , \GcB \} .
$$
As a consequence
$$
	\lim_{\tau \to 0^+} \frac{\ell_n ( \tau) -  \ell_n(0)}{\tau } \le  \lim_{\tau \to 0^+} \frac{ f^2(\tau) - f^2(0)  } {\tau} \frac{ \E_n (L_0) - \E_n ( \ell_n(\tau))}{ \min \{ \GcA  , \GcB \} } =  0 ,
$$
where we used the continuity of the energy $\E_n$ and the fact that $f \in C^1([0,T])$ with $f(0)=0$. Now, consider $t>0$. If $G_n (t, \ell_n(t)) < G_n^c( \ell_n(t))$ then by continuity of $G_n ( \cdot, l)$ we get $G_n ( \tau , \ell_n(t)) <  \Gc_n (\ell_n(t))$ for $\tau$ in a (sufficiently small) right neighborhood of $t$. By \eqref{e.ell+} it follows that $\ell_n (\tau) \le \ell_n (t)$; by monotonicity $\ell_n$ is constant in a right neighborhood of $t$.

\separe

Let us consider iii). Let $t \in J (\ell_n) $. By right continuity $\ell_n(t) > L_0$. By definition $\ell_n (t) = \inf \{ l \in [L_0, L] : G_n (t , l) <  \Gc_n(l) \}$, hence $ G_n (t, l ) \ge  \Gc_n (l)$ for  every $L_0 \le l < \ell_n (t)$ and in particular for $l \in [ \ell_n^- (t) , \ell_n (t) )$.

\separe

To conclude, let us prove uniqueness. Assume by contradiction that there exists another right continuous evolution $\ell_* \neq \ell_n$ which satisfies Griffith's criterion. First, we claim that there exists a time $t \not\in J ( \ell_*) $ such that $\ell_* (t) \neq \ell_n (t)$; indeed, if $\{ t \in [0,T] : \ell_* (t) \neq \ell_n (t) \} \subset  J ( \ell_*)$ then by right continuity of $\ell_*$ and $\ell_n$ we would get $\ell_*  = \ell_n$ everywhere in $[0,T)$. So, let $t>0$ be a continuity point of $\ell_*$ with $\ell_* (t) \neq \ell_n (t)$.

If $\ell_* (t) > \ell_n (t)$ then, by the definition \eqref{e.ell+} of $\ell_n$, there exists $l_* \in ( \ell_n (t) , \ell_*(t))$ with $l_* \not\in \Lambda_n$ such that $G_n ( t , l_* ) <  \Gc_n (l_*)$. Then, by continuity of the energy release, there exists a neighborhood $[l',l'']$ of $l_*$, with $L_0 \le l'$ and $l'' < \ell_*(t)$, such that $G_n ( t , l ) <  \Gc_n (l)$ for every $l \in [l',l'']$. Let us see that this interval is a ``barrier'' for the evolution $\ell_*$. First, note that by monotonicity in time, $G_n ( \tau , l ) \le G_n (t , l ) <  \Gc_n (l)$ for every $\tau \le t$ and every $l \in [l',l'']$.  Moreover, $\ell_*$ takes all the values in the interval $[l',l'']$, indeed if $ l \in [ \ell_*^- ( \tau ) , \ell_*( \tau))$ for some $\tau \le t$ and some $l \in [l', l'']$ then, by condition iii) of Griffith's criterion, $G_n (\tau, l) \ge  \Gc_n (l)$, which is a contradiction with $G_n ( \tau , l ) <  \Gc_n (l)$. Therefore we can define $t_\1 = \inf \{ \tau : l'  \le  \ell_* ( \tau) \}$ and $t_\2 = \sup \{ \tau : \ell_* (\tau) \le l'' \}$.  It follows that $ G_n ( \tau , \ell_* (\tau) )  < \Gc_n (l) $ for every $ \tau \in [ t_\1 , t_\2 ) $ and thus, by condition ii) of Griffith's criterion, $\ell_*$ is constantly equal to $l'$ in $[t_\1, t_\2)$. If $t_\2 = t$ then (being $t$ is a continuity point) we get $\ell_* (t) = l'  < \ell_* (t)$ . If $t_\2< t$ then $\ell_*^- (t_\2) = l'$ while $\ell_*^+ (t_\2) \ge l''$ and thus $t_\2 \in J ( \ell_*)$, this is again a contradiction since Griffith's criterion implies that $G_n (t_\2 , l) \ge  \Gc_n (l)$ for every $l \in [l' , l'')$. 

Let $\ell_* (t) < \ell_n (t)$. By the definition of $\ell_n$ we know that $G_n (t , l ) \ge  \Gc_n (l) > 0$ for every $ l \in [ \ell_*(t) , \ell_n (t) )$. 
By monotonicity $G_n (\tau , l) >  \Gc_n(l)$ for every $\tau > t$ and every $ l \in [ \ell_*(t) , \ell_n (t) )$. If $\tau \to t^+$ 
then  $\ell_* (\tau) \to \ell_*(t)$ and thus $ \ell_* (\tau) \in [ \ell_*(t) , \ell_n (t) )$ for every time $\tau$ sufficiently close to $t$; then we have $G_n (\tau , \ell_*(\tau) ) >  \Gc_n (\ell_*(\tau))$, which is a contradiction with Griffith's criterion. \qed


\separe

The arguments in the proof of Proposition \ref{p.existsn} allow to prove also the following result, which slightly improves condition ii) for $t>0$. 

\begin{corollary} \label{c.elldot} Let $t >0$ such that $G_n (t, \ell_n (t)) < \Gc_n (\ell_n(t))$ then $\ell_n$ is constant in a (sufficiently small) right neighborhood of $t$. 
\end{corollary}

In order to prove the energy identity we will need also the following corollary, which refines condition ii) for  $\ell_n (t) \not\in \Lambda_n$. 

\begin{corollary} \label{c.dotell} If  $\ell_n(t) \not\in \Lambda_n$ and $G_n ( t , \ell_n (t)) < \Gc_n ( \ell_n (t))$ then $\ell_n$ is constant in a (sufficiently small) neighborhood of $t$. 
\end{corollary}

\proof Let $\ell_n(t) \in (l_{n,k} , l_{n,k+\lambda})$. The energy $E_n$ is of class $C^1$ in $[0,T] \times (l_{n,k} , l_{n,k+\lambda})$ and thus the energy release $G_n$ is continuous. Moreover, $\Gc_n$ is constant in $(l_{n,k} , l_{n,k+1})$.  Hence, there exists $\delta >0$ (sufficiently small) such that $G_n ( \tau , l) < \Gc_n (l)$ whenever $| \tau - t | < \delta$ and $| l - \ell_n(t) | < \delta$. 
By condition iii) in Proposition \ref{p.existsn} it follows that $t \not \in J(\ell_n)$. By continuity of $\ell_n$ let us choose $\tau' < t$ such that $| \tau' - t | < \delta$ and $| \ell_n(\tau')  - \ell_n(t) | < \delta$. Then, $G_n ( \tau , \ell_n (\tau') ) < \Gc_n (\ell_n (\tau'))$ for every $\tau \in [ \tau' , t]$. 
By Proposition \ref{p.existsn} it follows that $\ell_n( \tau ) = \ell_n(t)$ is the unique solution for $\tau \in [ t , t + \delta]$.   The same argument applies if $\ell_n(t) \in (l_{n,k+\lambda} , l_{n,k+1})$. 
\qed 

\noindent {\bf Proof of Corollary \ref{c.unica}.} By uniqueness it is enough to check that the right continuous representative of $\lambda_n$, here denoted by $\hat\lambda_n$, satisfies Griffith's criterion, in the sense of Proposition \ref{p.existsn}. 

Arguing by contradiction it is easy to check that $\lambda_n$ is right continuous in $0$; indeed, if $\lim_{\tau \to 0^+} \lambda_n (\tau) > L_0$ then $0 \in J ( \lambda_n)$ and thus $G_n ( 0 , l) \ge \Gc_n (l) > 0$, however $G_n ( 0 , l ) = 0$. Hence $ \hat\lambda_n (0) = L_0$.

Let us check condition i). 
For $t=0$ there is nothing to prove. For $t>0$ we have $G_n ( \tau , \lambda_n (\tau) ) \le \Gc_n ( \lambda_n (\tau) )$ for every $\tau >t$ and thus the regularity of energy release and toughness imply that  
$$
	G_n (t ,  \hat\lambda_n (t) ) \le \liminf_{\tau \to t^+} G_n (\tau , \lambda_n (\tau)) \le \lim_{\tau  \to t^+} \Gc_n ( \lambda_n (\tau)) = \Gc_n (  \hat\lambda_n (t))  .
$$

Let us consider iii). Clearly, $J ( \lambda_n) = J ( \hat\lambda_n)$ and $ \hat\lambda^\pm_n (t) = \lambda^\pm_n (t)$. Hence the condition $G_n (t , l ) \ge \Gc_n (l)$ for every $l \in [ \hat\lambda^-_n(t) ,  \hat\lambda^+_n(t) )$ holds. 

Finally, let us check ii). If $t \not\in J ( \lambda_n) = J ( \hat\lambda_n)$ then $ \lambda_n (t) =  \hat\lambda_n (t)$ and $\dot{\lambda}_n^+ (t) = 0$ coincides with the right derivative of $\hat{\lambda}_n(t)$, hence there is nothing to check.  If $ t \in J ( \lambda_n) = J ( \hat\lambda_n)$ we argue in the following way. We have $G_n (t ,  \hat\lambda_n (t)) < \Gc_n (  \hat\lambda_n (t))$, then for $\eps > 0$ sufficiently small there exists a sequence $l_k \searrow \hat\lambda_n (t)$ such that $l_k \not\in \Lambda_n$, $\Gc_n (l_k) = \Gc_n (  \hat\lambda_n (t)) $ is constant  and 
$$
    G_n (t, l_k) \le \Gc_n (l_k) - \eps \quad \text{for every $k \in \mathbb{N}$.}
$$
Writing $G_n ( t , l_k) = f^2 (t) \G_n (l_k)$, by time continuity for every $0 < \eps' < \eps$ there exists $\delta>0$ such that 
$$   G_n (\tau, l_k) < \Gc_n (l_k) - \eps' \quad \text{for every $k \in \mathbb{N}$ and every $\tau \in ( t , t + \delta)$.} $$ 
Since $l_k \not\in \Lambda_n$ the energy release $\G_n$ is continuous, and thus for every $k \in \mathbb{N}$ there exists $l_k^\flat < l_k < l_k^\sharp < l_{k+1}$ such that 
$$
      G_n (\tau, l) < \Gc_n (l) \quad \text{for every $l \in (l_k^\flat, l_k^\sharp)$ and every $\tau \in ( t , t + \delta)$.} 
      $$ 
As in proof of Proposition \ref{p.existsn} the evolution $\lambda_n$ cannot take value in the  intervals $(l^\flat_k, l_k^\sharp)$ for $\tau \in ( t , t + \delta)$. Indeed, there are no jumps, moreover, if $\lambda_n (\tau') \in (l^\flat_k, l_k^\sharp)$ then 
condition ii) implies that the right derivative vanishes and thus $\lambda_n$ is constant in $(\tau', t + \delta)$. 
Hence by monotonicity $\lambda_n (\tau) \le l_k^\sharp$ for every $k$ and every $\tau \in ( t , t + \delta)$; in conclusion $\hat\lambda_n (\tau)  \le \inf l_k^\sharp = \hat\lambda_n (t) $ for every $\tau \in ( t , t + \delta)$. \qed 

\separe

\begin{remark} \normalfont \label{r.extend} 
For $l \in \Lambda_n$ the extension \eqref{e.extGc} of the toughness seems the most natural since 
$$	\Gc_n (l) = \lim_{s \to l^+} \Gc_n (s) . $$ 
The extension \eqref{e.extGn}, i.e., 
$$
	\G_n (l) =  \liminf_{ s \to l^+}  \G_n (s) 
$$
turns out to be the only one compatible with Griffith's criterion. In principle, since the existence of $\partial^+_l \E_n (l)$ is not known when $l \in \Lambda_n$ any value $\G_n (l)$ with 
\begin{equation} \label{e.bau}
	\liminf_{s \to l^+}  \frac{\E_n (s) - \E_n (l)}{s - l}  \le - \G_n (l) \le \limsup_{s \to l^+}  \frac{\E_n (s) - \E_n (l)}{s - l} 
\end{equation}
could be a good candidate for the extension. To fix the ideas assume that $L_0 \in \Lambda_n$ and let $\tilde\G_n(L_0)$ be another extension of $\G_n$ in $L_0$. First of all, \eqref{e.bau} implies that $\tilde{\G}_n(L_0) \ge \G_n ( L_0)$, indeed, by the mean value theorem we have 
$$
	\G_n (L_0) =  \liminf_{ s \to L_0^+}  \G_n (s)  \le \liminf_{s \to L_0^+}  - \frac{\E_n (s) - \E_n (L_0)}{s - L_0} = - \limsup_{s \to L_0^+}  \frac{\E_n (s) - \E_n (L_0)}{s - L_0} \le \tilde{\G}_n(L_0) .
$$
Let us consider 
$$	\tilde{\G}_n (L_0) > \G_n (L_0) =  \liminf_{ s \to L_0^+}  \G_n (s) .  $$ 
Clearly, this condition makes sense only when $\liminf_{ s \to L_0^+}  \G_n (s) < +\infty$.  Hence $ 0 <  \tilde{\G}_n (L_0) \le +\infty$.
Let 
$$    \tilde{G}_n ( 0 , L_0) = 0  ,\quad \tilde{G}_n ( t , L_0 ) = f^2(t) \, \tilde\G_n(L_0) \ \text{for $t>0$.}  $$ 
Hence $\tilde{G}_n (t, L_0) > G_n (t, L_0)$ for $t > 0$. Assume that there exists a quasi-static evolution $\tilde\ell_n$ which satisfies Griffith's criterion with $\tilde{G}_n$ and $\Gc_n$. We denote by $t_*$ and $t^*$ the points $t_* = \max \, \{ t \ge 0 : \tilde{G}_n ( t , L_0) \le \Gc_n (L_0) \}$ and $t^* = \max \, \{ t \ge 0 : G_n ( t , L_0) \le \Gc_n(L_0) \}$. We have $t_*  <  t^*$. Fix $t' \in ( t_* , t^*)$. By the definition of $\G_n (L_0)$ there exists $l_m \searrow L_0$ such that $l_m \not\in \Lambda_n$ and $\G_n (l_m) \to \G_n ( L_0) < \tilde{\G}_n (L_0)$.  Upon extracting a subsequence, not relabeled, we can assume that there exists $t_m < t'$ such that either $\tilde\ell_n (t_m) = l_m$ for every $m \in \mathbb{N}$ or $l_m \in [ \tilde\ell^-_n (t_m ) , \ell_n (t_m ) )$ for every $m \in \mathbb{N}$.

In the former case, for every $t \in (t_m , t')$ we have 
$$
	G_n ( t , l_m) \le G_n ( t' , l_m) \to G_n ( t', L_0) <  \Gc_n ( L_0)  = \Gc_n (l_m) .
$$
Hence for $m \gg 1$ we have $G_n ( t , l_m) < \Gc ( l_m)$ for every $t \in (t_m , t')$. As a consequence $\tilde\ell_n (t) = l_m$ in $(t_m , t')$. Repeating the argument for $l_{m+1}$ leads to $\tilde\ell_n (t) = l_{m+1}$ in $(t_{m+1} , t')$ which is a contradiction since $t_{m+1} < t_m$ and $l_{m+1} \neq l_m$.

In the other case, i.e., when $l_m \in [ \tilde\ell^-_n (t_m ) , \ell_n (t_m ) )$ for every $m \in \mathbb{N}$, we argue as follows.
As above, 
$$
	G_n (t_m , l_m) \le G_n (t' , l_m) \to G_n ( t' ,L_0) <  \Gc_n ( L_0)  = \Gc_n (l_m) 
$$
and thus for $m \gg 1$ we have $G_n ( t_m , l_m) < \Gc ( l_m)$ which contradicts condition iii).

In summary, the right lower semi-continuity of the energy release $\G_n$ is necessary for the existence of a quasi-static evolution in the sense of Proposition \ref{p.existsn}. 

Finally, note that the representation \eqref{e.ell+} is independent of the extension of $\G_n$ and $\Gc_n$ in $\Lambda_n \setminus \{ L \}$; indeed, \eqref{e.ell+} can be written also as 
$$
 \ell_n (t) = \inf \{ l \in [L_0,L] \setminus \Lambda_n  \text{ or } l =L  : G_n (t , l) < \Gc_n (l) \} . 
$$
\end{remark}

\subsection{Energy identity}

In the sequel the energy balance will be crucial to explain toughening. To this end,  we introduce the potential energy $F : [0,T] \times [L_0,L] \to \R$ given by
$$ F_n (t, l) = E_n ( t, l) + D_n (l ) , $$ 
where $D_n :[L_0,L] \to \R$ is the dissipated energy given by 
$$
	D_n (l) =  \int_{L_0}^l \Gc_n (s) \, ds  .
$$
Note that the energy $F_n ( \cdot, l)$ is differentiable in $[0,T]$ while $F_n ( t, \cdot)$ is differentiable only in $[L_0,L] \setminus \Lambda_n$.

\begin{corollary} \label{c.existsn} Let $\ell_n$ be the quasi-static evolution provided by Proposition \ref{p.existsn}. For every $ t \in [0,T)$ the following energy identity holds:
\begin{equation} \label{e.enid}
	F_n ( t , \ell_n(t)) =  \int_0^t \partial_t F_n  ( \tau , \ell_n (\tau)) \, d\tau  +  \hspace{-12pt} \sum_{\tau \, \in \, J(\ell_n) \cap [0,t]} \llbracket F_n ( \tau , \ell_n (\tau)) \rrbracket ,
\end{equation}
where  $\llbracket F_n ( \tau , \ell_n (\tau)) \rrbracket =  F_n ( \tau , \ell_n (\tau)) -  F_n ( \tau , \ell_n^- (\tau)) \le 0$.  
\end{corollary}

\separe

We recall that $F_n ( 0 , L_0 )=0$. The energy identity \eqref{e.enid} can be written also as 
\begin{equation} \label{e.enid2} 
	E_n (t, \ell_n(t)) =  \int_0^t P_n ( t, \ell_n (t)) \, dt - D_n ( \ell_n(t))   + \hspace{-12pt} \sum_{\tau \, \in \, J(\ell_n) \cap [0,t]}  \llbracket F_n ( \tau , \ell_n (\tau))  \rrbracket ,
\end{equation}
where $P_n ( t, \ell_n (t))  = \partial_t F ( t , l) = \partial_t E ( t, l)$ is the power of external forces. 
Therefore, part of the energy (supplied by the work of external forces) is stored in the elastic energy, part is dissipated by the crack, while part is dissipated in the discontinuity points. In this respect, we recall that the rate-independent evolutions $\ell_n$ can be obtained also from rate-dependent evolutions (see e.g.~\cite{KneesMielkZan08}) by vanishing viscosity; in this approach,  the energy dissipated in the jumps turns out to be the limit of the energy dissipated in the rate-dependent processes.

\separe

Clearly, by right continuity the energy identity holds also in every time interval $[t_\1, t_\2]$, in the form
$$
	F_n ( t_\2 , \ell_n(t_\2)) =  F_n ( t_\1 , \ell_n(t_\1))  + \int_{t_1}^{t_2} \partial_t F_n  ( \tau , \ell_n (\tau)) \, d\tau  +  \hspace{-12pt} \sum_{\tau \, \in \, J(\ell_n) \cap (t_1,t_2]} \llbracket F_n ( \tau , \ell_n (\tau)) \rrbracket .
$$

\separe

\proof In general the energy $F_n ( t, \cdot)$ is not differentiable in $\Lambda_n$, thus the idea is to consider the sub-intervals in which the evolution takes values in $[L_0,L] \setminus \Lambda_n$, taking care of the possible discontinuities. 
In this way we also make reference to the stick-slip and re-nucleation effect at the interfaces, see \cite{HossainHsuehBourdinBhattacharya_JMPS14}. 

Let $l_{n,k} = k l_n = k L /n$. Let $k_0$ and $k_I$ such that 
$$
	 l_{n,k_0} \le L_0 < l_{n,k_1} , \qquad l_{n,k_I-1} < \sup \ell_n \le l_{n,k_I} .
$$
If $I=0$ we have $\ell_n (t)  \equiv L_0$, thus the energy identity \eqref{e.enid} boils down to 
$$
	F_n (t, \ell_n(t)) =  \int_0^t \partial_t F_n  ( \tau , \ell_n (\tau)) \, d\tau  ,
$$
which is clearly true by the time regularity of the energy. If $I >0$, consider the points $l_{n,k} \in \Lambda_n$ for $k=k_0, ..., k_I$ and define the times
$$
	t_{n,k_0} = \sup \{ t : \ell_n (t) \le L_0 \} , 
	\quad
	t^\flat_{n,k_I} = \sup \{ t : \ell_n (t) < l_{n,k_I} \} , 
	\quad 
	t^\sharp_{n,k_I} = T , 
$$
$$
	t_{n,k_i}^\flat = \sup \{ t : \ell_n (t) < l_{n,k_i} \} , 
		\quad
			t_{n,k_i}^\sharp = \inf \{ t : \ell_n (t) > l_{n,k_i} \} 
			\quad \text{for $0 < i < I$ (if any).}
$$
In this way we have a finite partition of $[0,T]$ given by
$$
	0 < t_{n,k_0} \le t_{n,k_\1}^\flat \le t_{n,k_\1}^\sharp \le ... \le t_{n,k_I}^\flat \le t_{n,k_I}^\sharp = T ;
$$
note that some of these points may coincide. We will prove \eqref{e.enid} by induction, for $i=0,...,I$.

\separe

In the interval $[0,t_{n,k_0})$ we have $\ell_n (t) = L_0$ .  Thus for every $t \in [0,t_{n,k_0})$ we can write 
$$
	F_n (t ,\ell_n(t) ) 
	= \int_0^t \partial_t F_n (\tau , \ell_n (\tau)) \, d \tau .
$$
Passing to the limit as $t \nearrow t_{n,k_0}$ we get
$$
	F_n  (t_{n,k_0} ,\ell_n^- (t_{n,k_0}) ) =  \int_0^{t_{n,k_0}} \partial_t F_n ( \tau , \ell_n (\tau)) \, d \tau 
$$
and then, by right continuity, 
\begin{equation} \label{e.prima}
       F_n  (t_{n,k_0} ,\ell_n (t_{n,k_0}) ) =  \int_0^{t_{n,k_0}} \partial_t F_n ( \tau , \ell_n (\tau)) \, d \tau 
       + \jump{ F_n ( t_{n,k_0} , \ell_n (t_{n,k_0} ))}  .
\end{equation}

\separe

If $t^\flat_{n,k_1} > t_{n,k_0}$ then $\ell_n (t) \in ( L_0 , l_{n,k_1} ) \subset ( l_{n,k_0} , l_{n,k_1})$ for $ t \in ( t_{n,k_0} , t^\flat_{n,k_1})$. In the interval $(L_0, l_{n,k_1})$ the elastic energy $\E_n ( t, \cdot)$ is of class $C^1$ while the dissipated energy $D_n$ is affine, since $G_n^c$ is constant. Thus in every subinterval $ t_{n,k_0} < t' < t <  t^\flat_{n,k_1}$ we can apply the chain rule in $BV$ to get 
\begin{align*}
	F_n (t , \ell_n (t ) ) &  = F_n (t', \ell_n(t')) \, + \int_{t'}^{t} \partial_t F_n (\tau, \ell_n (\tau) ) \, d \tau \, + \int_{t'}^{t}  \partial_l F_n ( \tau , \ell_n (\tau) ) \, d_{LC} \ell_n (\tau)  \\ & \quad + \sum_{ \tau \in J (\ell_n) \cap (t', t ]}  \llbracket F_n ( \tau , \ell_n (\tau)) \rrbracket ,
\end{align*}
where $d_{LC} \ell_n$ denotes the sum of the absolutely continuous and Cantor part of the measure $d \ell_n$. Note that $d_{LC} \ell_n$ is concentrated on the set of continuity points of $\ell_n$. We claim that 
$$
	\int_{t'}^{t}  \partial_l F_n ( \tau , \ell_n (\tau) ) \, d_{LC} \ell_n (\tau) = 0 .
$$
Indeed, for $l \in ( L_0 , l_{n,k_1} )$ we have $\partial_l F_n ( \tau , l ) = - G_n (\tau , l ) + G_n^c(l)$. By Proposition \ref{p.existsn} we known that $ G_n (\tau , \ell_n(\tau) ) \le  G_n^c(\ell_n(\tau))$. If equality holds then $\partial_l F_n ( \tau , \ell_n (\tau) ) = 0$. If $ G_n (\tau , \ell_n(\tau) ) < G_n^c(\ell_n(\tau))$ then 
by Corollary \ref{c.dotell} it follows that $\ell_n$ is constant in a neighborhood of $\tau$, and thus $d_{LC} \ell_n (\tau) = 0$. The claim is proved. Now, passing to the limit as $t'\to t_{n,k_0}$ and using \eqref{e.prima} we get 
\begin{align}
	F_n (t , \ell_n ( t ) ) 
		& = F_n ( t_{n,k_0} , \ell ( t_{n,k_0} ) ) + \int_{t_{n,k_0} }^{t} \partial_t F_n (\tau, \ell_n (\tau) ) \, d \tau 	
			 + \sum_{ \tau \in J (\ell_n) \cap (t_{n,k_0} , t ] }  \llbracket F_n ( \tau , \ell (\tau)) \rrbracket  \nonumber \\		
		& = \int_{0}^{t} \partial_t F_n (\tau, \ell_n (\tau) ) \, d \tau 	
			+ \sum_{ \tau \in J (\ell_n) \cap [0 , t ] }  \llbracket F_n ( \tau , \ell (\tau)) \rrbracket  .  \label{e.seconda}
\end{align}
Passing to the limit as $t \to t_{n,k_1}^\flat$ and taking into account the possible jump in $t_{n,k_1}^\flat$ yields
\begin{align}
	F_n (t_{n,k_1}^\flat , \ell_n ( t_{n,k_1}^\flat ) ) = \int_{0}^{t_{n,k_1}^\flat} \partial_t F_n (\tau, \ell_n (\tau) ) \, d \tau 	
			+ \sum_{ \tau \in J (\ell_n) \cap [0 , t_{n,k_1}^\flat ] }  \llbracket F_n ( \tau , \ell (\tau)) \rrbracket  .  \label{e.terza}
\end{align}
If $t^\flat_{n,k_1} = t_{n,k_0}$ there is nothing to prove since the previous identity boils down to \eqref{e.prima}. 

\separe

If $t^\flat_{n,k_1} <  t^\sharp_{n,k_1}$ then $\ell_n = l_{n,k_1} $ in the interval $(t^\flat_{n,k_1} , t^\sharp_{n,k_1} )$, thus for every $t \in ( t^\flat_{n,k_1} , t^\sharp_{n,k_1} )$ we can write
\begin{align}
F (t , \ell_n ( t)  ) = 
& \int_{0}^{t} \partial_t F_n (\tau, \ell_n (\tau) ) \, d \tau 	
			+ \sum_{ \tau \in J (\ell_n) \cap [0 , t  ] }  \llbracket F_n ( \tau , \ell (\tau)) \rrbracket  .  \label{e.quarta} 
\end{align}
If $t^\flat_{n,k_1} =  t^\sharp_{n,k_1}$ then \eqref{e.terza} holds replacing $t^\flat_{n,k_1}$ with  $t^\sharp_{n,k_1}$. We proceed in this way up to any time $t <T$. \qed

\subsection{Compactness and convergence of evolutions}

By Helly's Theorem it is well known that up to non-relabeled subsequences the evolutions $\ell_n$ converge to a certain limit $\ell$. The goal of the next proposition is to provide a first characterization of the limit $\ell$ independently of the convergence of the energies. 

\begin{proposition}  \label{p.ell+-1} Let $\ell_n$ be the sequence of quasi-static evolutions given by Proposition \ref{p.existsn}. There exists a subsequence (not relabeled) such that $\ell_n \to \ell$ pointwise in $[0,T)$. For $t\in [0,T)$ the right continuous representative of the limit $\ell$ is characterized by 
\begin{equation} \label{e.ell+-1}
 \ell^+ (t) =  \inf \left\{ l \in [L_0,L] :  \Gamma\text{-}\liminf_{n \to \infty} \frac{\G_n}{ \Gc_n} (l) < \frac{1}{f^2 (t)} \right\} ,
\end{equation}
where $1 / f^2 (0) = +\infty$. 
\end{proposition}

\proof We recall that 
$$
     \Gamma\text{-}\liminf_{n \to \infty} \frac{\G_n}{ \Gc_n} (l) = \inf \Big\{   \liminf_{n \to \infty} \frac{\G_n (l_n)}{ \Gc_n (l_n)}  : l_n \to l \Big\} .
$$
For the theory of $\Gamma$-convergence we refer the reader to \cite{DalMaso93}. 

Let $t > 0$. By Proposition \ref{p.existsn} we know that 
$G_n (t , \ell_n(t) ) \le  \Gc_n ( \ell_n(t) )$, hence $G_n (t , \ell_n(t) )$ is finite. Writing $G_n (t , \ell_n(t) ) = f^2 (t) \, \mathcal{G}_n ( \ell_n (t ))$ we get
$$
	\frac{1}{f^\2 (t)} \ge \frac{\G_n (\ell_n(t))}{ \Gc_n (\ell_n( t) )} .
$$
The above inequality holds also in generalized sense for $t=0$.
Since $\ell_n (t) \to \ell(t)$ pointwise for every $t \in [0,T)$ we have
$$
	\frac{1}{f^2 (t)} \ge \liminf_{n \to \infty} \frac{ \G_n (\ell_n(t))}{ \Gc_n (\ell_n( t) )} \ge \inf \Big\{   \liminf_{n \to \infty} \frac{\G_n (l_n)}{ \Gc_n (l_n)}  : l_n \to \ell (t) \Big\} = \Gamma\text{-}\liminf_{n \to \infty} \frac{\G_n}{ \Gc_n} (\ell(t)) .
$$

Let $\tau > t \in [0,T)$. By the monotonicity of $f$ we can write 
$$
   \frac{1}{f^2(t)} > \frac{1}{f^2 (\tau)} \ge \Gamma\text{-}\liminf_{n \to \infty} \frac{\G_n}{ \Gc_n} (\ell(\tau)) . 
$$
Hence 
$$  
	 \ell(\tau) \in  \left\{ l \in [L_0,L] : \Gamma\text{-}\liminf_{n \to \infty} \frac{\G_n}{ \Gc_n} (l) < \frac{1}{f^2 (t)} \right\} 
	 \quad \text{ for every $\tau > t$} 
	 $$ 
and then 
$$
	\ell^+ (t) \ge \inf \left\{ l \in [L_0,L] : \Gamma\text{-}\liminf_{n \to \infty} \frac{\G_n}{ \Gc_n} (l) < \frac{1}{f^2 (t)} \right\}  . 
$$

\separe

Let us prove the opposite inequality. If $\ell^+ (t) = L_0$ there is nothing to prove. Otherwise, let $L_0 \le l < \ell^+(t) \le \ell (\tau)$ for $\tau > t$. Let $\tau_n \searrow t^+$ such that $\ell_n (\tau_n) \to \ell^+(t)$ and let $l_n$ be a ``recovery sequence'' for $l$, i.e., $l_n \to l$ and 
$$	\liminf_{n \to \infty}  \frac{\G_n (l_n)}{  \Gc_n (l_n) } =  \Gamma\text{-}\liminf_{n \to \infty} \frac{\G_n}{ \Gc_n} (l) .
$$
From the representation \eqref{e.ell+} of the evolution $\ell_n$ we know that
\begin{align*}  
	\ell_n (\tau_n) 
	& = \inf \{ l \in [L_0,L] : G_n ( \tau_n , l) <  \Gc_n (l) \}  \\ 
	& = \inf \left\{ l \in [L_0,L] : \frac{1}{f^2 (\tau_n)} > \frac{\G_n ( l) }{ \Gc_n (l)} \right\} .
 \end{align*}
Being $l_n < \ell_n (\tau_n)$ for $n \gg1$ (since $l < \ell^+(t)$) the above representation implies that 
$$
	\frac{1}{f^2 (\tau_n)} \le  \frac{\G_n ( l_n) }{ \Gc_n (l_n)} .
$$
Taking the liminf as $n \to \infty$ and remembering that $l_n$ is a recovery sequence, we get 
$$
	\frac{1}{f^2 (t)} \le \liminf_{n \to \infty} \frac{\G_n (l_n)}{ \Gc_n (l_n)} =  \Gamma\text{-}\liminf_{n \to \infty} \frac{\G_n}{ \Gc_n} (l) .
$$
Thus, every $l < \ell^+ (t)$ does not belong to the set 
$$
\left\{ l \in [L_0,L] :  \Gamma\text{-}\liminf_{n \to \infty} \frac{\G_n}{ \Gc_n} (l) < \frac{1}{f^2 (t)} \right\} 
$$
and then
$$  \ell^+(t) \le \inf  \left\{ l \in [L_0,L] :  \Gamma\text{-}\liminf_{n \to \infty} \frac{\G_n}{ \Gc_n} (l) < \frac{1}{f^2 (t)} \right\} , $$
which concludes the proof.  \qed

\begin{remark} \label{r.sub-seq} \normalfont The characterization  \eqref{e.ell+-1} seems to be the natural choice combined with the pointwise convergence of a subsequence (extracted  in Proposition \ref{p.ell+-1}) of the quasi-static evolutions $\ell_n$. 
However, if there exists 
\begin{align}   \label{e.Glimex}
\Gamma\text{-}\hspace{-4pt}\lim_{n \to \infty} \frac{\G_n^{\phantom{I}}}{\Gc_n} (l) 
\end{align}
then the whole sequence $\ell_n \weakstarto \ell$ in $BV(0,T)$. Indeed, each function $\ell_{n}$ is monotone and bounded in $[L_0,L]$, hence the sequence $\ell_n$ is weakly* pre-compact in $BV(0,T)$. Given any subsequence $\ell_{n_k}$ of $\ell_n$ there exists a further subsequence, denoted by $\ell_{n_i}$, converging to some limit $\ell$ weakly in $BV (0,T)$ and pointwise in $(0,T)$. 
By Proposition \ref{p.ell+-1} the limit $\ell$ (possibly depending on the subsequence) is characterized by 
$$
	  \ell^+ (t) =  \inf \left\{ l \in [L_0,L] :  \Gamma\text{-}\liminf_{n \to \infty} \frac{\G_{n_i}}{ \Gc_{n_i}} (l) < \frac{1}{f^2 (t)} \right\} 
	  	= \inf \left\{ l \in [L_0,L] :  \Gamma\text{-}\hspace{-4pt}\lim_{n \to \infty} \frac{\G_{n}}{ \Gc_{n}} (l) < \frac{1}{f^2 (t)} \right\} .
$$
Therefore, the evolution $\ell$ is actually independent of the subsequence $\ell_{n_k}$. 
Clearly, this argument requires the existence of the $\Gamma$-limit in \eqref{e.Glimex} while the $\Gamma$-liminf in \eqref{e.ell+-1} is always defined. 
\end{remark}

\subsection{Homogenization of the elastic energy}

In this section we will first recall a few results on the homogenization of the elastic energy. Let us introduce the homogenized stiffness matrix $\C_\hom$ given by 
$$
 \C_{\hom} =  \left( \begin{matrix}  \mu_{\hom, \1} & 0  \\ 0 & \mu_{\hom, \2}  \end{matrix}  \right) ,
$$ 
where $1 / \mu_{\hom,\1} = \lambda / \mu_{\A,\1} + (1-\lambda) / \mu_{\B,\1}$ is the weak* limit of $1/\mu_{n,\1}$ and 
$\mu_{\hom, \2} =  \lambda \mu_{\A,\2} + ( 1- \lambda ) \mu_{\B,\2} $
 is the weak* limit of $\mu_{n,\2}$. Note that the coefficients $\mu_{\hom,i}$ are constant and independent of $l \in (0,L)$. We denote by $\W_{l,\hom}$ the homogenized elastic energy and by $\E_\hom$ its static condensation. By classical results on the homogenization of elliptic problems (see Appendix \ref{a.B}) and on the energy release (see e.g.~\cite{N_AMO17} and the references therein) we infer the following proposition. 

\begin{proposition} \label{p.homen-vert} If $l_n \to l$ then $\E_n (  l_n) \to \E_\hom (l)$. The energy $\E_\hom$ is of class $C^1 (0,L)$. 
\end{proposition}


Thanks to the regularity of $\E_\hom$ we can define the energy release 
$$
	\G_\hom (l) = - \partial_l^+ \E_\hom (l) =  - \partial_l \E_\hom (l)  \quad \text{for $l \in (0,L)$,} \qquad \G_\hom (L) = 0 , 
$$
which is continuous in $(0,L)$. We remark that in general $\G_\hom$ is neither the a.e.~limit nor the $\Gamma$-liminf of $\G_n$ (an example is given in \S\ref{explicit}). However, thanks to the regularity and the convergence of the energies we have the following results, which will be useful in the sequel. 

\begin{lemma} \label{l.GntoG_1} Let $\G_{n_i}$ be a subsequence of $\G_n$. For a.e.~$l \in (0,L)$ it holds
$$
	\liminf_{i \to \infty} \G_{n_i} (l) \le \G_\hom (l)
	\le \limsup_{n \to \infty} \G_{n_i} (l)  .
$$
(Note that the set where the above estimate holds may depend on the subsequence).
\end{lemma}

\proof Let $0 < l_1 < l_2 < L$. By Proposition \ref{p.homen-vert} and by the regularity of the energies we can write
$$
	\int_{l_1}^{l_2}  \G_{n_i} (s) \, ds  = \E_{n_i} (l_1) - \E_{n_i} (l_2) 
		\quad  \to \quad 
	\E_\hom (l_1) - \E_\hom (l_2) = \int_{l_1}^{l_2}  \G_\hom (s) \, ds 	. 
$$	
Hence, by Fatou's Lemma we also have
$$
	\int_{l_1}^{l_2}  \liminf_{i \to \infty} \G_{n_i} (s) \, ds  \le \liminf_{i \to \infty} \int_{l_1}^{l_2}   \G_{n_i} (s) \, ds 
	= \int_{l_1}^{l_2}  \G_\hom (s) \, ds  ,
$$
which proves the lower inequality by the arbitrariness of $l_1$ and $l_2$. 

\separe

For the upper estimate we can argue in a similar way. \qed

\begin{lemma} \label{l.GntoG_2} For every $l \in (0,L)$ there exist
\begin{itemize}
\item a sequence $l^\sharp_n \to l^+$ such that  $ \G_\hom (l) \ge \limsup_{n \to \infty} \G_n (l^\sharp_n)$,
\item a sequence $l^\flat_n \to l^+$ such that  $ \G_\hom (l) \le \liminf_{n \to \infty} \G_n (l^\flat_n)$.
\end{itemize}
\end{lemma}

\proof  By the continuity of $\G_\hom$, for every $\delta >0$ there exists $\eps>0$ such that $\G_\hom (s) < \G_\hom (l) + \delta $ for every $s \in [l, l+\eps)$. 

We claim that there exists $m \in \mathbb{N}$ such that for every $n > m$ we can find  $l_n \in ( l , l + \eps )$ with $\G_n (l_n) < \G_\hom (l) + \delta$; from the claim the existence of $l^\sharp_n$ follows.  By contradiction, assume that for every $m \in \mathbb{N}$ there exists $n_m > m$   such that   $\G_{n_m} (s) \ge \G_\hom (l) + \delta$ for every  $s \in ( l  , l + \eps )$, then we can find a subsequence such that 
$ \liminf_{m \to \infty} \G_{n_m} (s)  \ge \G_\hom (l) + \delta > \G_\hom (s)$ in $ ( l  , l + \eps )$, which is a contradiction with Lemma \ref{l.GntoG_1}.

The proof for $\l^\flat_n$ is analogous. \qed 

From the previous lemma we get the following corollary.

\begin{corollary} \label{c.gamma} 
If $\G_\hom \le\Gamma\text{-}\liminf_{n \to \infty} \G_n$ then $\G_\hom = \Gamma\text{-}\lim_{n \to \infty} \G_n$.
\end{corollary}

\proof By Lemma \ref{l.GntoG_2} for every $l \in (0,L)$ we have
$$
	\Gamma\text{-}\limsup_{n \to \infty } \G_n (l) \le \limsup_{n \to \infty} \G_n (l_n^\sharp) \le \G_\hom (l) \le\Gamma\text{-}\liminf_{n \to \infty} \G_n (l) , 
$$
which concludes the proof.  \qed
 
As in the previous sections we consider the energy $E_\hom (t, l) = f^2(t) \E_\hom (l)$. Most important, we define $G_\hom (t, l) = \partial_l E_\hom (t, l) = f^2(t) \G_\hom (l)$ for every $t \in [0,T]$ and $l\in (0,L)$; note that $\G_\hom (l)$ is finite for $l \in (0,L)$. Finally, we define $G_\hom ( t, L)=0$ for every $t \in [0,T]$.

\subsection{Effective toughness and homogenization of Griffith's criterion}

The aim of this section is to provide an {\it effective toughness} $\Gceff$ in such a way that the limit evolution $\ell$ (given by Proposition \ref{p.ell+-1}) satisfies Griffith's criterion for $\Ghom$ and $\Gceff$. Formally, regardless of the regularity of the energies involved, we may define $G^c_\eff$ as the function which makes the following identity true:
\begin{equation} \label{e.myformula}
     \frac{\G_\hom}{\Gc_\eff} (l) = \Gamma\text{-}\liminf_{n \to \infty}  \frac{\G_n}{ G^c_n} (l)  . 
\end{equation}
More precisely, we define $\Gc_\eff : [0,L] \to [ 0, +\infty]$ as
\begin{align}  \label{e.Ghomgen}
\Gceff (l) = \begin{cases}   
                        {\displaystyle   \frac{\G_\hom (l) }{\displaystyle \Gamma\text{-}\liminf_{n \to \infty} \frac{\G_n^{\phantom{I}}}{\Gc_n} (l) }  } &   \text{ if $\G_\hom (l) \neq 0 $,} \\
                        + \infty &   \text{ if $\G_\hom (l) = 0 $}.
			\end{cases}
\end{align}
In the above definition we assume $\Gceff (l)=+\infty$ when $ \displaystyle \Gamma\text{-}\liminf_{n \to \infty} \frac{\G_n^{\phantom{I}}}{\Gc_n} (l) = 0$; by the proof of Lemma \ref{l.propGhom} (see below) it follows that $\displaystyle \Gamma\text{-}\liminf_{n \to \infty} \frac{\G_n^{\phantom{I}}}{\Gc_n} (l) < + \infty$. As far as Griffith's criterion is concerned, the definition of $\Gc_\eff (l)$ when $\G_\hom (l) = 0$ is in some sense arbitrary, indeed, if the evolution takes, at time $t$,  a value $l$ where $\G_\hom (l)=0$ then the inequality $G_\hom (\tau, l) \le \Gc_\eff (l)$ is satisfied for every $\tau \ge t$ and every non-negative value of $\Gc_\eff$, further, the evolution will be constantly equal to $l$  for $\tau \ge t$ (see Remark \ref{r.defGeff} for a different definition of $\Gc_\eff$ when $\G_\hom$ vanishes).

Note that $ \Gceff$ is independent of time and of the evolution, however in principle it may depend on every parameter and every datum of the problem; in mechanical terms, $\Gceff$ is an R-curve (R standing for resistance, see e.g.~\cite{And95}). In \S\ref{explicit} we will give an example in which $ \Gceff$ is actually a  constant, depending on the volume fraction $\lambda$, the ratio between the shear moduli $\mu_{\A,i}$ and $\mu_{\B,i}$ and the toughness $\Gc_\A$ and $\Gc_\B$. Before proceeding, we provide this lemma.

\begin{lemma} \label{l.propGhom} The effective toughness $ \Gceff$ is upper semi-continuous (and thus Borel measurable). Moreover, $\Gceff \ge \min \{ \GcA , \GcB \}$.  \end{lemma} 

\proof It is enough to consider the set $\{ \G_\hom \neq 0 \}$ since, by continuity, the set $\{ \G_\hom = 0 \}$ is closed. The upper semi-continuity follows from the continuity of $\G_\hom$ and the lower semi-continuity of the $\Gamma\text{-}\liminf$. Next, write
\begin{align*}
	\Gamma\text{-}\liminf_{n \to \infty} \frac{\G_n^{\phantom{I}}}{\Gc_n} (l)  
	& = 
	\inf \Big\{  \liminf_{n \to \infty} \frac{ \G_n (l_n)^{\phantom I}}{\Gc_n (l_n)^{\phantom I}}  : l_n \to l \Big\}   \\
	& \le \frac{1}{\min \{ \GcA , \GcB \}} \inf \left\{  \liminf_{n \to \infty}  \G_n (l_n) : l_n \to l \right\}
	\\ & \le \frac{\Gamma\text{-}\liminf_{n \to \infty} \G_n (l)}{\min \{ \GcA , \GcB \}} .
\end{align*}
By Lemma \ref{l.GntoG_1} we infer 
\begin{equation}  \label{e.23se}
 \Gamma\text{-}\liminf_{n \to \infty} \G_n (l) \le \G_\hom (l)
 \end{equation} 
and thus
\begin{equation}  \label{e.upx}
	\Gamma\text{-}\liminf_{n \to \infty} \frac{\G_n^{\phantom{I}}}{\Gc_n} (l)   
	\le \frac{ \G_\hom (l)}{\min \{ \GcA , \GcB \}} .
\end{equation}
Hence, $\Gceff (l) \ge \min \{ \GcA , \GcB \}$.  \qed

A local upper bound for $\Gceff$ is given in Corollary \ref{c.--}, as a byproduct of the following theorem.

\begin{theorem} \label{t.genconv} Let $\ell$ be the limit evolution given by Proposition \ref{p.ell+-1}. 
Then $\ell^+ (t) = \ell_\hom (t)$, where $\ell_\hom: [0,T) \to [L_0, L]$ is the unique non-decreasing, right continuous function which satisfies the initial condition $\ell_\hom (0) = L_0$ and Griffith's criterion: 
\begin{itemize}
\item[i)] $\Ghom (t, \ell_\hom(t) ) \le   \Gceff (\ell_\hom(t))$ for every time $t \in [0,T)$; 
\item[ii)]  if $G_\hom ( t, \ell_\hom (t) ) <  \Gceff (\ell_\hom (t)) $ then $\ell_\hom$ is right differentiable in $t$ and $\dot{\ell}_\hom^+ (t)=0$; 
 \item[iii)] if $t \in J ( \ell_\hom)$ then $G_\hom (t , s ) \ge  \Gceff (s)$ for every $s \in [\ell^-_\hom(t) , \ell_\hom(t) )$.
\end{itemize}
\end{theorem}

\proof First, we check that for every $t \in [0,T)$ we have 
\begin{align} \label{e.zuz}
	\{ l \in [L_0,L] : G_\hom (t, l) <  \Gceff (l) \} & = 
      \left\{ l \in [L_0,L] :  \Gamma\text{-}\liminf_{n \to \infty} \frac{\G_n}{ \Gc_n} (l) < \frac{1}{f^2 (t)}  \right\}  . 
\end{align}
If $t=0$ then $G_\hom ( 0 , l ) = 0$ and thus by Lemma \ref{l.propGhom} we have $G_\hom ( 0 , l ) <  \min \{ \Gc_\A, \Gc_\B \} \le \Gc_\eff (l)$ for every $l \in [L_0,L]$. By \eqref{e.upx} we can write 
$$
	\Gamma\text{-}\liminf_{n \to \infty} \frac{\G_n}{ \Gc_n} (l) < +\infty = \frac{1}{f^2 (0)} 
$$
for every $l \in [L_0,L]$. If $t>0$ we write $G_\hom ( t, l) = f^2 (t) \G_\hom (l)$. If $\G_\hom (l) \neq 0$, using the definition of $\Gceff$ we get 
$$
	G_\hom (t, l) <  \Gceff (l)  \quad \Leftrightarrow \quad
	f^2 (t) <  
      \frac{1}{\displaystyle \Gamma\text{-}\liminf_{n \to \infty} \frac{ \G_n^{\phantom{I}} }{\Gc_n} (l) } .
$$ 
If  $\G_\hom (l) = 0$ then $G_\hom (t, l) = 0$ and thus $G_\hom (t, l) <  \min \{ \Gc_\A, \Gc_\B \} \le \Gceff (l)$. 
Moreover by \eqref{e.upx} we have 
$$	\Gamma\text{-}\liminf_{n \to \infty} \frac{\G_n}{ \Gc_n} (l) = 0 <  \frac{1}{f^2 (t)} .  $$

By \eqref{e.zuz} we can now invoke Proposition \ref{p.ell+-1} which gives the representation
\begin{equation} \label{e.1507}
	\ell^+ ( t) = \inf \{ l \in [L_0,L] : G_\hom (t, l) <  \Gceff (l) \} . 
\end{equation}
At this point it is enough to follow the proof of Proposition \ref{p.existsn}, with a few minor differences, in order to show that the function $\ell^+$ defined by \eqref{e.1507} satisfies Griffith's criterion. 
The proof is concluded. \qed

\begin{remark}  \normalfont
In general, the limit evolution $\ell$ does not satisfy Griffith's criterion in terms of the energy release $G_\hom$, of the homogenized energy $E_\hom$, and the homogenized toughness $\Gc_\hom = \lambda \Gc_\A + (1-\lambda) \Gc_\B$. Starting from the \eqref{e.ell+-1}, which is independent of the limit energies, there are actually different ways to proceed.  In our approach we define an effective toughness $\Gc_\eff$  in such a way that the limit $\ell$ satisfies Griffith's criterion for $G_\hom$ and $\Gc_\eff$. This idea seems very natural in terms of convergence of the energies and is consistent with the approach followed in mechanics, see e.g.~\cite{HossainHsuehBourdinBhattacharya_JMPS14}. On the other hand it could be interesting also to fix $\Gc_\hom$ and define an effective energy release $\G_\eff$ taking the limit of $\G_n$. However, even if $\G_n$ has a ``nice'' integral representation (see Lemma \ref{l.En}) it seems difficult to study its limit, the main issue being the fact that the support of the auxiliary functions $\phi$ shrinks with $n$. In \S\ref{ss.conv} we will consider a different setting in which the representation of $\G_n$ holds for a function $\phi$ independent of $n$ (see Lemma \ref{l.Enhoriz}); in that case the limit of $\G_n$ can be computed, and actually equals $\G_\hom$. 
\end{remark}

\subsubsection[Evolutionary $\Gamma$-convergence]{Evolutionary {$\boldsymbol \Gamma$}-convergence}  \label{a.EvGconv}

In this section we give a slight generalization of Theorem \ref{t.genconv}, in order to provide an evolutionary $\Gamma$-convergence result for our rate independent system, in the spirit of \cite{Mielke_LNAMM16}. To this end, let us introduce the rate independent dissipation functionals 
$$
	\mathcal{R}_n ( l , \dot{l} ) = \begin{cases} \Gc_n ( l ) \, \dot{l} & \dot{l} \ge 0 , \\
	+ \infty & \text{otherwise,}  \end{cases}
	\qquad
	\mathcal{R}_\eff ( l , \dot{l} ) = \begin{cases} \Gc_\eff (l) \, \dot{l}  & \dot{l} \ge 0 , \\
	+ \infty & \text{otherwise.}  \end{cases}
$$
We refer the reader to Appendix \ref{a.A} for the notion of balanced viscosity solution for the rate independent system $(E_n , \mathcal{R}_n)$, in the language of \cite{MielkeRossiSavare_COCV12}; actually (see again Appendix \ref{a.A}), these solutions  coincide with the quasi-static evolutions of Corollary \ref{c.unica}. Therefore
the right continuous representative of any balanced viscosity solution of the rate independent system $(E_n , \mathcal{R}_n)$ coincides with $\ell_n$, i.e., the unique  right continuous evolution provided by Proposition \ref{p.existsn}.

\begin{corollary} \label{c.evol-G-conv} The rate independent systems $( E_n , \mathcal{R}_n)$ evolutionary $\Gamma$-converge to $( E_\hom , \mathcal{R}_\eff)$. 
\end{corollary}

\proof By \cite[Definition 2.10]{Mielke_LNAMM16} the above statement is equivalent to the following: let $L_n \in (0,L)$ such that $L_n \to L_0 \in (0,L)$  and let $\gamma_n$ be a solution of the rate independent system $(E_n , \mathcal{R}_n)$ 
with intial condition $L_n$, 
then $\gamma_n$ converge to $\gamma_\hom$ where $\gamma_\hom$ is a solution of the rate independent system $(E_\hom , \mathcal{R}_\eff)$ with initial condition $L_0$. 

By Helly's Theorem there exists a subsequence (not relabeled) such that $\gamma_n \to \gamma$ pointwise in $[0,T)$. In order to characterize the limit $\gamma$ it is not restrictive to consider the convergence of the sequence $\ell_n$ of the right continuous representatives of $\gamma_n$. At this point, in order to apply Theorem \ref{t.genconv} it is enough to ``fix the initial condition''. There are at least a couple of alternatives. The first is to follow, step by step, the proof of Theorem \ref{t.genconv} employing the explicit representation 
$$
\ell_n (t) = \inf \{ l \in [L_n,L] : G_n (t , l) < \Gc_n (l) \} .  
$$
The second chance is instead to employ the change of variable $\psi_n: [0,L] \to [0,L]$ given by 
$$
      \psi_n (l) = \begin{cases} l  L_0 / L_n  & l < L_n ,  \\ (l - L_n)  (L - L_0)  / (L-L_n)  + L_0  & l \ge L_n ,   \end{cases}
$$
and then consider the evolutions $\tilde\ell_n (t) = \psi_n \circ \ell_n (t)$, which satisfy $\tilde{\ell}_n (0)=L_0$.  The (tedious) check of convergence is omitted. \qed

\subsubsection{An upper bound} 

Thanks to Griffith's criterion, we can provide also a simple a posteriori estimate on the effective toughness. 

\begin{corollary} \label{c.--} Assume that \,$\sup \{ \ell_\hom (t) : t \in [0,T) \} > L_0$. Then, the effective toughness $ \Gceff$ is locally bounded, more precisely
$$
	\Gceff (l ) \le f^2(T)   \G_\hom (l) \frac{ \max \{ \GcA, \GcB \}  } {\min \{ \GcA, \GcB \} } 
$$
for every $l \in \mathrm{co} \, \ell_\hom ([0,T))$, the convex envelope of the range $\ell_\hom ([0,T))$.

\end{corollary}

\proof  The convex envelope coincides with the set $\{ l \in [L_0, L] : l \le \ell_\hom (t) \text{ for some $t \in [0,T)$} \}$ and is used to ``fill the gaps'' in the discontinuity points.

If $l \in [ \ell^-_\hom (t) , \ell_\hom (t) )$ for some $t \in J ( \ell_\hom)$ then $ \Gceff (l) \le G_\hom (t, l ) \le f^2(T)   \G_\hom (l)$. 

Otherwise, let $L_0 < l' < \ell_\hom (t)$ for some $t \in (0,T)$.  For any sequence $l'_n \to l'$, we argue as follows. Clearly, $l'_n < \ell_n (t)$ for $n \gg 1$. 
Remember that $ \ell_n (t) = \inf \{ l \in [L_0, L]  : G_n (t, l) <  \Gc_n (l)  \} $, hence  $G_n (t, l'_n ) \ge  \Gc_n (l'_n) \ge \min \{ \GcA, \GcB \}$ and thus $\G_n ( l'_n) \ge \min \{ \GcA, \GcB \} / f^2(t) $. As a consequence
\begin{align*}
	\Gamma\text{-}\liminf_{n \to \infty} \frac{\G_n^{\phantom{I}}}{\Gc_n} (l')  
	& = 
	\inf \Big\{  \liminf_{n \to \infty} \frac{ \G_n (l_n)^{\phantom I}}{\Gc_n (l_n)^{\phantom I}}  : l'_n \to l' \Big\}   \\
	& \ge \frac{1}{\max \{ \GcA , \GcB \}} \inf \left\{  \liminf_{n \to \infty}  \G_n (l_n) : l'_n \to l' \right\}
	\\ & \ge \frac{\min \{ \GcA , \GcB \}}{ f^2(t) \max \{ \GcA , \GcB \}} .
\end{align*}
It follows that
$$
	 \Gceff (l') \le \frac{ \max \{ \GcA, \GcB \} f^2(t) \G_\hom (l') } {\min \{ \GcA, \GcB \} }  \le f^2 (T)  \G_\hom (l') \frac{ \max \{ \GcA, \GcB \}  } {\min \{ \GcA, \GcB \} } .
$$
We proved the inequality for any $l' \in (L_0, \sup \{ \ell_\hom (t) : t \in [0,T) \} )$. We conclude by the upper semi-continuity of $\Gc_\eff$. \qed

In \S\ref{explicit} we will see an example in which the effective toughness is strictly greater than $\max \{ \GcA, \GcB \}$.

\begin{remark} \normalfont \label{r.defGeff} Joining the lower bound of Lemma \ref{l.propGhom} and the upper bound of Corollary \ref{c.--} it is clear that the set $\mathrm{co} \, \ell_\hom ([0,T))$ and the set $\{ \G_\hom = 0\}$ are disjoint. Therefore, to the purpose of proving Theorem \ref{t.genconv} it is enough to define $\Gc_\eff$ in the set $\{ \G_\hom \neq 0\}$. In other terms, we could provide alternative definitions of $\Gc_\eff$, which still makes Griffith's criterion true.  For instance, let 
$$
\Gc_\infty (l) =  \frac{\G_\hom (l) }{\displaystyle \Gamma\text{-}\liminf_{n \to \infty} \frac{\G_n^{\phantom{I}}}{\Gc_n} (l) }  
$$
be defined for $l \in [L_0, L]$ such that $\G_\hom (l) \ne 0$.
Then we could employ as effective toughness the function $\tilde{G}^{\hspace{-2.5pt}\text{\sl \fontsize{6}{6}  c}}_\eff (l) = \inf  \{ \Gc (l) \}$ where the infimum is taken over all the upper semi-continuous functions $\Gc : [L_0,L] \to \mathbb{R}$ with $\Gc \ge \Gc_\infty$ in the set $\{ \G_\hom \ne 0  \}$.
\end{remark}

\subsubsection{Effective toughness under convergence of the energy release}

In this subsection we consider a couple of cases in which the energy release $\G_n$ converge (in a suitable sense) to the energy release $\G_\hom$ and, as a consequence, the effective toughness $\Gc_\eff$ is computed explicitly and equals the maximum between $\GcA$ and $\GcB$. 

\begin{corollary} \label{c.elln2ell-a} Assume that $\G_n \to \G_\hom$ locally uniformly in $(0,L)$. If $ \G_\hom (l) \neq 0$ then 
$$   \Gceff (l) 
 = \max \{ \GcA  , \GcB \}  .
 $$ 
\end{corollary}

\proof We recall that 
$$
\Gceff (l) = \frac{\G_\hom (l) }{\displaystyle \Gamma\text{-}\liminf_{n \to \infty} \frac{\G_n^{\phantom{I}}}{\Gc_n} (l) } .
$$
By uniform convergence $\G_n (l_n) \to \G_\hom (l)$ for every $l_n \to l$, hence 
\begin{align*}
	\Gamma\text{-}\liminf_{n \to \infty} \frac{\G_n^{\phantom{I}}}{\Gc_n} (l) & = \inf \Big\{  \liminf_{n \to \infty} \frac{ \G_n (l_n)^{\phantom I}}{\Gc_n (l_n)^{\phantom I}}  : l_n \to l \Big\}   = \G_\hom (l)  \inf \Big\{  \liminf_{n \to \infty} \frac{1}{\Gc_n (l_n)}  : l_n \to l \Big\}  \\ &  = \frac{\G_\hom (l) }{ \max \{ \GcA  , \GcB \}  } ,
\end{align*}
which concludes the proof.  \qed

The previous result can be slightly generalized as follows.

\begin{corollary} \label{c.elln2ell-b} Let $\G_\hom (l) \neq0$. Assume that $\G_\hom (l) \le \Gamma\text{-}\liminf_{n \to \infty} \G_n (l)$ and that there exists a ``joint recovery sequence'' for $\G_\hom (l)$ and $\Gc_\eff(l)$, i.e., $l_n \to l$ such that
$$
	\limsup_{n \to \infty}  \Gc_n (l_n) =  \max \{  \Gc_\A ,  \Gc_\B \} , 
	\qquad
	\limsup_{n \to \infty} \G_n ( l_n ) = \G_\hom(l)  .
$$
Then 
$  \Gceff (l) 
= \max \{  \Gc_\A ,  \Gc_\B \}  $.
\end{corollary}

\proof Clearly $\Gc_n \le \max \{  \Gc_\A ,  \Gc_\B \}$ and $\G_\hom (l) = \Gamma\text{-}\lim_{n \to \infty} \G_n (l)$ by Corollary \ref{c.gamma}. Hence
\begin{align*}
	\Gamma\text{-}\liminf_{n \to \infty} \frac{\G_n^{\phantom{I}}}{\Gc_n} (l) & \ge \frac{1}{\max \{  \Gc_\A ,  \Gc_\B \}} \inf \big\{  \liminf_{n \to \infty} \G_n (l_n) : l_n \to l \big\}   
	= \frac{\G_\hom (l) }{ \max \{ \GcA  , \GcB \} } . 
\end{align*}
By $\Gamma$-convergence we know that $\G_n (l_n) \to \G_\hom (l)$. Hence 
$$
\Gamma\text{-}\liminf_{n \to \infty} \frac{\G_n^{\phantom{I}}}{\Gc_n} (l) \le \liminf_{n \to \infty} \frac{\G_n (l_n)}{\Gc_n(l_n)} = \frac{\G_\hom (l) }{ \max \{ \GcA  , \GcB \} } .
$$
Hence, 
$$
	\frac{\G_\hom (l) }{ \max \{ \GcA  , \GcB \} } =  \Gamma\text{-}\liminf_{n \to \infty} \frac{\G_n^{\phantom{I}}}{\Gc_n} (l) ,
$$
which concludes the proof.  \qed

\separe

\subsection{Energy identity}

\noindent Given \eqref{e.Ghomgen} the effective dissipated energy is defined by
$$
	D_\eff ( l) = \int_{L_0}^l  \Gceff (s) \, ds 
$$
and the potential energy by $F (t, l) = E_\hom (t, l) + D_\eff (l)  $. By Corollary \ref{c.--} it turns out that $D_\eff \in AC ( L_0, \ell_\hom (T) )$.

\separe

\begin{corollary} \label{c.enidhom} For every $ t \in [0,T)$ the following energy identity holds:
\begin{equation} \label{e.enidhom}
	F ( t , \ell_\hom(t)) = \int_0^t \partial_t F  ( \tau , \ell_\hom(\tau)) \, d\tau + \hspace{-12pt} 
	\sum_{\tau \, \in \, J_\smahom \cap [0,t]} \llbracket F ( \tau , \ell_\hom (\tau)) \rrbracket ,
\end{equation}
where $J_\hom = J ( \ell_\hom)$ and $\llbracket F ( \tau , \ell_\hom (\tau)) \rrbracket  = F  ( \tau , \ell_\hom (\tau)) - F ( \tau , \ell_\hom^- (\tau)) \le 0$.
\end{corollary}

\proof In this case we cannot  employ the proof of Corollary \ref{c.existsn}. Clearly, if $\ell_\hom (t) \equiv L_0$ the proof is trivial by time regularity. 
As $\ell_\hom \in BV(0,T)$ we can write 
$$
	d \ell_\hom = \dot{\ell}_\hom d\tau + d_{C} \ell_\hom + d_{J} \ell_\hom , 
$$
where $\dot{\ell}_\hom \in L^1 (0,T)$ denotes the density of the absolutely continuous part (with respect to the Lebesgue measure $d\tau$), $d_{C} \ell$ denotes the Cantor part, while
$$
      d_{J} \ell_\hom  = \sum_{\tau \in J_\smahom} \jump{  \ell_\hom (\tau) } \, \delta_{\tau} .
$$
For convenience, denote $d_{LC} \ell_\hom = \dot\ell_\hom d\tau + d_{C} \ell$.

\separe

Being $E_\hom$ of class $C^1$ by the chain rule in $BV$, see e.g.~\cite{AmbrosioFuscoPallara_00}, it follows that 
\begin{align}
	E_\hom (t , \ell_\hom (t)) & = E_\hom ( 0 , \ell_\hom (0)  ) + \int_{0}^{t}  \partial_t E_\hom ( \tau , \ell_\hom (\tau)    ) \, d\tau + 
	\int_{0}^{t}  \partial_l E_\hom ( \tau , \ell_\hom (\tau)    ) \, d_{LC}  \ell_\hom  \nonumber  \\
	& \phantom{=}  +  \hspace{-10pt} \sum_{\tau \, \in \, J_\smahom \cap [0, t] } \hspace{-10pt} \jump{  E_\hom ( \tau ,  \ell_\hom (\tau)  )  }  . \label{e.eb1}
\end{align}
It remains to prove that a similar representation holds also for the effective dissipated energy, i.e, 
\begin{align}
	D_\eff ( \ell_\hom (t) ) & = D_\eff ( \ell_\hom (0) )  + \int_{0}^{t}   \Gceff ( \ell_\hom (\tau) ) \, d_{LC} \ell_\hom + 
	 \hspace{-12pt} \sum_{\tau \, \in \, J_\smahom \cap [0,t] } \hspace{-6pt} \jump{  D ( \tau ,  \ell_\hom (\tau)  )  }  .  \label{e.eb2}
\end{align}
Remember that $D_\eff$ is not of class $C^1$, however, by definition we can write 
\begin{equation}  \label{e.eb3}
	D_\eff (l_\2) = D_\eff (l_\1) + \int_{l_1}^{l_2}  \Gceff (l) \, dl ,
\end{equation}
for every $0 < l_\1 < l_\2 < L$ ($ \Gceff$ is Borel measurable since it is upper semi-continuous). Hence
$$
	  D_\eff ( \ell_\hom (t) )  = \int_{L_0}^{\ell_\smahom (t)} \Gceff (l) \, dl = \int_{L_0}^{\ell_\smahom (t)} \Gceff (l) \, dl_{|R} +  \int_{L_0}^{\ell_\smahom (t)} \Gceff (l) \, dl_{|U} , 
$$
where $R = \ell_\hom ( [0,T) )$ while $U = \mathrm{co} \, \ell_\hom ([0,T))  \setminus R$. We will check that 
\begin{equation} 
	 \int_{L_0}^{\ell_\smahom(t)}  \Gceff (l)  \, dl_{| U}  = \sum_{\tau \, \in \, J_\smahom \cap [0, t ] } \hspace{-6pt} \jump{  D (\ell_\hom (\tau)  )  } 
	\label{e.eb6}
\end{equation}
and that
\begin{align}
	\int_{L_0}^{\ell_\smahom(t)}  \Gceff (l)  \, dl_{|R} = \int_{0}^{t}   \Gceff ( \ell_\hom(\tau))  \, d_{LC} \ell_\hom 
			\label{e.eb5} ,
\end{align}
which, put together, imply \eqref{e.eb2}. 

\separe

By the right continuity of $\ell_\hom$ 
\begin{equation} \label{e.eb4}
	U = \bigcup_{\tau \in J_\smahom} U_\tau ,
\end{equation}
where  $U_\tau$ denotes an interval of the form $(\ell^-_\hom(\tau) , \ell_\hom (\tau) )$ or $[\ell^-_\hom(\tau) , \ell_\hom (\tau) )$ (the fact that $\ell_\hom (\tau)$ is excluded is due to the right continuity of $\ell_\hom$). Note that the sets $U_\tau$ are pairwise disjoint and remember that $J_\hom$ is a countable set. Then 
\begin{align*}
\int_{L_0}^{\ell_\smahom(t)}  \Gceff (l)  \, dl_{| U} 
	& = \sum_{\tau \in J_\smahom \cap [0,t]} \int_{\ell^-_\smahom (\tau)}^{\ell_\smahom (\tau)}  \Gceff (l)  \, dl =  \sum_{\tau \in J_\smahom \cap [0,t]}  D ( \ell_\hom (\tau) ) - D (\ell^-_\smahom (\tau) )  .
\end{align*}

In order to prove \eqref{e.eb5}, we will first check that  $dl_{|R}= d^{\,\sharp}_{LC} \ell_\hom$, 
where $d^{\,\sharp}_{LC} \ell_\hom$ denotes the push forward of the measure $d_{LC} \ell_\hom$; then, by the change of variable formula for the push-forward, see e.g.~\cite{AmbrosioFuscoPallara_00}, we get  
\begin{align*}
	\int_{L_0}^{\ell_\smahom(t)}  \Gceff (l)  \, dl_{|R} = \int_{L_0}^{\ell_\smahom(t)}  \Gceff (l)  \, d^{\,\sharp}_{LC} \ell_\hom = \int_{0}^{t}   \Gceff ( \ell_\hom(\tau))  \, d_{LC} \ell_\hom . 
\end{align*}

\separe

Note that both $dl_{|R}$ and  $d^\sharp_{LC} \ell_\hom$ are Borel measure which do not contain atoms. Thus, it is enough to check that they coincide on closed intervals $[l_\1, l_\2] \subset \mathrm{co} \, \ell_\hom ([0,T))  $.

\separe

To evaluate the push-forward measure of $[l_\1, l_\2]$ we should compute the $(d_{LC} \ell_\hom)$-measure of the set $ (\ell_\hom)^{-1} ([l_\1, l_\2])$. If $(\ell_\hom)^{-1} ( [l_\1, l_\2] ) = \emptyset$ the measure vanishes.
Otherwise, by monotonicity of $\ell_\hom$ it turns out that $(\ell_\hom)^{-1} ( [l_\1, l_\2] )$ is again an interval; more precisely, it is of the form $[t_\1, t_\2)$ or $[t_\1, t_\2]$, where $t_\1 = \min \{ (\ell_\hom)^{-1} ( [l_\1, l_\2]) \}$ and $t_\2 = \sup \{ (\ell_\hom)^{-1} ([l_\1, l_\2]) \}$. Note that the point $t_\1$ is always included by the right continuity of $\ell_\hom$, moreover,  the $d_{LC} \ell_\hom$-measure of $[t_\1, t_\2)$ or $[t_\1,t_\2]$ 
actually coincides with the $(d_{LC} \ell_\hom)$-measure of the open interval $(t_\1, t_\2)$ since $d_{LC} \ell_\hom = \dot\ell_\hom d\tau + d_{C} \ell_\hom$ does not contain atoms. Then
\begin{align}
	d_{LC} 
	\ell_\hom (t_\1, t_\2) 
		& = d \ell_\hom (t_\1, t_\2) - d_{J} \ell_\hom (t_\1,t_\2) \nonumber \\
		& = \ell^-_\hom (t_\2) - \ell_\hom (t_\1) -  \sum_{ \tau \,\in \, J_\smahom\cap \, (t_1, t_2)} \jump{\ell_\hom (\tau)}  \label{e.op} .
\end{align}

\separe

Let us turn to the measure $dl_{|R}$. If $(\ell_\hom)^{-1} ( [l_\1, l_\2] ) = \emptyset$ then $[l_\1, l_\2] \subset U$ and thus its $dl_{|R}$-measure vanishes.  Otherwise,  by \eqref{e.eb4} we can write 
$$
	dl_{|R} [ l_\1, l_\2 ] = ( l_\2 - l_\1)  - dl_{|U} [ l_\1, l_\2] = (l_\2 - l_\1) - \sum_{\tau \in J_\smahom} | U_\tau \cap [l_\1 , l_\2] |  .  
$$
Let us study the measure $| U_\tau \cap [l_\1 , l_\2] |$ as a function of $\tau \in J_\hom$. Note that $\tau >0$ since $\ell_\hom$ is right continuous. 


\separe

If $\tau < t_\1$ then $\ell_\hom  (\tau) < l_\1$ and thus $| U_\tau \cap [l_\1 , l_\2] |  = 0$.  If $ \tau > t_\2$ then 
$\ell^-_\hom ( \tau) \ge l_2$ (since $\ell_\hom (t) > l_2$ for $ t_\2 < t < \tau $) and  thus $| U_{\tau} \cap [l_\1 , l_\2] |  = 0$.

\separe

If $\tau = t_\1$ then 
$l_1 \le \ell_\hom (t_\1) \le l_\2$  while  $\ell^-_\hom (t_\1) \le l_1$   (since $\ell_\hom (t) < l_1$ for any $t< t_\1$). 
It follows that  
$| U_{t_1} \cap [l_\1 , l_\2] |  =  \ell_\hom (t_\1) - l_\1$. 
Hence, we can write 
$$
	 - l_\1 - | U_{t_1} \cap [l_\1 , l_\2] |  = - \ell_\hom (t_\1) . 
$$

\separe

If $\tau \in (t_\1, t_\2)$ then $l_1 \le \ell_\hom (t_\1) \le l_\2$ while $\ell^-_\hom (\tau)  \ge l_\1$ (since $\ell_\hom (t) \ge l_1$ for  $t_\1 \le t< \tau$). 
Hence $ [ \ell^-_\hom  (\tau) , \ell_\hom  (\tau) ] \subset [ l_1, l_2 ] $ and 
$| U_\tau \cap [l_\1 , l_\2] |   = \jump { \ell_\hom (\tau) }  $. 

\separe

If $\tau = t_\2$ and $t_2 \not\in   (\ell_\hom)^{-1} ( [l_\1, l_\2] )$ then $\ell_\hom (t_\2) > l_2$ while $l_1 \le \ell_\hom^- (t_2) \le l_2$ (since $l_1 \le \ell_\hom^- (t) \le l_2$ for every $t \in [t_1, t_\2)$). Then $| U_{t_2} \cap [l_\1 , l_\2] |  =  l_\2 - \ell^-_\hom (t_\2) $ and we can write
$$
	 l_\2 -  | U_{t_2} \cap [l_\1 , l_\2] |  = \ell^-_\hom (t_\2) . 
$$

\separe

If $\tau = t_\2$ and $t_2 \in   (\ell_\hom)^{-1} ( [l_\1, l_\2] )$ then $\ell_\hom (t_2) = l_2$; indeed, assume by contradiction that $\ell_\hom (t_2) < l_2 \le \ell_\hom(T)$ then $t_\2 < T$ and (by right continuity) $\ell_\hom (t) < l_2$ in a sufficiently small right neighborhood of $t_\2$. Since $l_1 \le \ell_\hom^- (t_2) \le l_2$ (as above) we get again
$$
	 l_\2 -  | U_{t_2} \cap [l_\1 , l_\2] |  = \ell^-_\hom (t_\2) . 
$$

\separe

\separe

In conclusion 
\begin{equation} \label{e.po}
	dl_{|R} [ l_\1, l_\2 ] = \ell^-_\hom (t_\2) -  \ell_\hom (t_\1) - \hspace{-12pt} \sum_{\tau \in J_\smahom \cap (t_1, t_2) } \jump{ \ell_\hom (\tau) }  . 
\end{equation}
By \eqref{e.po} and \eqref{e.op} the measures $dl_{|R}$ and $d^\sharp_{LC} \ell_\hom$ do coincide.  \qed

\subsubsection{Toughening and micro-instabilities} \label{tougher}

In this section we study the toughening effect resulting from our convergence result. As we will see toughening turns out to be the macroscopic effect of microscopic instabilities in the evolutions $\ell_n$. 

Let $ \Gchom = \lambda \GcA  + (1-\lambda) \GcB$ be the homogenized (or averaged) toughness, that is the weak\hspace{1pt}$^*$ limit of $ \Gc_n$. 

For simplicity, let us consider a time interval $[t_a, t_b]$ such that $\ell_\hom (t_a) < \ell_\hom (t_b)$ and  $J ( \ell_\hom ) \cap [t_a , t_b] = \emptyset$. We will show that 
$$  \Gchom \le  \Gceff ( l ) \ \text{for every $l\in [\ell_\hom (t_a) , \ell_\hom (t_b)]$}. 	$$ 
Let us rewrite the energy identity for $\ell_n$ in the form 
\begin{equation} 
	E_n ( t , \ell_n(t)) + \int_{L_0}^{\ell_n(t)}   \Gc_n (l) \, dl = \int_0^t \partial_t F_n  ( \tau , \ell_n (\tau)) \, d\tau  + d_{J} F_n ( [ 0, t] ) , 
\end{equation}
where we used the shorthand notation
$$
	d_{J} F_n ( [ 0, t] ) = \sum_{\tau \in J(\ell_n) \cap [0,t]} \llbracket F_n ( \tau , \ell_n (\tau)) \rrbracket  .
$$
Thus, for every $t_a \le t_\1 < t_\2 \le t_b$ we can write 
\begin{align}
	E_n ( t_\2 , \ell_n(t_\2))  -  
		  E_n ( t_\1 , \ell_n(t_\1))  & -  \int_{t_1}^{t_2} \partial_t F_n  ( \tau , \ell_n (\tau)) \, d\tau  \nonumber \\ & = - \int_{\ell_n (t_1)}^{\ell_n(t_2)}   \Gc_n (l) \, dl   +  d_{J} F_n ( (t_\1, t_\2] ) . 
\label{e.ebn}
\end{align}
Since $\ell_n \to \ell_\hom$ pointwise in $[0,T)$ by Proposition \ref{p.homen-vert} we get 
$$
	E_n ( t_i , \ell_n(t_i) )  \to  E_\hom ( t_i , \ell_\hom (t_i) ) \qquad \text{for $i=1,2$.} 
$$
Remember that $\partial_t F_n ( \tau , l ) = 2 f (\tau) \dot{f} (\tau) \, \E_n ( l)$, hence $\partial_t F_n ( \tau , \ell_n (\tau) ) \to \partial_t F ( \tau , \ell_\hom (\tau) ) $ and by dominated convergence we get 
$$
	 \int_{t_1}^{t_2} \partial_t F_n  ( \tau , \ell_n (\tau)) \, d\tau
	\to
	 \int_{t_1}^{t_2} \partial_t F_\hom   ( \tau , \ell_\hom (\tau)) \, d\tau .
$$
In summary, the left hand side of \eqref{e.ebn} converges to 
\begin{align*}
	E_\hom ( t_\2 , \ell_\hom (t_\2))  -  
		  E_\hom  ( t_\1 , \ell_n(t_\1))  & -  \int_{t_1}^{t_2} \partial_t F_\hom  ( \tau , \ell_n (\tau)) \, d\tau  \nonumber \\ & = - \int_{\ell_\smahom (t_1)}^{\ell_\smahom (t_2)}   \Gceff (l) \, dl   ,
\end{align*}
where the identity follows from Corollary \ref{c.enidhom} and from the assumption $J ( \ell_\hom ) \cap [t_a , t_b] = \emptyset$. 

\separe

As a consequence the right hand side of \eqref{e.ebn} converges, i.e., 
$$
	 \int_{\ell_n (t_1)}^{\ell_n(t_2)}   \Gc_n (l) \, dl   + d_{J} F_n ( (t_\1, t_\2] ) 
	\to 
   \int_{\ell_\smahom (t_1)}^{\ell_\smahom (t_2)}   \Gceff (l) \, dl   .
$$
As $ \Gc_n \weakstarto  \Gchom = \lambda \GcA  + (1-\lambda) \GcB$, 
$$
	\int_{\ell_n (t_\1)}^{\ell_n(t_\2)}   \Gc_n (l) \, dl  \to  \int_{\ell_\smahom (t_1)}^{\ell_\smahom (t_2)}   \Gchom (l) \, dl  .
$$
Hence
$$
	d_{J}  F_n ( (t_\1, t_\2] ) \to  \int_{\ell_\smahom (t_1)}^{\ell_\smahom (t_2)}  \Gchom (l) -   \Gceff (l) \, dl  .
$$
Since $d_{J} F_n ( (t_\1, t_\2] )  \le 0$ and  it follows that 
$$
     \int_{\ell_\smahom (t_1)}^{\ell_\smahom (t_2)}  \Gchom (l) -   \Gceff (l) \, dl  \le 0 . 
$$
By the arbitrariness of $t_\1$ and $t_\2$ it follows that $ \Gchom (l) \le  \Gc_\eff ( l )$  a.e.~in $[\ell_\hom (t_a) , \ell_\hom (t_a)]$. Since $\Gc_\hom$ is constant and $\Gc_\eff$ is upper semi-continuous the inequality actually holds everywhere in $[\ell_\hom (t_a) , \ell_\hom (t_b)]$.

\separe

\subsection{An explicit computation of the effective toughness \label{explicit}}

In this section we consider a specific case which highlights some peculiar features of the effective toughness, and in particular its dependence on the data. Besides the hypotheses of the previous section, we assume that $\mu_{\B,\1} \, \mu_{\B,\2} = \mu_{\A,\1} \, \mu_{\A,\2}$. No further assumptions are made on the toughness, thus it may occur $\GcA \neq \GcB$ or  $\GcA = \GcB$. 
These assumptions include for instance the case in which material $A$ is anisotropic and material $B$ is obtained by a $\frac{\pi}{2}$ rotation of material $A$. 

\separe

For convenience, let us recall the results of Proposition \ref{p.homen-vert} and a few notations. The homogenized stiffness matrix is given by 
$$
 \C_{\hom} =  \left(\begin{array}{cc}  \mu_{\hom, \1} & 0  \\ 0 & \mu_{\hom, \2}  \end{array}  \right) ,
$$ 
where $1 / \mu_{\hom, \1} = \lambda / \mu_{\A,\1} + (1-\lambda) / \mu_{\B,\1}$ and 
$\mu_{\hom, \2} =  \lambda \mu_{\A,\2} + ( 1- \lambda ) \mu_{\B,\2} $.
The (condensed) elastic energy $\E_\hom$ is of class $C^1 (0,L)$ and the energy release is denoted by $\G_\hom (l)$. 
We will prove that 
\begin{equation} \label{e.Gchom}
		\Gceff   =  \lambda  \max \left\{  \GcA  \, , \,  \GcB \, \frac{\mu_{\B,\1}}{\mu_{\A,\1}} \right\}   
		+ ( 1 - \lambda) \max \left\{  \GcA \frac{\mu_{\A,\1}}{\mu_{\B,\1}}  \, , \,  \GcB \right\}  .
\end{equation}
In particular, the effective toughness is constant and depends also on the elastic contrast $\mu_{\A,\1} / \mu_{\B,\1}$. 
More precisely, we will prove that $\Gc_\eff$ takes the above value where $\G_\hom \neq 0$, which is enough for our purpose in the light of Remark \ref{r.defGeff}.

\begin{example}  \label{e.bigG} Assume that $\mu_{\B,\1} < \mu_{\A,\1}$ and $\GcB \le \GcA$. Then 

$$
	 \Gceff = \lambda \GcA  +   ( 1 - \lambda) \GcA \frac{\mu_{\A,\1}}{\mu_{\B,\1}}  > \lambda \GcA  +   ( 1 - \lambda) \GcA = \GcA = \max \{ \GcA , \GcB \} .
$$

\end{example}

\separe

\noindent Next theorem states the main convergence result. 

\begin{theorem} \label{p.evhom} 
Let $\ell_n$ be the quasi-static evolutions given by Proposition \ref{p.existsn}. Then $\ell_n \weakstarto \ell_\hom$ in $BV(0,T)$, where $\ell_\hom$ is the unique non-decreasing, right continuous function which satisfies the initial condition $\ell_\hom (0) = L_0$ and Griffith's criterion in the following form:
\begin{itemize}
\item[i)] $G_\hom (t, \ell_\hom(t) ) \le \Gceff $ for every time $t \in [0,T)$;
\item[ii)]  if $G_\hom ( t, \ell_\hom (t) ) < \Gceff $ then $\ell_\hom$ is differentiable in $t$ and $\dot{\ell}(t) =0$; 
\item[iii)]  if $t \in J ( \ell_\hom)$ then $G_\hom (t , l ) \ge \Gceff $ for every $l \in [\ell^-_\hom(t) , \ell_\hom(t) )$.
\end{itemize}
\end{theorem}

Note that in the above statement it is not necessary to extract any subsequence of $\ell_n$. As the effective toughness is constant, the effective dissipated energy boils down to 
$D_\eff (l)  =    \Gceff ( l - L_0)$;
by Corollary \ref{c.existsn}  we get the energy identity for $F (t, u) = E_\hom (t, u) + D_\eff (l)$, stated hereafter.

\begin{corollary} \label{c.existsnhom} For every $ t \in [0,T)$ the following energy identity holds, 
\begin{equation} 
	F ( t , \ell_\hom  (t)) = \int_0^t \partial_t F  ( \tau , \ell_\hom (\tau)) \, d\tau  +  \hspace{-12pt} \sum_{\tau \, \in \, J_\smahom \cap [0,t]} \llbracket F ( \tau , \ell_\hom (\tau)) \rrbracket .
\end{equation}
\end{corollary}

\begin{remark} \normalfont Note that the limit energy $F (t, \cdot) $ does not coincide with the $\Gamma$-limit of the energies $F_n (t,\cdot) = E_n (t, \cdot) + D_n (\cdot)$, which is given by 
$$
F_\hom (t,l) = E_\hom (t, l) +  \Gchom (l - L_0) \ \text{ where } \  \Gchom = \lambda \GcA + (1-\lambda) \GcB .
$$ 
\end{remark} 

Before proceeding we state, for the reader's convenience, the following elementary change of variable, which will be often used in the sequel. 

\begin{lemma} \label{l.gigetto} Let $O$  be an open set in $\R^2$ and ${\boldsymbol C} \in \R^{2 \times 2}$. Let $\Phi : O \to \widehat O$ be a bi-Lipschitz map. For $z \in H^1(O)$ denote $\hat{z} = ( z \circ \Phi^{-1}) \in H^1(\widehat O)$, then 
\begin{equation} \label{e.chvar}
	\int_O \nabla z \, {\boldsymbol C} \, \nabla z^T dx = \int_{\widehat O} \nabla \hat{z} \,  \widehat{\boldsymbol C} \nabla \hat{z}^T d\hat{x}
	\quad \text{where} \quad  \widehat{\boldsymbol C} =   (D\Phi ) {\boldsymbol C} (D\Phi )^{T} ( \mathrm{det}D\Phi) ^{-1} .
\end{equation}
In particular, if $\Phi ( x_\1 , x_\2 ) = ( \alpha x_\1 , \beta x_\2)$ and $\C = \mu_\1 \hat{e}_\1 \otimes \hat{e}_\1 + \mu_\2 \hat{e}_\2 \otimes \hat{e}_\2$ then 
$$
	\widehat{\boldsymbol C} = \left(   \begin{array}{cc}  \mu_{\1} \alpha / \beta   &  0 \\  0 &  \mu_{\2} \beta / \alpha   \end{array} \right) . 
$$
\end{lemma} 

\separe

\noindent {\bf Proof of Theorem \ref{p.evhom}.} The Theorem follows from 
\begin{align} \label{e.Gchom-hor}
\frac{\G_\hom (l) }{\displaystyle \Gamma\text{-}\hspace{-4pt}\lim_{n \to \infty} \frac{\G_{n}^{\phantom{I}}}{\Gc_{n}} (l) }  = \lambda  \max \left\{  \GcA  \, , \,  \GcB \, \frac{\mu_{\B,\1}}{\mu_{\A,\1}} \right\}   
		+ ( 1 - \lambda) \max \left\{  \GcA \frac{\mu_{\A,\1}}{\mu_{\B,\1}}  \, , \,  \GcB \right\} ,
\end{align}
for $\G_\hom (l) \neq 0$, after employing  Remark \ref{r.sub-seq} and Remark \ref{r.defGeff}.

{\bf I.} Let $\alpha = \mu_{\,\B,\1} / \mu_{\,\A,\1}$. For $S =  \lambda L \alpha + (1-\lambda )L$ let $\widehat\Omega = (0, S)  \times  (-H,H)$. 
We consider a bi-Lipschitz piecewise affine map $\Phi_n : \Omega \to \widehat\Omega$ of the form $\Phi_n ( x_\1, x_2 ) = (\phi_n ( x_\1) , x_2 )$ where 
$$	\phi_n (x_1) = \int_0^{x_1} \alpha_n \, ds \quad \text{and } 
\quad
	\alpha_n  = \begin{cases} \alpha   & \text{in $K_{n,\A}$} \\ 1  & \text{in $K_{n,\B}$.}  \end{cases}  
$$ 

We employ the notation $\widehat\Omega_{n,\A}$ to denote the set $\Phi_n ( \Omega_{n,\A} )$ and similarly for all the sets introduced in \S\ref{s.2-1}. We denote $\hat{g}_n = g \circ \Phi_n^{-1}$ and consider the spaces
$$	\widehat\U_{s} = \{ \hat{u} \in H^1(\widehat\Omega \setminus K_s) : \hat{u}  = \hat{g}_n \text{ in } \partial_D \widehat\Omega \} 
\quad \text{for $s \in (0,S]$.  } 
$$ 
Let $s = \phi_n(l)$. For $u \in \U_{\l}$ and $ \hat{u} = u \circ \Phi_n^{-1} \in \widehat\U_{s}$ we have 
and 
$$   \W_{l,n} (u ) 
 = \widehat \W_{s,n} ( \hat{u} ) = \tfrac12 \int_{\widehat\Omega \setminus \widehat{K}_{s} } \nabla\hat{u} \, \widehat{\boldsymbol C}_n \nabla\hat{u}^T \, dy . $$ 
Under the assumption $\mu_{\B,\1} \, \mu_{\B,\2} = \mu_{\A,\1} \, \mu_{\A,\2}$ and with the above choice of $\alpha = \mu_{\,\B,\1} / \mu_{\,\A,\1}$, it turns out that $\widehat{\boldsymbol C}_n$ is constant in $\widehat\Omega$; indeed we have 
	$$
	\widehat{\boldsymbol C}_n =  \left(   \begin{array}{cc}   \alpha \, \mu_{\A,\1}  &  0 \\  0 & \mu_{\A,\2}  / \alpha  \end{array} \right)  =  \left(   \begin{array}{cc}   \mu_{\B,\1}  &  0 \\  0 & \mu_{\B,\2}  \end{array} \right) 
	\ \text{in $\widehat\Omega_{A,n}$} ,  \qquad 
	\widehat{\boldsymbol C}_n = \left(   \begin{array}{cc}   \mu_{\B,\1}  &  0 \\  0 & \mu_{\B,\2}  \end{array} \right) 
	\ \text{in $\widehat\Omega_{B,n}$.}  
	$$
Note also that $\widehat{\boldsymbol C}_n$ is actually independent of $n$; therefore in the sequel we will skip the dependence on $n$ in the notation. The reduced energy $\widehat\E$ is then of class $C^1 (0,S)$; in particular the energy release $\widehat\G$ is continuous in $(0,S)$.
Clearly, 
by the change of variable 
\begin{gather}  \label{e.su_n}
 \widehat{\E} (\phi_n (l) ) = \E_n (l)  \quad \text{and thus} \quad 
\widehat{\G} (\phi_n (l)) \,\phi_n' (l) = \G_n ( l ) \  \text{for $ l \in (0,L) \setminus \Lambda_n$.}
\end{gather}

{\bf II.} Let $\alpha_{\hom} =  \alpha \lambda + ( 1- \lambda)$.  Let $\Phi : \Omega \to \widehat\Omega$ be the linear map $\Phi ( x_\1, x_2 ) = ( \phi (x_\1) , x_2 )$ for $\phi(x_1) = \alpha_\hom x_1$.  With this change of variable, for $s = \phi(l) = \alpha_\hom l $, we write 
$$
     \W_{l,\hom} (u)  =  \widehat{\W}_{s,\hom} ( \hat{u} )  = 
     \tfrac12 \int_{\widehat\Omega \setminus \widehat K_s} \nabla \hat{u}  \, \widehat{\mathbf{C}}_\hom \, \nabla \hat{u}^T \, dy  ,  
$$
where $$\widehat{\boldsymbol C}_\hom = 
\left( \begin{array}{cc}  \alpha_{\hom} \, \mu_{\hom,\1}  & 0  \\  0   & \mu_{\hom, \2} /  \alpha_{\hom} \end{array} \right)  . $$
As expected, $\widehat{\mathbf{C}}_\hom =  \widehat{\mathbf{C}}$; indeed, 
$$
	 \frac{1}{\mu_{\hom,\1}} = \frac{\lambda}{\mu_{\A,\1}} + \frac{1-\lambda}{\mu_{\B,\1}} = 
	 \frac{\lambda \mu_{\B,\1} + (1-\lambda) \mu_{\A,\1}}{\mu_{\A,\1}\mu_{\B,\1}} ,
	  \qquad
	  \alpha_\hom = \frac{\mu_{\B,\1}}{\mu_{\A,\1}} \lambda + ( 1- \lambda) = \frac{\lambda \mu_{\B,\1} + (1-\lambda) \mu_{\A,\1}}{\mu_{\A,\1}}  ,
$$
and thus $ \alpha_\hom \, \mu_{\hom,\1} = \mu_{\B, \1}$, moreover, recalling that $\mu_{\B,\1} \, \mu_{\B,\2} = \mu_{\A,\1} \, \mu_{\A,\2}$, 
$$
	\mu_{\hom, \2} =  \lambda \mu_{\A,\2} + ( 1- \lambda ) \mu_{\B,\2} ,
	\qquad
	\alpha_\hom = \frac{\mu_{\A,\2}}{\mu_{\B,\2}} \lambda + ( 1- \lambda) = \frac{\lambda \mu_{\A,\2} + ( 1- \lambda ) \mu_{\B,\2}}{\mu_{\B,\2}} ,
$$
and then $ \mu_{\hom,\2} / \alpha_\hom  = \mu_{\B, \2}$. In conclusion, 
\begin{equation} \label{e.su-hom} 
\widehat\E ( \phi (l) ) = \E_\hom (l) \quad \text{and thus} \quad  \widehat\G ( \phi ( l ) ) \alpha_\hom = \G_\hom (l)  
\  \text{for $ l \in (0,L)$.} 
\end{equation}

{\bf III.} By \eqref{e.su_n} we can write 
$$
	\frac{\G_{n} (l) }{\Gc_{n} (l) }  = 
	\frac{\widehat\G (\phi_n (l) ) \, \phi'_n (l) }{\Gc_{n}(l)} . 
$$
Note that $\phi_n \to \phi$ uniformly in $[0,L]$. Then, being $\widehat\G$ continuous, $\widehat\G \circ\phi_n \to \widehat\G \circ \phi$ locally uniformly in $(0,L)$. As a consequence
\begin{equation} \label{e.gofam}
	\Gamma\text{-}\hspace{-4pt}\lim_{n \to \infty} 
	\frac{\G_{n}}{\Gc_{n}} (l)  = 
	\widehat\G ( \phi (l) )  \, \Gamma\text{-}\hspace{-4pt}\lim_{n \to \infty} 
	\frac{\phi'_n }{\Gc_{n}} (l) . 
\end{equation} 
Writing explicitly
$$
	\frac{\, \phi'_n (l) }{\Gc_{n}(l)}  = \begin{cases}
	\alpha / \GcA    & \text{ in $K_{n,\A}$, } \\
	1 / \GcB & \text{ in $K_{n,\B}$, }
	\end{cases}
$$
we get $\Gamma\text{-}\lim_{n \to \infty}  {\displaystyle \frac{\phi'_n}{\Gc_{n}} } (l) = \min \{   \alpha / \GcA  , 1 / \GcB \}$.  If $\G_\hom (l) \neq 0$ then $\widehat\G ( \phi (l) ) \neq 0$ and by \eqref{e.su-hom}  and \eqref{e.gofam} 
\begin{align*} 
\Gc_\eff (l) & = \frac{\G_\hom (l) }{\displaystyle \Gamma\text{-}\hspace{-4pt}\lim_{n \to \infty} \frac{\G_{n}^{\phantom{I}}}{\Gc_{n}} (l) }  = 
\frac{ \alpha_\hom \widehat\G ( \phi (l) )  } {  \widehat\G ( \phi (l) )  \min \{   \alpha / \GcA  , 1 / \GcB \} } =  \big( \alpha \lambda + ( 1- \lambda) \big) \max \{  \GcA / \alpha  , \GcB \}     \\
& =  \lambda  \max \left\{  \GcA  \, , \,  \GcB \, \frac{\mu_{\B,\1}}{\mu_{\A,\1}} \right\}   
		+ ( 1 - \lambda) \max \left\{  \GcA \frac{\mu_{\A,\1}}{\mu_{\B,\1}}  \, , \,  \GcB \right\} .
\end{align*}
The proof of \eqref{e.Gchom-hor}  is concluded.  \qed 

\begin{remark} \normalfont In general $\G_\hom$ is neither the pointwise nor the $\Gamma$-limit of $\G_n$. Indeed, by \eqref{e.su_n} and \eqref{e.su-hom} we have 
$$
	\widehat{\G} (\phi_n (l)) \,\phi_n' (l) = \G_n ( l ) 
	\ \text{ and  } \ 
	\widehat\G ( \phi ( l ) ) \alpha_\hom = \G_\hom (l)  .
$$
We know that $\widehat{\G} \circ \phi_n \to \widehat\G \circ \phi$ locally uniformly but the pointwise limit of $\phi'_n$ does not exists. Moreover, $\Gamma\text{-}\lim_{n \to \infty} \phi'_n (l) = \min \{ \alpha , 1 \}$, while $\alpha_\hom = \lambda \alpha + ( 1- \lambda) = \min \{ \alpha , 1 \}$ if and only if $\alpha =1$; the latter case is trivial, since $\mu_{\B, \1} = \mu_{\A,\1}$ and $\mu_{\B,\1} \, \mu_{\B,\2} = \mu_{\A,\1} \, \mu_{\A,\2}$  imply that the stiffness matrix ${\boldsymbol C}_n $ is constant. 
\end{remark}

\section{Horizontal layers with horizontal crack} \label{s.oro}

\subsection{Energy and energy release}

As in the previous section we assume that the uncracked reference configuration is the set $\Omega = (0,L) \times (-H,H)$ and that the crack set is of the form $K_l = [0, l] \times \{ 0 \}$ for $l \in (0,L]$. On the contrary, here we assume that $\Omega$ is made periodic horizontal layers of thickness $h_n = H/n$ for $n \in 2 \mathbb{N}$ (the fact that $n$ is even implies that $K_l$ is contained in the interface between two layers). For $\lambda \in (0,1)$ each layer is itself made of a horizontal layer of material $A$, with thickness $\lambda h_n $, and an horizontal layer of material $B$, with thickness $(1-\lambda) h_n$.

For functional spaces, boundary condition, energies and energy release we employ the same assumptions and notation of the previous section, apart from toughness which here is simply denoted by $\Gc$, since it is independent of $n$.
Following the proof of Lemma \ref{l.En} it not difficult to show the following result. 

\begin{lemma} \label{l.Enhoriz} $\E_n : [0,L] \to \R$ is decreasing and continuous.
Moreover, it is of class $C^1 (0,L)$ and
\begin{equation} \label{e.E'=DEhor}
      \E'_n ( l) = \int_{\Omega \setminus K_l} 
      \nabla u_{l,n} \boldsymbol{C}_n \boldsymbol{E} \nabla u_{l,n}^{T} \, \phi' dx 
	 \quad \text{for} \quad 
    \boldsymbol{E} =   \left( \begin{matrix} -1 &   0 \\ 0 & 1 \end{matrix} \right) , 
\end{equation}
where $\phi \in W^{1,\infty} (0,L)$ with $\mathrm{supp} ( \phi)  \subset (0,L)$, $0 \le \phi \le 1$ and $ \phi(l) =1$.
\end{lemma}

\begin{remark} Comparing with Lemma \ref{l.En} note that here $\phi$ is independent of $n$, in particular its support does not shrink with $n$. 
Finally, note that $\E_n$ is of class $C^1$ in the whole $(0,L)$. \end{remark}

\subsection{Quasi-static evolution by Griffith's criterion}

Following the proofs of Proposition \ref{p.existsn} and Corollary \ref{c.existsn} we obtain the following results.

\begin{proposition} \label{p.existsnhor} There exists a unique non-decreasing, right continuous function $\ell_n : [0,T) \to [L_0, L]$ which satisfies the initial condition $\ell_n(0) = L_0$ and Griffith's criterion, in the following form:
\begin{itemize}
\item[i)] $G_n (t, \ell_n(t) ) \le  \Gc$ for every time $t \in [0,T)$;
\item[ii)] if $G_n ( t, \ell_n (t) ) < \Gc$ then $\ell_n$ 
is differentiable in $t$ and $\dot \ell_n  (t) = 0$;
\item[iii)] if $t \in J ( \ell_n)$ then $G_n (t , l ) \ge \Gc $ for every $l \in [\ell^-_n(t) , \ell_n(t) )$.
\end{itemize}
\end{proposition} 

\begin{corollary} \label{c.betterhor} 
If $G_n ( t, \ell_n (t) ) < \Gc $ then $\ell_n$ is constant in a neighborhood of $t$. 
\end{corollary} 

\begin{corollary} \label{c.existsnhor} For every $ t \in [0,T)$ the following energy identity holds, 
\begin{equation}
	F_n ( t , \ell_n (t)) = \int_0^t \partial_t F_n  ( \tau , \ell_n (\tau)) \, d\tau  + \sum_{\tau \in J(\ell_n) \cap [0,t]} \llbracket F_n ( \tau , \ell_n (\tau)) \rrbracket  .
\end{equation}
\end{corollary}

\subsection{Homogenization of Griffith's criterion}

Here the homogenized stiffness matrix is 
$$
\C_{\hom} =  \left(\begin{matrix}  \mu_{\hom,\1} & 0  \\ 0 & \mu_{\hom,\2}  \end{matrix}  \right) ,
$$ 
where $\mu_{\hom,\1} =  \lambda \mu_{\A,\1} + ( 1- \lambda ) \mu_{\B,\1} $
 is the weak* limit of $\mu_{n,\1}$ while   $1 / \mu_{\hom,\2} = \lambda / \mu_{\A,\2} + (1-\lambda) / \mu_{\B,\2}$ is the weak* limit of $1/\mu_{n,\2}$. 
 
Accordingly, we introduce the homogenized elastic energy $\W_{l,\hom}$ 
and the reduced homogenized energy $\E_\hom$. 
Arguing as in Lemma \ref{l.En},  it turns out that the energy $\E_\hom$ is of class $C^{1} (0,L)$, hence the energy release $\G_\hom (l) = - \partial_l \E_\hom (l)$ is well defined in $(0,L)$. 
We set $\G_\hom (L)=0$.

\separe

In this setting it turns out that $ \Gceff =  \Gc$; indeed, we have the following result. 

\begin{theorem} \label{t.evhomhor} Let $\ell_n$ be the quasi-static evolutions given by Proposition \ref{p.existsnhor}. There exists a subsequence (not relabeled) such that $\ell_n \to \ell$ pointwise in $[0,T)$. Then $\ell^+ = \ell_\hom$, where $\ell_\hom$ is  the unique non-decreasing, right continuous function which satisfies the initial condition $\ell_\hom (0) = L_0$ and Griffith's criterion: 
\begin{itemize}
\item[i)] $G_\hom (t, \ell_\hom(t) ) \le \Gc$ for every time $t \in [0,T)$;
\item[ii)] if $G_\hom ( t, \ell_\hom (t) ) < \Gc$ then $\dot \ell_\hom  (t) = 0$;
\item[iii)] if $t \in J ( \ell_\hom)$ then $G_\hom (t , l ) \ge \Gc$ for every $l \in [\ell^-_\hom(t) , \ell_\hom(t) )$.
\end{itemize}
\end{theorem}

\begin{corollary} \label{c.enidhom2} For every $ t \in [0,T)$ the following energy identity holds, 
\begin{equation} 
	F_\hom ( t , \ell_\hom (t)) = \int_0^t \partial_t F_\hom  ( \tau , \ell_\hom (\tau)) \, d\tau + \sum_{\tau \in J(\ell_\smahom) \cap [0,t]} \llbracket F_\hom ( \tau , \ell_\hom (\tau)) \rrbracket  .
\end{equation}
\end{corollary}

\separe

The proof of Theorem \ref{t.evhomhor} is a consequence of Corollary \ref{c.elln2ell-a} and Remark  \ref{r.defGeff}, together with the following result on the convergence of the energy release.

\begin{theorem} \label{t.Gconv} $\G_n$ converge to $\G_\hom$ locally uniformly in $(0,L)$. 
\end{theorem}

\noindent The proof of this theorem is contained in the next subsection. 

\medskip


\subsection{Convergence of the energy release} \label{ss.conv}

\separe

\begin{lemma} \label{l.lerr} $\G_n$ is uniformly bounded in $C^{\,0, 1/2}( [ L',L'' ])$ for every $0 < L' < L'' < L$. \end{lemma}

\proof By minimality $\W_{l,n} ( u_{l,n} ) \le \W_{l,n} (g) \le C$ where $C$ is independent of $n$. Hence, by uniform coercivity, $ \| \nabla  u_{l,n} \|_{L^2} \le C $ for every $l \in [0,L]$. 

\separe

Now, following \cite{N_AMO17} we prove that there exists $C>0$ such that for every {\normalfont $0 \le l_\1 < l_\2 \le L$} and every $n \in 2 \mathbb{N}$ it holds 
\begin{equation} \label{e.AMO}
\| \nabla   u_{l_1,n} - \nabla  u_{l_2,n} \|_{L^2} \le C | \E_n ( l_\1 ) -  \E_n (l_\2) | ^{1/2}  .
\end{equation}
By the variational formulation we have 
\begin{gather*}
     \int_{ \Omega \setminus K_{l_1}}  \nabla  {u}_{l_1,n} \C_n \nabla   v_{\1}^T  \, dy  = 0  \qquad \text{for every $  v_{\1} \in  \V_{\l_1}$,} \\
     \int_{ \Omega \setminus K_{l_2}}  \nabla  {u}_{l_2,n} \C_n \nabla   v_{\2} ^T \, dy  = 0  \qquad \text{for every $  v_{\2} \in  \V_{\l_2}$.} 
\end{gather*}
Since $( {u}_{l_1,n} -   {u}_{l_2,n} ) \in  \V_{\l_2}$ we get 
$$
     \int_{ \Omega \setminus K_{l_2}}  \nabla  {u}_{l_2,n} \C_n \nabla (  {u}_{l_1,n} -   {u}_{l_2,n} )^T \, dy  = 0  . \qquad 
$$
Then, by monotonicity and uniform coercivity of the energy we get 
\begin{align*}
	|   \E_n ( l_\1 ) -   \E_n (l_\2) | & =   \E_{n} ( l_\1 ) -   \E_n (l_\2) 
	\\ & 
	= \tfrac12 \int_{ \Omega } \nabla  {u}_{l_1,n} \C_n \nabla  {u}_{l_1,n}^T  - \nabla  {u}_{l_2,n} \C_n \nabla  {u}_{l_2,n}^T \, dy  \\ & 
	= \tfrac12 \int_{ \Omega } \nabla (  {u}_{l_1,n} +  {u}_{l_2,n}) \C_n \nabla (  {u}_{l_1,n} -   {u}_{l_2,n})^T  \, dy  \\ &
	= \tfrac12 \int_{ \Omega } \nabla   {u}_{l_1,n}  \C_n \nabla (  {u}_{l_1,n} -   {u}_{l_2,n})^T  \, dy  \\ &
	= \tfrac12 \int_{ \Omega } \nabla (  {u}_{l_1,n} -  {u}_{l_2,n}) \C_n \nabla (  {u}_{l_1,n} -   {u}_{l_2,n})^T  \, dy
\ge C \| \nabla   {u}_{l_1,n} -  \nabla  {u}_{l_2,n} \|_{L^2}^2 . 
\end{align*}

\separe

Next, we show that given $0< L' < L''< L$ there exists $C >0$ s.t. 
$$ \G_n (l) \le C \quad \text{for every $l \in [L',L'']$ and every $n \in 2 \mathbb{N}$.}  $$ 
 Let $\phi \in C_c^{\infty} (0,L)$ with $0 \le \phi \le 1$ and $\phi=1$ in $[L',L'']$ such that 
$$
     \E'_n (l) = \int_{  \Omega} 
      \nabla   u_{l,n} \C_n \boldsymbol{E} \, \nabla  {u}_{l,n}^T \,  \phi' \, dy
	 \quad \text{for} \quad 
    \boldsymbol{E} =   \left( \begin{matrix} -1 &  0 \\ 0 & 1 \end{matrix} \right)  .
$$
Hence $ \G_n (l) \le C $ (note that $C$ depends on $\phi'$ and thus on $L'$ and $L''$). It follows that the energies $ \E_n$ are uniformly Lipischitz continuous in $[L',L'']$ and thus, by \eqref{e.AMO}
$$
	\| \nabla   u_{l_1,n} - \nabla   u_{l_2,n} \|_{L^2} \le C | l_1 -  l_2  | ^{1/2} . 
$$
As a consequence, using the integral representation of $\G_n = -   \E'_n$ we get 
$$
	| \G_n ( l_\1 ) - \G_n ( l_\2 ) | \le C | l_\1 - l_\2 |^{1/2} , 
$$
which concludes the proof.  \qed

By Ascoli-Arzel\`a Theorem, with the aid of a diagonal argument, Lemma \ref{l.lerr} implies the following result.

\begin{corollary} \label{c.Gconv} There exists a subsequence (non relabeled) and a limit $\G$ such that $\G_n \to \G$ locally uniformly in $(0, L)$.
\end{corollary}

We claim that $\G =\G_\hom$. To this end, we will need a couple of lemmas on the local convergence of the energy.

\begin{lemma} \label{l.locenconv} For $0 < a < b < L$, let $\Omega_{(a,b)} = (a,b) \times (- H , H)$. Then 
$$
	\int_{\Omega_{(a,b)} \setminus K_l } \nabla  {u}_{l,n} \C_n \nabla  {u}_{l,n}^T  \,dx \ \to \ \int_{\Omega_{(a,b)} \setminus K_l} 
	\nabla  {u}_{l,\hom} \C_{\hom} \nabla  {u}_{l,\hom}^T   \,dx .
$$
\end{lemma}

\proof  For $\eps >0$ let $\theta_\eps \in C^\infty_c (0,L)$ with $0 \le \theta_\eps \le 1$, $\theta_\eps =1$ in $(a,b)$, and such that 
$$
	\int_{ \Omega \setminus K_l} \nabla  {u}_{l,\hom} \C_\hom \nabla  {u}_{l,\hom}^T \theta_\eps \, dx \le \eps + \int_{ \Omega_{(a,b)} \setminus K_l} 
	\nabla  {u}_{l,\hom} \C_\hom \nabla  {u}_{l,\hom}^T \, dx . 
$$
Then, by the properties of $\theta_\eps$ and since $  {u}_{l,n} \theta_\eps  \in  \V_l$ by the variational formulation we can write
\begin{align*}
	\int_{ \Omega_{(a,b)} \setminus K_l} \nabla  {u}_{l,n} \C_n \nabla  {u}_{l,n}^T \, dx  & \le \int_{ \Omega \setminus K_l} \nabla  {u}_{l,n} \C_n \nabla  {u}_{l,n}^T \theta_\eps \, dx \\
		& = \int_{ \Omega \setminus K_l}\nabla  {u}_{l,n} \C_n  \left( \nabla (u_{l,n} \theta_\eps) -  {u}_{l,n} \nabla \theta_\eps  \right)^T dx \\
			& = - \int_{ \Omega \setminus K_l}\nabla  {u}_{l,n} \C_n  \nabla \theta_\eps^T  {u}_{l,n}  \, dx .
\end{align*}
As $\C_n \nabla  {u}_{l,n}^T \weakto \C_\hom \nabla  {u}_{l,\hom}^T$ in $L^2 (  \Omega \setminus K_l; \R^2)$ and since $ {u}_{l,n} \to  {u}_{l,\hom}$ in $L^2( \Omega \setminus K_l)$ the last term converge to 
\begin{align*}
	 - \int_{ \Omega \setminus K_l}\nabla  {u}_{l,\hom}  \C_\hom \nabla \theta_\eps^T  {u}_{l,\hom} \, dx & = 
	  - \int_{ \Omega \setminus K_l}\nabla  {u}_{l,\hom} \C_\hom  \left( \nabla ( {u}_{l,\hom} \theta_\eps) - \nabla  {u}_{l,\hom} \theta_\eps  \right)^T dx \\
	  & = \int_{ \Omega \setminus K_l} \nabla  {u}_{l,\hom} \C_\hom \nabla  {u}_{l,\hom}^T \theta_\eps \, dx  \\ & \le \eps +  \int_{ \Omega_{(a,b)} \setminus K_l} \nabla  {u}_{l,\hom} \C_\hom \nabla  {u}_{l,\hom}^T \, dx .
\end{align*}
Hence
$$
	\limsup_{n \to \infty} \int_{ \Omega_{(a,b)} \setminus K_l} \nabla  {u}_{l,n} \boldsymbol{C}_n \nabla  {u}_{l,n}^T  \,dx \le \eps + \int_{ \Omega_{(a,b)} \setminus K_l} \nabla  {u}_{l,\hom} \C_\hom  \nabla  {u}_{l,\hom}^T \, dx  . 
$$
\separe

For $\eps>0$ let $\psi_\eps \in C^\infty_c ( 0, L )$ with $ 0 \le \psi_\eps \le 1$, $\psi_\eps = 0$ in $(0,L) \setminus (a,b) $,  and such that 
$$
	\int_{ \Omega \setminus K_l}\nabla  {u}_{l,\hom} \C_\hom \nabla  {u}_{l,\hom}^T \psi_\eps \, dx \ge - \eps +  \int_{ \Omega_{(a,b)} \setminus K_l} 
	\nabla  {u}_{l,\hom} \C_\hom  \nabla  {u}_{l,\hom}^T \, dx . 
$$
Then, arguing as above, 
\begin{align*}
	\int_{ \Omega_{(a,b)} \setminus K_l} \nabla  {u}_{l,n}  {\boldsymbol{C}}_n \nabla  {u}_{l,n}^T  \,dx & \ge \int_{ \Omega \setminus K_l}  \nabla  {u}_{l,n}  {\boldsymbol{C}}_n \nabla  {u}_{l,n}^T \psi_\eps \, dx  \\ 
		& \hspace{65pt} \downarrow \\ 
			& \phantom{=}  \int_{ \Omega \setminus K_l}  \hspace{-5pt} \nabla  {u}_{l,\hom} \C_\hom  \nabla  {u}_{l,\hom}^T \psi_\eps \, dx \ge - \eps + \int_{ \Omega_{(a,b)} \setminus K_l}  \nabla  {u}_{l,\hom} \C_\hom  \nabla  {u}_{l,\hom}^T \,  dx . 
\end{align*}
Hence
$$
	\liminf_{n \to \infty} \int_{ \Omega_{(a,b)} \setminus K_l} \nabla  {u}_{l,n}  {\boldsymbol{C}}_n \nabla  {u}_{l,n}^T  \,dx \ge - \eps + \int_{ \Omega_{(a,b)} \setminus K_l} \nabla  {u}_{l,\hom} \C_\hom  \nabla  {u}_{l,\hom}^T \, dx  . 
$$
We conclude by the arbitrariness of $\eps >0$. \qed 

\begin{lemma} \label{l.conv-comp} For $0 < a < b < L$, let $\Omega_{(a,b)} = (a,b) \times (- H , H)$. Then
\begin{equation} \label{e.zio}
	\int_{\Omega_{(a,b)} \setminus K_l} \mu_{n,\1} \, | \partial_{x_1} u_{l,n} |^2  \,dx
	\ \to \
	\int_{\Omega_{(a,b)} \setminus K_l} \mu_{\hom,\1} \, | \partial_{x_1} u_{l,\hom} |^2 \,dx . 
\end{equation}
\end{lemma}

\proof We employ a change of variable. Let $\widehat\Omega = (0, L)  \times (- S  , S )$ where 
$ S =  H \left(  \lambda / \mu_{\B,\1}  +  (1-\lambda )  /  \mu_{\A,\1} \right)$. 
Denote by $\Pi_{n,\A}$ the projection of $\Omega_{n,\A}$ on the verical axis, i.e., 
$$
	\Pi_{n,\hspace{1pt}\A} =  \{ x_2 \in (-H,H) : (x_1, x_2) \in \Omega_{n,\hspace{1pt}\A}  \text{ for $x_1 \in (0,L)$}   \}  ,
$$
and similarly for $\Pi_{n,\hspace{1pt}\B}$. Let $\Phi_n : \Omega \to \widehat\Omega$ be the bi-Lipschitz piecewise affine map $\Phi_n ( x_\1, x_2 ) = ( x_\1 , \phi_n (x_2) )$ where 
$$	\phi_n (x_2) = \int_0^{x_2} \beta_n \, ds \quad \text{and } 
\quad
	\beta_n  = \begin{cases} 1 / \mu_{\B,\1}   & \text{in $\Pi_{n,\hspace{1pt}\A}$, } \\ 1/ \mu_{\A,\1} & \text{in $\Pi_{n,\hspace{1pt}\B}$.}  \end{cases}  
$$ 
Denote $\hat{g}_n = g \circ \Phi_n^{-1}$ and 
$$    \widehat\U_{\l,n} = \{ \hat{u} \in H^1 (\widehat\Omega \setminus K_l)  : \hat{u} = \hat{g}_n \text{ in $ \partial_D \widehat\Omega$} \} \qquad 
\widehat\V_l = \{ \hat{v} \in H^1 (\widehat\Omega \setminus K_l)  : \hat{v} = 0 \text{ in $ \partial_D \widehat\Omega$} \} . $$ 
For $u \in \U_{\l}$  and $ \hat u= ( u \circ \Phi_n^{-1} ) \in \widehat\U_{\l,n}$ we have 
$$   \W_{l,n}  (u ) 
= \widehat{\W}_{l,n} ( \hat u ) =  \tfrac12 \int_{\widehat\Omega \setminus K_{l} } \nabla \hat u  \, \widehat{\boldsymbol C}_n \nabla \hat u  ^T \, dy  $$
where (see Lemma \ref{l.gigetto}) 
	$$
	\widehat{\boldsymbol C}_n = \left(   \begin{matrix}   \hat\mu_{n,\1}  &  0 \\  0 &  \hat\mu_{n,\2}   \end{matrix} \right)  
	\quad  \text{with}  \quad 
	\hat\mu_{n,\1} = \mu_{\A,\1} \, \mu_{\B,\1} 
	\quad \text{and} \quad
	\hat\mu_{n,\2}    =  \begin{cases}  \mu_{\A,\2} / \mu_{\B,\1}  & \text{in $\widehat\Omega_{n, \hspace{1pt}\A}$,}   \\   \mu_{\B,\2} / \mu_{\A,\1}  & \text{in $\widehat\Omega_{n,\hspace{1pt}\B}$.}   \end{cases}
	$$
Note that $\hat{\mu}_{n, \1}$ is constant and independent of $n \in 2 \mathbb{N}$. We denote $\hat{\mu}_{\A,\2} = \mu_{\A,\2} / \mu_{\B,\1} $ and $\hat{\mu}_{\B,2} = \mu_{\B,\2} / \mu_{\A,\1} $.
Note that the relative size of the layers $\widehat\Omega_{n,\hspace{1pt}\A}$ and $\widehat\Omega_{n,\hspace{1pt}\B}$ is respectively
$$
	\hat \lambda 	= \frac { \lambda  / \mu_{\B,\1}  } {  \lambda / \mu_{\B,\1}  +  ( 1 - \lambda) / \mu_{\A,\1} }  
				= \frac { \lambda \mu_{\A,\1}  } {  \lambda \mu_{\A,\1}  +  ( 1 - \lambda) \mu_{\B,\1}} , 
	\qquad
	( 1 - \hat\lambda ) = \frac { ( 1 - \lambda ) \mu_{\B,\1}  } {  \lambda \mu_{\A,\1}  +  ( 1 - \lambda) \mu_{\B,\1}} .
$$
Therefore, by Proposition \ref{p.Ehom} the homogenized stiffness matrix in the rescaled domain $\widehat\Omega$ is given by 
$$
	\widehat{\boldsymbol C}_\hom = \left(   \begin{matrix}   \hat\mu_{\hom,\1}  &  0 \\  0 &  \hat\mu_{\hom,\2}   \end{matrix} \right)  
$$
where  $\hat\mu_{\hom,\1} = \mu_{\A,\1} \, \mu_{\B,\1}$ while
$$
	\frac{ 1 } { \hat\mu_{\hom,\2} }  = \frac{ \hat \lambda } { \hat\mu_{\A,\2}  }   + \frac{ 1 - \hat \lambda } { \hat \mu_{\B,\2}  }  = 
	\frac{ \mu_{\A,\1} \, \mu_{\B,\1}  } { \lambda \mu_{\A,\1}  +  ( 1 - \lambda) \mu_{\B,\1}  } \left (  
	\frac{\lambda  }{ \mu_{\A,\2} } +  \frac{(1-\lambda) }{ \mu_{\B,\2} } 
	\right )  .
$$
Next, let $\Phi : \Omega \to \widehat\Omega$ be the bi-Lipschitz affine map  $\Phi ( x_\1, x_2 ) = ( x_\1 , \phi   (x_2) )$ where  
$$	\phi (x_2) = \beta  x_2  
 \quad \text{ for }  \quad  \beta  =  \lambda / \mu_{\B,\1}	+ (1-\lambda) / \mu_{\A,\1} 
 = \frac{\lambda \mu_{\A,\1} + (1-\lambda) \mu_{\B,\1}}{\mu_{\A,\1} \mu_{\B,\1}} . $$
Note that $\beta_n \weakstarto \beta $. 
With this change of variable the stiffness matrix $\C_\hom$ becomes $\widehat{\boldsymbol C}_\hom$, indeed (see Lemma \ref{l.gigetto}) 
\begin{gather*}
	 \mu_{\hom,\1}  = \lambda \mu_{\A,\1} + (1-\lambda) \mu_{\B,\1}   \quad \text{ and thus }  \quad  \mu_{\hom,\1} / \beta = \mu_{\A,\1} \, \mu_{\B,\1}   =  \hat\mu_{\hom,\1}  	 ,
	 \\[4pt]
	 \frac{1}{\mu_{\hom,\2}} = \frac{\lambda  }{ \mu_{\A,\2} } +  \frac{(1-\lambda) }{ \mu_{\B,\2} }  
\quad \text{ and thus }  \quad
	 \frac{1} { \beta \mu_{\hom,\2}  }  =    	\frac{ \mu_{\A,\1} \, \mu_{\B,\1}  } { \lambda \mu_{\A,\1}  +  ( 1 - \lambda) \mu_{\B,\1}  } \left (  
	\frac{\lambda  }{ \mu_{\A,\2} } +  \frac{(1-\lambda) }{ \mu_{\B,\2} } 
	\right )  = \frac{1} { \hat\mu_{\hom,\2}  }  .
\end{gather*} 
Since $\beta_\hom$ is the weak$^*$ limit of $\beta_n$, $\phi_n \to \phi$ uniformly. Let $\hat{g} = g \circ \Phi^{-1}$. Note that $\hat{g}_n$ is a bounded sequence in $H^1(\widehat\Omega \setminus K_L)$, moreover, $\hat{g}_n \to \hat{g}$ in $L^2(\partial_D \widehat\Omega)$. 

Let  $\hat{u}_{l,n}  \in \mathrm{argmin} \{ \widehat{\W}_{l,n} ( \hat{u} ) : \hat{u} \in \widehat\U_{l,n} \}$ and $u_{l,\hom} \in \mathrm{argmin} \{ \W_{l,\hom} ( u ) : u \in \U_\l \}$. 
We know that $\hat{u}_{l,n} \weakto \hat{u}_{l,\hom}$ in $H^1 ( \widehat\Omega \setminus K_l)$. Since $\hat{\mu}_{n,\1} = \hat{\mu}_{\hom,\1} = \mu_{\A, \1} \, \mu_{\B,\1}$ we get
$$
	\int_{\widehat\Omega_{(a,b)} \setminus K_l} \hat{\mu}_{\hom,\1} \, | \partial_{y_1} \hat{u}_{l,\hom} |^2 \,dy 
	\le \liminf_{n \to \infty }  
	\int_{\widehat\Omega_{(a,b)} \setminus K_l} \hat{\mu}_{n,\1} \, | \partial_{y_1} \hat{u}_{l,n} |^2  \,dy , 
$$
where $\widehat\Omega_{(a,b)} = (a,b) \times (-S, S)$. Clearly,  $\hat{u}_{l,n} = u_{l,n} \circ \Phi_n^{-1}$, where $u_{l,n} \in \mathrm{argmin} \{ \W_{l,n} ( u ) : u \in \U_\l \}$ and $\hat{u}_{l,\hom} =  u_{l,\hom} \circ \Phi^{-1}$ where $u_{l,\hom} \in \mathrm{argmin} \{ \W_{l,\hom} ( u ) : u \in \U_\l \}$. Therefore, applying the changes of variable we get 
$$
	\int_{\Omega_{(a,b)} \setminus K_l} \mu_{\hom,\1} \, | \partial_{x_1} u_{l,\hom} |^2 \,dx 
	\le \liminf_{n \to \infty }  
	\int_{\Omega_{(a,b)} \setminus K_l} \mu_{n,\1} \, | \partial_{x_1} u_{l,n} |^2  \,dx . 
$$
Using a different change of variable we can apply the above argument in such a way that 
$$
	\int_{\Omega_{(a,b)} \setminus K_l} \mu_{\hom,\2} \, | \partial_{x_2} u_{l,\hom} |^2 \,dx 
	\le \liminf_{n \to \infty}  
	\int_{\Omega_{(a,b)} \setminus K_l} \mu_{n,\2} \, | \partial_{x_2} u_{l,n} |^2  \,dx . 
$$
By Lemma \ref{l.locenconv} we have 
$$
	\int_{\Omega_{(a,b)} \setminus K_l} \mu_{n,\1} \, | \partial_{x_1} u_{l,n} |^2 +  \mu_{n,\2} \, | \partial_{x_2} u_{l,n} |^2  \,dx
	\ \to \
	\int_{\Omega_{(a,b)} \setminus K_l} \mu_{\hom,\1} \, | \partial_{x_1} u_{l,\hom} |^2 +  \mu_{\hom,\2} \, | \partial_{x_2} u_{l,\hom} |^2 \,dx . 
$$
and thus
\begin{align*}
    \limsup_{n \to \infty }  \int_{\Omega_{(a,b)} \setminus K_l} \mu_{n,\1} \, | \partial_{x_1} u_{l,n} |^2  \,dx & = \int_{\Omega_{(a,b)} \setminus K_l} \mu_{\hom,\1} \, | \partial_{x_1} u_{l,\hom} |^2 +  \mu_{\hom,\2} \, | \partial_{x_2} u_{l,\hom} |^2 \,dx  \\ & \quad -  \liminf_{n \to \infty}  
	\int_{\Omega_{(a,b)} \setminus K_l} \mu_{n,\2} \, | \partial_{x_2} u_{l,n} |^2  \,dx \\ &  \le 
	\int_{\Omega_{(a,b)} \setminus K_l} \mu_{\hom,\1} \, | \partial_{x_1} u_{l,\hom} |^2  \,dx .
\end{align*}
The proof is concluded.  \qed

\noindent {\bf Proof of Theorem \ref{t.Gconv}.}  It is enough to show that $\G_n \to \G_\hom$ pointwise in $(0,L)$.

Let $r>0$ such that $0 < 2r < l < L - 2r$.  Let $\phi_r : [0,L] \to [0,1]$ be the Lipschitz map defined by $\phi_r = 0$ in $(0, r)$ and in $(L- r,L)$, $\phi_r =1$ in $(2r,L-2r)$, $\phi_r$ is affine in $( r,2 r)$ and in $(L-2r, L- r)$. By Lemma \ref{l.Enhoriz} we can write  
$$
      \E'_n ( l) = \int_{\Omega \setminus K_l} 
      \nabla u_{l,n} \boldsymbol{C}_n \boldsymbol{E} \nabla u_{l,n}^{T} \, \phi'_r dx 
	 \quad \text{for} \quad 
    \boldsymbol{E} =   \left( \begin{matrix} -1 &   0 \\ 0 & 1 \end{matrix} \right) = \mathbf{I} +  \left( \begin{matrix} - 2 &   0 \\ 0 & 0 \end{matrix} \right) , 
$$
where $\boldsymbol{I}$ denotes the identity matrix; a similar representation holds also for $\G_\hom$. Then, splitting $\boldsymbol{E}$ as above, 
\begin{align*}
	\G_n (l) 
	& = - \E'_n (l) = - \int_{\Omega \setminus K_l}  \nabla \hat{u}_{l,n} \C_n \nabla \hat{u}^T_{l,n} \, \phi'_r   \, dx + 
	2 \int_{\Omega \setminus K_l} \mu_{n,\1} | \partial_{x_1} {u}_{l,n} |^2  \phi'_r \, dx  .
\end{align*}
By the definition of $\phi_r$, we can write 
\begin{align*}
	\int_{\Omega \setminus K_l}  \nabla {u}_{l,n} \C_n \nabla {u}^T_{l,n} \, \phi'_r \, dx 
		& = r^{-1} \int_{\Omega_{(r,2r)} \setminus K_l }   \hspace{-12pt} \nabla {u}_{l,n} \C_n \nabla {u}^T_{l,n}  \, dx
			- r^{-1} \int_{\Omega_{(L-2r,L-r)} \setminus K_l } \hspace{-12pt}  \nabla {u}_{l,n} \C_n \nabla {u}^T_{l,n} \, dx . 
\end{align*}
Hence by Lemma \ref{l.locenconv} we get 
$$
	\int_{\Omega \setminus K_l}  \nabla {u}_{l,n} \C_n \nabla {u}^T_{l,n} \, \phi'_r  \, dx
	\ \to \ 
	\int_{\Omega \setminus K_l}  \nabla {u}_{l,\hom} \C_\hom \nabla {u}^T_{l,\hom} \, \phi'_r   \, dx .
$$
In the same way we get 
$$
	\int_{\Omega \setminus K_l} \mu_{n,\1} | \partial_{x_1} {u}_{l,n} |^2  \phi'_r \, dx 
	\ \to\ 
	\int_{\Omega \setminus K_l} \mu_{\hom,\1} | \partial_{x_1} {u}_{l,\hom} |^2  \phi'_r \, dx .
$$
Hence $\G_n (l) \to \G_\hom (l)$.  \qed

\section{Generalizations} 

In this section we present a couple of simple generalizations. 

\subsection{Non-monotone boundary conditions \label{s.non-mono}}

Consider a Dirichlet boundary condition of the type $f (t) g$ where $f \in C^1( [0,T])$ with $f(0)=0$ ($f$ is non necessarily monotone). We will provide an existence and convergence result. 

Let $\hat\ell_n$ and $\hat\ell_\hom$ be the unique solutions obtained with the monotone boundary condition $u = s g $ and $s \in [0, S]$, where $S = \max \{ | f ( t) | :  t \in [0,T] \} $. 
By Theorem  \ref{t.genconv} we know that $\hat\ell_n \to \hat\ell$ pointwise in $[0,S)$ and that $\hat\ell^+ = \hat\ell_\hom$.

Following \cite[\S 12]{Negri_ACV10} we introduce the non-decreasing function $ \bar{f} (t) = \max \{ | f (\tau) | : \tau  \in [0,t]   \} $
and we define $\ell_n = \hat\ell_n \circ \bar{f} $. 
We will prove that:
\begin{itemize}
\item[1)] $\ell_n$ is a quasi-static evolution in the sense of Proposition \ref{p.existsn}, i.e., $\ell_n : [0,T) \to [L_0,L]$ is 
non-decreasing, right continuous, satisfies the initial condition $\ell_n(0) = L_0$ and Griffith's criterion, in the following form: 
\begin{itemize}
\item[i)] $G_n (t, \ell_n(t) ) \le \Gc_n (\ell_n(t))$ for every time $t \in [0,T)$;
\item[ii)]  if $G_n ( t, \ell_n (t) ) < \Gc_n (\ell_n (t)) $ then $\ell_n$ is right differentiable in $t$ and $\dot{\ell}^+_n (t) =0$; 
\item[iii)]  if $t \in J ( \ell_n)$ then $G_n (t , l ) \ge \Gc_n (l)$ for every $l \in [\ell^-_n(t) , \ell_n(t) )$;
\end{itemize}
\item [2)]  $\ell_n = \hat\ell_n \circ \bar{f}$ converge to $\ell= \hat\ell \circ \bar{f}$ pointwise in $[0,T)$ (up to subsequences); the evolution $\ell$ is non-decreasing, satisfies the initial condition and Griffith's criterion, in the following form: 
\begin{itemize}
\item[i)] $G_\hom (t, \ell (t) ) \le \Gc_\eff (\ell (t))$ for every time $t \in [0,T)$;
\item[ii)]  if $G_\hom ( t, \ell  (t) ) < \Gc_\eff (\ell (t)) $ then $\ell $ is right differentiable in $t$ and $\dot{\ell}^+ (t) =0$; 
\item[iii)]  if $t \in J ( \ell )$ then $G_\hom (t , l ) \ge \Gc_\eff (l)$ for every $l \in [ \ell^- (t) , \ell^+ (t) )$.
\end{itemize}
\end{itemize}

Let us prove 1). 
Clearly $\ell_n$ is non-decreasing, since both $\hat\ell_n$ and $\bar{f}$ are non-decreasing, and $\ell_n (0) = \hat\ell_n \circ f (0) = \hat\ell_n (0) = L_0$. The right continuity of $\hat\ell_n$ together with the monotonicity and continuity of $\bar{f}$ imply the right continuity of $\ell_n$. 
It remains to check Griffith's criterion, but in this case the representation formula \eqref{e.ell+} does not hold. 
Before proceeding, note that for $s = \bar{f} (t)$  we have
\begin{equation} \label{e.mono09}
    \hat\ell_n^- (s) \le \ell_n^- (t) \le  \ell_n^+(t) \le \hat\ell_n (s) .
\end{equation}

\begin{itemize}
\item[i)] We know that $G_n (s , \hat\ell_n(s) ) \le \Gc_n (\hat \ell_n(s))$ for every $s \in [0,S)$; writing  $s = \bar{f}(t)$ and $\ell_n = \hat\ell \circ \bar{f}$ yields
$$ \bar{f}^{\,2}(t) \G_n ( \ell_n (t) ) = \bar{f}^{\,2}(t)  \G_n ( \hat\ell_n \circ \bar{f} (t) ) = G_n ( \bar{f} (t) , \hat\ell_n \circ \bar{f} (t))   \le \Gc_n (\hat\ell_n \circ \bar{f} (t) ) =  \Gc_n (\ell_n (t) )  . $$
Since $f^2 \le \bar{f}^{\,2} $ we get 
$$
	G_n ( t, \ell_n (t) ) = f^2 (t) \G_n ( \ell_n (t) ) \le \bar{f}^{\,2} (t) \G _n ( \ell_n(t))  \le \Gc_n (\ell_n (t) )  . 
$$

\item[ii)] If $| f(t) | = \bar{f} (t)$ then $G_n ( t, \ell_n (t) ) < \Gc_n (\ell_n (t)) $ reads
$$
	\bar{f}^{\,2} (t) \G_n ( \hat\ell_n \circ \bar{f} (t) )  = f^{2} (t) \G_n ( \hat\ell_n \circ \bar{f} (t) ) < \Gc_n (\hat\ell_n \circ \bar{f} (t))
$$
which, for $s = \bar{f}(t)$, gives $G_n (s ,  \hat\ell_n(s) ) <  \Gc_n(\hat\ell_n (s))$. Since $\bar{f}$ is continuous and monotone non-decreasing, it turns out that 
$$ 0 \le \dot\ell^+_n (t) \le  \lim_{h \to 0^+} \frac{\hat\ell_n(s+h) - \hat{\ell}_n(s) }{h}  = 0 .  $$ 

If $| f(t) | < \bar{f} (t)$ then $\bar{f} (t) = f ( t^*)$ for some $t^* < t$. Moreover, by the continuity of $f$ there exists $\delta>0$ such that $\bar{f} (t') = f ( t^*)$ for every $| t' - t | < \delta$. Hence $\ell_n (t') = \hat\ell_n \circ \bar{f} (t') = \ell_n \circ f (t^*)$ for every $| t' - t | < \delta$; hence $\ell_n$ is constant in a neighborhood of $t$. 

 \item[iii)] Let $t   \in J ( \ell_n)$. First, note that $ \bar{f} ( t ) = | f(t )|$; indeed, if $ \bar{f} ( t ) > | f(t )|$ then, repeating the argument above, $\ell_n$ would be constant in a neighborhood of $t $ and thus $t \not\in J (\ell_n)$. 
 By \eqref{e.mono09} 
 it turns out that $ s    \in J (\hat\ell_n)$. 
Hence, $G_n ( s , l ) \ge  \Gc_n (l)$ for every $l \in [ \hat\ell^-_n (s ) , \hat\ell^+_n (s  ) )$; the substitution $s = \bar{f}(t )$ and \eqref{e.mono09} lead to 
 $$   G_n (t, l) = f^2 (t) \G_n (l) = \bar{f}^{\,2}(t) \G_n (l ) = G_n (s, l) \ge \Gc_n (l)    $$ 
 for every $l \in [\ell^-_n (t) , \ell^+_n(t) )$.
\end{itemize}

\separe


Let us prove 2). Before proceeding, note that in general $ \ell^+ = (\hat\ell \circ \bar{f})^+ \neq \hat\ell^+ \circ \bar{f} = \hat\ell_\hom \circ \bar{f}$, therefore we cannot employ directly the properties of $\hat\ell_\hom$. Instead, we first check Griffith's criterion for $\hat\ell$.  We know that $\hat\ell^+ = \hat\ell_\hom$ and thus $J ( \hat\ell ) = J ( \hat\ell_\hom)$, moreover 
$\hat\ell^- (s) = \hat\ell_\hom^- (s)$ and $\hat\ell^+ (s) = \hat\ell_\hom (s)$ for every $s \in J ( \hat\ell )$, while $\hat\ell (s)  = \hat\ell_\hom(s)$ for every  $s \in [0,S) \setminus J ( \hat\ell)$. 
\begin{itemize}
\item[iii)]  By Griffith's criterion for $\hat\ell_\hom$, if $s \in J ( \hat\ell_\hom ) = J (\hat\ell) $ then $G_\hom (s , l ) \ge \Gc_\eff (l)$ for every $l \in [ \hat\ell^-_\hom (s) , \hat\ell_\hom (s) ) =  [ \hat\ell^- (s) , \hat\ell^+(s) )$.
\item[ii)]  If  $G_\hom ( s , \hat\ell  (s) ) < \Gc_\eff ( \hat\ell (s)) $ then either $s \not\in J ( \hat\ell)$ or $s \in J(\hat\ell)$ and $\hat\ell (s) = \hat\ell_\hom (s)$ (otherwise the opposite inequality would hold). In both the cases $\hat\ell (s) = \hat\ell_\hom(s)$ and their right derivatives coincide. Therefore, $\hat\ell $ is right differentiable in $s$ and its right derivative vanishes. 
\item[i)] We employ the definition of $\Gc_\eff$. If $\G_\hom ( \hat\ell (s) ) =0$ or $\Gc_\eff ( \hat\ell (s) ) = +\infty$ there is nothing to prove. Otherwise, by pointwise convergence
$$
	\frac{ G_\hom ( s , \hat\ell (s)) }{\Gc_\eff ( \hat\ell(s)) }  \le \liminf_{n \to \infty} \frac{ G_n ( s , \hat\ell_n (s)) }{\Gc_n ( \hat\ell_n (s)) } \le 1 , 
$$
which implies that $G_\hom (s , \hat\ell  (s) ) \le \Gc_\eff ( \hat\ell (s))$ for every $s \in [0,S)$.
\end{itemize}

At this point, employing the change of variable $s = \bar{f}(t)$ and arguing as in point 1) it follows that Griffith's criterion holds for $\ell = \hat\ell \circ \bar{f}$.

\subsection{Plane-strain setting} \label{s.pss}

The abstract analysis developed in \S\ref{s.e} can be easily extended to the plane strain setting with minor changes in the proofs. On the contrary, the generalization of the example and of the results provided in \S\ref{explicit} and \S\ref{s.oro} is not straightforward, indeed, the computations provided in that sections are tailored to the anti-plane setting, being based on the specific form of the homogenized energy and of the energy release. 

\medskip
We give a brief outline. Consider the same geometry for the reference domain, the cracks and the layers. Consider a boundary condition of the form $f(t) g$ where $g \in H^1(\partial_D \Omega ; \mathbb{R}^2)$. The condensed elastic energy $\E_n$ is again of class $C^1( (0,L) \setminus \Lambda_n)$, thus, as in Proposition \ref{p.existsn}, there exists a unique, right-continuous evolution $\ell_n$ which satisfies Griffith's criterion  and which can be characterized by $\ell_n (t) = \inf \{ l \in [L_0,L] : G_n (t , l) < \Gc_n (l) \}$. 

Next, by classical results on H-convergence, see e.g.~\cite{FrancfortMurat_ARMA86}, the linear elastic energy $E_n$ converge to a linear elastic energy $E_\hom$ with constant coefficients, in particular the condensed (or reduced) energy $\E_\hom$ is of class $C^1(0,L)$. 

By Helly's Theorem, upon extracting a subsequence, there exists a limit evolution $\ell$.
After defining the effective toughness, as in \eqref{e.Ghomgen}, and arguing as in Theorem \ref{t.genconv} it follows that $\ell^+ = \ell_\hom$, where $\ell_\hom$ is again the unique  right-continuous evolution which satisfies Griffith's criterion.

\appendix 
\section{Balanced viscosity solutions and Griffith's criterion}  \label{a.A}

For sake of completeness, we provide the notion of balanced viscosity solutions for the rate independent system $(\E_n ,\mathcal{R}_n)$, making reference to \cite{MielkeRossiSavare_COCV12} with a few differences due to irreversibility and lack of regularity in $\Lambda_n$. Moreover we briefly show that the evolution $\ell_n$ (provided by Proposition \ref{p.existsn}) is a balanced viscosity solution and, vice versa, that each right continuous balanced viscosity solution satisfies Griffith's criterion. A similar reasoning applies to $\ell_\hom$. 

Following \cite{MielkeRossiSavare_COCV12} a right continuous evolution $\gamma_n$ is a balance viscosity evolution for the rate independnent system $(\E_n , \mathcal{R}_n)$ if it satisfies the following conditions for every  time $t \in [0,T]$: 
\begin{equation} \label{e.BV1}
	\partial_{l'} \mathcal{R}_n ( \gamma_n(t) , 0 ) + \partial_{l} E_n ( t , \gamma_n (t)) \ni 0 
\end{equation}
\begin{equation} \label{e.BV2}
	E_n (t, \gamma_n(t)) =  E_n ( 0, L_0)  + \int_0^t P_n ( t, \gamma_n (t)) \, dt - \mathrm{Var}_{\mathfrak{f}_n} ( \gamma_n ; [0,t] ) . 
\end{equation} 
We recall that the state dependent dissipation $\mathcal{R}_n : [L_0, L] \times \mathbb{R} \to [0,+\infty]$ is given by 
$$
	\mathcal{R}_n ( l , l' ) =  \begin{cases} G^c_n (l) \, l'   &  l' \ge  0  , \\  +\infty & l' < 0 .  \end{cases}
$$
We denote by $\partial _{l'}\mathcal{R}_n $ the subdifferential of the dissipation $\mathcal{R}_n ( l , \cdot)$ and by $\partial_l E_n ( t , l)$ the partial derivative of the reduced elastic energy. Then $\partial_{l'} \mathcal{R}_n ( l , 0 )  = ( -\infty , G^c_n(l)]$ while $\partial_l E_n ( t , l) = - G_n ( t , l)$, if $l \not\in \Lambda_n$. If $l \in \Lambda_n$ the energy $\E_n $ is not differentiable, therefore following \eqref{e.extGn} we set $\partial_l E_n ( t , l) = - G_n ( t , l)$.

Next, let us introduce the dissipation ``distance'' $\mathsf{dist}_n : [L_0,L] \times [L_0,L]  \to [0, +\infty] $ given by
$$	\mathsf{dist}_n  ( l_0 , l_\1)  = \inf \left\{ \int_0^1 \mathcal{R}_n ( l(r) , l'(r) ) \, dr :  l \in AC (0,1) , \, l (0) = l_0 , \, l(1)=l_\1   \right\}  	$$
and the Finsler dissipation cost $\Delta_{\mathfrak{f}_n } : [0,T] \times [L_0,L] \times [L_0,L]   \to [0,+\infty]$ given by 
$$
	\Delta_{\mathfrak{f}_n } ( t, l_0 , l_\1)  =  \inf \left\{ \int_0^1 \mathfrak{f}_n ( t , l(r) , l'(r) ) \, dr :  l \in AC (0,1) , \, l (0) = l_0 , \, l(1)=l_\1   \right\}  , 
$$
where
\begin{align*}
	\mathfrak{f}_n ( t , l(r) , l'(r) ) & =  \mathcal{R}_n ( l(r) , l'(r) ) + \mathrm{dist} \big(   \!-\! \partial_l E ( t , l (r) )  , \partial_{l'} \mathcal{R}_n ( l (r) , 0 ) \big)  | \l' (r) |  .
\end{align*}
Finally, let us introduce the notation
$$\mathrm{Var}_{\mathfrak{f}_n} ( \gamma_n ; [0, t]) = \mathrm{Var}_{\mathsf{dist}_n} ( \gamma_n ; [0, t]) - \mathrm{Jmp}_{\mathsf{dist}_n} ( \gamma_n ; [0, t]) + \mathrm{Jmp}_{\mathfrak{f}_n}  ( \gamma_n ; [0, t]) ,  $$ 
where 
\begin{align*}
\mathrm{Var}_{\mathsf{dist}_n} ( \gamma_n ; [0, t]) & =  \sup \left\{ \sum_{i=0}^{m-1}  \mathsf{dist_n} (  \gamma_n (t_{i+1}) , \gamma_n ( t_i )   ) : 0 = t_0 \le ... \le t_m = t  \right\} , \\
\mathrm{Jmp}_{\mathsf{dist}_n} ( \gamma_n ; [0, t]) & =  \sum_{\tau \in J ( \gamma_n) \cap [0,t]} \mathsf{dist}_n ( \gamma_n (\tau) , \gamma_n^-(\tau))  , \\
\mathrm{Jmp}_{\mathfrak{f}_n} ( \gamma_n ; [0, t]) & =  \sum_{\tau \in J ( \gamma_n) \cap [0,t]} \Delta_n ( \tau , \gamma_n (t) , \gamma_n^-(t))  .
\end{align*}
(Remember that here $\gamma_n$ is right continuous). 

First, let us check that the evolution $\ell_n$ provided by Proposition \ref{p.existsn} is a balanced viscosity solution. 
The inequality $ G_n ( t, \ell_n (t) ) \le G^c_n ( \ell_n(t)) $ reads 
$$
	- \partial E_n ( t , \ell_n (t)) \in (- \infty , \Gc_n ( \ell_n(t)) ] = \partial_{l'} \mathcal{R}_n ( \ell_n(t) , 0 )  ,
$$
which is \eqref{e.BV1}. Next, let us check that the energy identity \eqref{e.enid2} implies \eqref{e.BV2}, showing that 
\begin{align*}
	\mathrm{Var}_{\mathfrak{f}_n} ( \ell_n ; [0,t] ) & = \mathrm{Var}_{\mathsf{dist}_n} ( \ell_n ; [0, t]) - \mathrm{Jmp}_{\mathsf{dist}_n} ( \ell_n ; [0, t]) + \mathrm{Jmp}_{\mathfrak{f}_n}  ( \ell_n ; [0, t]) ,  \\
	& = D_n ( \ell_n(t))  - \hspace{-12pt} \sum_{\tau \, \in \, J(\ell_n) \cap [0,t]}  \llbracket F_n ( \tau , \ell_n (\tau)) \rrbracket .
\end{align*}
First of all, note that in our irreversible setting it is not restrictive to consider $l' \ge 0$ in the definition of both $\mathsf{dist}_n$ and 
$\Delta_{\mathfrak{f}_n}$. 
By explicit computations, it is easy to check that
$$
\mathsf{dist}_n ( l_0 , l_\1)  = \begin{cases}
	D_n ( l_\1) - D_n (l_0)  & \text{if } l_0 \le l_\1 , \\
	+\infty & \text{otherwise.} 
	\end{cases}
$$
As a consequence, by monotonicity and right continuity of $\ell_n$, 
\begin{align*}
	\mathrm{Var}_{\mathsf{dist}_n} ( \ell_n ; [0, t]) & = D_n ( \ell_n (t) )  ,  \\
	\mathrm{Jmp}_{\mathsf{dist}_n} ( \ell_n ; [0, t]) & =  \sum_{\tau \in J ( \ell_n) \cap [0,t]} \mathsf{dist}_n ( \ell_n (t) , \ell_n^-(t)) = \sum_{\tau \in J ( \ell_n) \cap [0,t]}  \jump{ D_n (\ell_n (\tau)) }  .
\end{align*}
Next, note that 
$$
	\mathrm{dist} \big(   \!-\! \partial_l E ( t , l)  , \partial_{l'} \mathcal{R}_n ( l , 0 ) \big) = [ G (t, l) - \Gc_n (l) ]_+ 
$$
By condition iii) of Griffith's criterion for $\tau \in J ( \ell_n) $ we have $G_n (\tau , l ) - \Gc_n (l) \ge 0$ for every $l \in [ \ell_n ^-(\tau) , \ell_n (\tau))$ and thus
$$ 	\mathfrak{f}_n ( \tau , l(r) , l'(r) )  =  G_n (\tau , l (r) ) \,   \l' (r)    $$ 
for every monotone $l \in AC(0,1)$ with $l(0) = \ell_n^- (\tau)$ and $l(1) = \ell_n (\tau)$. 
It follows easily that 
$$
 	\Delta_n ( \tau , \ell_n (\tau) , \ell_n^-(\tau)) = -  \jump{ E_n (\tau, \ell_n (\tau)) } 
$$
and thus 
$$
\mathrm{Jmp}_{\mathfrak{f}_n} ( \ell_n ; [0, t])  = \sum_{\tau \in J ( \ell_n) \cap [0,t]}  -  \jump{ E_n (\tau, \ell_n (\tau)) } .
$$ 
Conversely, let $\gamma_n$ be a (right continuous) balanced viscosity solution. If $\gamma_n$ satisfies \eqref{e.BV1} then 
$$
    (- \infty , \Gc_n ( \gamma_n(t) ) ] - G_n ( t , \gamma_n (t)) \ni 0 ,
$$
that is $G_n ( t , \gamma_n (t)) \le \Gc_n (\gamma_n(t))$. Let us check that $G_n ( \tau , l ) \ge \Gc_n (l)$ for $\tau \in J (\gamma_n)$ and every $ l \in [ \gamma_n^- (\tau) , \gamma_n (\tau)]$. In the interval $[ \tau -\delta , \tau + \delta]$ (for some $\delta>0$) the energy identity \eqref{e.BV2} reads
$$
E_n (\tau + \delta , \gamma_n(\tau+\delta)) -  E_n ( \tau -\delta , \gamma_n (\tau - \delta) ) = \int_{\tau - \delta}^{\tau + \delta} P_n ( t, \gamma_n (t)) \, dt - \mathrm{Var}_{\mathfrak{f}_n} ( \gamma_n ; [\tau-\delta,\tau+\delta] ) . 
$$
Hence, taking the limit as $\delta \to 0^+$, we get 
\begin{equation} \label{e.enidin}
      E_n ( \tau , \gamma_n^- (\tau) ) - E_n ( \tau , \gamma_n (\tau)) 
      = \int_{\gamma_n^- (\tau)}^{\gamma_n (\tau)} G_n ( \tau , s) \, ds  =  \lim_{\delta \to 0^+} \mathrm{Var}_{\mathfrak{f}_n} ( \gamma_n ; [\tau-\delta,\tau+\delta] ) .
\end{equation}
It is easy to check that
\begin{align*}
    \lim_{\delta \to 0^+} \mathrm{Var}_{\mathfrak{f}_n} ( \gamma_n ; [\tau-\delta,\tau+\delta] ) 
    			 & = \Delta ( \tau , \gamma_n^- (\tau) , \gamma_n (\tau)) = \int_{\gamma_n^- (\tau)}^{\gamma_n (\tau)} \Gc_n (  s) + [ G_n ( \tau , s) - \Gc (s)     ]_+  ds  .
\end{align*}
Hence  by   \eqref{e.enidin} 
$$
        \int_{\gamma_n^- (\tau)}^{\gamma_n (\tau)} G_n ( \tau , s) - \Gc_n ( s) \, ds  = \int_{\gamma_n^- (\tau)}^{\gamma_n (\tau)} [ G_n ( \tau , s) - \Gc (s)     ]_+  ds ,
$$
which in turn implies $ G_n ( \tau , s) - \Gc_n (s) \ge 0$ for a.e.~$s \in (\gamma_n^- (\tau) , \gamma_n (\tau))$. By the regularity of energy release and toughness we conclude that $ G_n ( \tau , s) \ge \Gc_n (s) $ holds for every $s \in [\gamma_n^- (\tau) , \gamma_n (\tau))$. It remains to show that $G_n (\tau , \gamma_n(\tau)) < \Gc_n ( \gamma_n(\tau))$ implies $\dot{\gamma}_n^+ (\tau) =0$. Let us first assume that $\ell_n (\tau) \not\in \Lambda_n$. Then, by continuity of $G_n$ and $\Gc_n$ it follows that $G_n ( t , l) <   \Gc_n (l)$  for every $(t,l)$ such that $| \tau -t | < \delta$ and $| l - \ell_n(\tau) | < \delta$. 
By right continuity we can further assume that $| \gamma_n (t) - \gamma_n (\tau) | < \delta $ for $ t \in ( \tau , \tau +\delta)$.  The fact that $G_n ( t , l) <   \Gc_n (l)$ implies that there are no jumps in $ t \in ( \tau , \tau +\delta)$, otherwise iii) (which we already proved) would be violated. Hence for $ t \in ( \tau , \tau +\delta)$ the energy identity reads
$$
      E_n (t, \gamma_n(t) ) = E_n ( \tau, \gamma_n (\tau)) + \int_{\tau}^{t} P_n (t, \gamma_n (t)) \, dt - \mathrm{Var}_{\mathsf{dist}_n} ( \gamma_n ; [\tau, t]) 
$$
Writing 
\begin{align*}
E_n (t, \gamma_n(t) ) - E_n ( \tau, \gamma_n (\tau)) & =  \int_{\tau}^{t} \partial_t E_n (t, \gamma_n (t)) \, dt + \int_\tau^t \partial_l E_n (t, \gamma_n (t))   \, d\gamma_n (t)  \\
	\int_{\tau}^{t} P_n (t, \gamma_n (t)) \, dt & = \int_{\tau}^{t} \partial_t E_n (t, \gamma_n (t))  \, dt \\
	\mathrm{Var}_{\mathsf{dist}_n} ( \gamma_n ; [\tau, t]) & = \int_{\tau}^t \Gc_n ( \gamma_n (t) ) \, d\gamma_n (t) 
\end{align*}
it follows that 
$$
        \int_\tau^t G_n (t, \gamma_n (t))   \, d\gamma_n (t)  =  \int_{\tau}^t \Gc_n ( \gamma_n (t) ) \, d\gamma_n (t) . 
$$
This is a contradiction with $G_n (t, \gamma_n (t)) < \Gc_n ( \gamma_n (t) ) $ unless $d\gamma_n = 0$, which implies the thesis. It $\gamma_n (\tau) \in \Lambda_n$ then there exists a sequence $l_m \searrow \gamma_n (\tau)$ such that $l_m \not\in \Lambda_n$ and $G_n (\tau, l_m ) \to G_n ( \tau, \gamma_n (\tau) ) < \Gc_n (\gamma_n (\tau) ) = \Gc_n (\l_m)$. Since $l_m \not\in \Lambda_n$, arguing as above it follows that $\gamma_n (t) < l_m$ in a right neighborhood of $\tau$, uniformly with respect to $m \in \mathbb{N}$. Hence, taking the infimum with respect to $m \in \mathbb{N}$, we get  $\gamma_n (t) \le \gamma_n (\tau)$; by monotonicity, $\gamma_n$ is constant in a right neighborhood of $\tau$.

\section{Homogenization of displacement, stress, and energy} \label{a.B}

This section contains classic results on the homogenization of Dirichlet energies (see e.g.~\cite{MuratTartar_97} or \cite{BraidesDefranceschi_98} and the references therein). We provide a short proof, adapted to the case of horizontal layers. We recall that
$$
 \C_{\hom} =  \left(\begin{matrix}  \mu_{\hom, \1} & 0  \\ 0 & \mu_{\hom, \2}  \end{matrix}  \right) ,
$$ 
where $1 / \mu_{\hom,\1} = \lambda / \mu_{\A,\1} + (1-\lambda) / \mu_{\B,\1}$ and 
$\mu_{\hom, \2} =  \lambda \mu_{\A,\2} + ( 1- \lambda ) \mu_{\B,\2} $. For $l \in [0,L]$ we denote 
$$	u_{l,n} \in \argmin \{ \W_{l,n} ( u ) :  u \in \U_l \}  , \qquad  u_{l,\hom}  \in \argmin \{ \W_{l,\hom} ( u ) :  u \in \U_l \} .  $$

\begin{proposition} \label{p.Ehom} Let $l \in [0,L]$. Then $u_{l,n} \weakto  {u}_{l,\hom}$ in $H^1( {\Omega} \setminus K_l)$ and $\C_n \nabla  {u}_{l,n}^T \weakto  \C_\hom \nabla  {u}_{l,\hom}^T$ in $L^2 (   {\Omega} \setminus K_l ; \R^2)$.
\end{proposition}

\medskip

\proof 
By uniform coercivity $ {u}_{l,n}$ is a bounded sequence in $H^1(  {\Omega} \setminus K_l)$ and thus  $ \C_n \nabla  {u}_{l,n}^T $ is a bounded sequence in $L^2 (  {\Omega} \setminus K_l ; \R^2)$.
Then, up to non-relabeled subsequences, $ {u}_{l,n} \weakto  {u}_{l,\infty}$ in $H^1 (   {\Omega} \setminus K_l)$ while $  \C_n \nabla  {u}_{l,n}^T  \weakto \stress_{l,\infty}$ in $L^2 (   {\Omega} \setminus K_l ; \R^2)$. 
We will prove that 
\begin{equation} \label{e.stresshom}
	\stress_{l,\infty} =  \C_\hom \nabla  {u}_{l,\infty}^T .
\end{equation} 
From this identity it follows that $ {u}_{l,\infty} =  {u}_{l,\hom}$; indeed, $u_{l,\infty} \in \U_l$  and, passing to the limit in the variational formulation 
$$
	\int_{  {\Omega} \setminus K_l} \nabla v  \C_n \nabla  {u}_{l,n} ^T\, dx = 0 
		\quad \text{for every $v \in   \V_l$, }
$$
we easily get that $ {u}_{l,\infty}$ is the unique solution of the variational problem 
$$
	\int_{  {\Omega} \setminus K_l} \nabla v  \C_\hom \nabla  {u}_{l,\infty} ^T\, dx = 0 
		\quad \text{for every $v \in   \V_l$.}
$$

The idea to prove \eqref{e.stresshom} is to show that $(\stress_{l,\infty})_i =   \mu_{\hom, i} \partial_{x_i}  {u}_{l,\infty}$ for $i=1,2$,  by passing to the limit in the constitutive relation, i.e.,  
 $(\stress_{l,n} )_{i} =   \mu_{n,i} \partial_{x_{i}}  {u}_{l,n}$ for $i=1,2$. 

\separe

{\bf I.} Let $i=1$. For every $\varphi \in C^\infty_c (   {\Omega} \setminus K_l)$ we have
$$
            \int_{  {\Omega} \setminus K_l} ( \stress_{l,n})_{\1} \varphi /   \mu_{n,\1} \, dx  = \int_{  {\Omega} \setminus K_l} ( \partial_{x_{1}}  {u}_{l,n} ) \varphi \, dx .
$$
By weak convergence we can pass to the limit in the right hand side, getting
$$
   \int_{  {\Omega} \setminus K_l} ( \partial_{x_{1}}  {u}_{l,n} ) \varphi \, dx  \to \int_{  {\Omega} \setminus K_l} ( \partial_{x_{1}}  {u}_{l,\infty} ) \varphi \, dx .
$$
Let us prove that 
\begin{equation} \label{e.vuen}
            \int_{  {\Omega} \setminus K_l} ( \stress_{l,n})_{\1} \varphi /   \mu_{n,\1} \, dx  \to 
            \int_{  {\Omega} \setminus K_l} ( \stress_{l,\infty})_{\1} \varphi /   \mu_{\hom, \1} \, dx . 
\end{equation}
Let $I= (-H,0) \cup (0,H)$ and write 
\begin{align} 
               \int_{  {\Omega} \setminus K_l} ( \stress_{l,n})_{\1} \varphi /   \mu_{n,\1} \, dx
            & = \int_{(0,L)}  \Big( \int_I ( \stress_{l,n}  )_{\1}  (x_{\1} , x_{\2} )  \varphi ( x_{\1} , x_{\2}) \, d x_{\2}  \Big) \Big/    \mu_{n,\1} ( x_{\1})   \, d x_{\1} \nonumber \\ & = \int_{(0,L)}  s_n (x_{\1}) /   \mu_{n,\1} ( x_{\1} ) \, d x_{\1} ,  \label{e.vueno}
\end{align}
where we have introduced the auxiliary function $s_n$ given by
$$
       s_n (x_{1}) 
       = \int_I ( \stress_{l,n}  )_{\1}  (x_{\1} , x_{\2} )  \varphi ( x_{\1} , x_{\2}) \, d x_{\2}  ,
       	\quad \text{for a.e.~$x_\1 \in (0,L)$.}
$$
We claim that $s_n \to s_\infty$ strongly in $L^1 (0,L)$, where 
$$
       s_\infty (x_{\1}) 
       = \int_I  (\stress_{l,\infty})_{\1} ( x_{\1} , x_{\2} ) \, \varphi ( x_{\1} , x_{\2} ) \, d x_{\2} 
$$
from which \eqref{e.vuen} follows from \eqref{e.vueno} because $1/   \mu_{\hom,\1}$ is the weak$^*$ limit of $1/  \mu_{n,\1}$ in $L^\infty(0,L)$.

Denote by $\sigma_n$ the (vector valued) function $\sigma_{n} (x_{\1} ) =  (\stress_{l,n})_{\1} ( x_{\1} , \cdot )$; since $\stress_n$ is bounded in $L^2(   {\Omega} \setminus K_l)$, it turns out that $\sigma_n$  is bounded in $L^2( 0, L; L^2(I) )$. We claim that $\sigma_n'$ is bounded in $L^2( 0, L ; H^{-1} (I) )$. Let us start writing 
\begin{align*}
\| \sigma_n' \|_{L^2(0,L;H^{-1}(I) )} = \sup \int_{(0,L)} \langle \sigma'_n , z \rangle  \, d x_{\1} ,
\end{align*}
where the supremum is taken over all the (vector valued) functions $z$ with $\| z \|_{L^2(0,L;H^1_0(I)  )} \le 1$ and $\langle \cdot , \cdot \rangle$ denotes the duality between $H^{-1} (I) $ and $H^1_0 (I) $. The set of functions of the form $z = \phi \psi $ for $\phi \in C^\infty_c (0,L)$ and $ \psi \in C^\infty_c ( I )$ is dense in $L^2(0,L;  H^1_0 (I) )$. For any such $z$ we can write 
\begin{align*}
	\int_{(0,L)} \langle \sigma'_n , z \rangle  \, d x_{\1} & = \int_{(0,L)} \langle \sigma'_n , \psi \rangle \phi \, dx_\1 
	= - \int_{(0,L)} \langle \sigma_n , \psi \rangle  \phi' \, dx_\1 \\ & = - \int_{(0,L)} \int_{I} \sigma_n (x_{\1}, x_{\2}) \psi (x_{\2})  \partial_{x_1} \phi (x_{\1}) \, dx_\2 \, dx_\1 \\
		& =  -\int_{ {\Omega} \setminus K_l} \stress_{l,n} \partial_{x_1} z \, dx \le \| \stress_{l,n} \|_{L^2 ( {\Omega} \setminus K_l)} .
\end{align*}
Then, by Schauder's theorem and Aubin-Lions lemma there exists a subsequence (not relabeled) such that $ \sigma_n \to \sigma_\infty$ (strongly) in $L^2(0,L; H^{-1} (I)  )$. As the stresses are $L^2$ functions, writing  
$$
       s_n (x_{\1}) =   \!\! \phantom{|}_{H^{-1}(I) } \langle \sigma_n  ( x_{\1} ) , \varphi ( x_\1 , \cdot  )  \rangle_{H^1_0(I) } \quad \text{and} \quad 
       s_\infty (x_{\1}) =   \!\! \phantom{|}_{H^{-1}(I) } \langle \sigma_\infty  ( x_{\1} ) , \varphi ( x_\1, \cdot  )  \rangle_{H^1_0(I) } 
$$
it follows that $s_n \to s$ in $L^1(0,L)$.

\separe

{\bf II.} Let $i=2$. Since $\mu_{n,\2}$ is independent of $x_\2$, for every $\varphi \in C^\infty_c (   {\Omega} \setminus K_l)$ we have
$$
            \int_{  {\Omega} \setminus K_l} ( \stress_{l,n})_{\2} \, \varphi \, dx  
            = \int_{  {\Omega} \setminus K_l} (  \partial_{x_{2}} (  \mu_{n,\2}   {u}_{l,n} ) ) \varphi \, dx = - \int_{  {\Omega} \setminus K_l}   \mu_{n, \2}   {u}_{l,n} \partial_{x_{2}}  \varphi \, dx .
$$
We have $   \mu_{n,\2}  \, \partial_{x_{2}} \varphi \weakstarto    \mu_{\hom, \2}  \, \partial _{x_{2}} \varphi$ in $L^\infty ( {\Omega} \setminus K_l)$ and $u_{l,n} \to  {u}_{l,\infty}$ in $L^1 ( {\Omega} \setminus K_l)$ (by compact embedding), therefore 
\begin{align*}
- \int_{  {\Omega} \setminus K_l}     \mu_{n,\2}    {u}_{l,n}   \partial_{x_{2}} \varphi  \, dx
\to
- \int_{  {\Omega} \setminus K_l}     \mu_{\hom,\2}  {u}_{l,\hom}  \partial_{x_{2}} \varphi  \, dx
& = 
\int_{  {\Omega} \setminus K_l}     \mu_{\hom,\2}  (\partial_{x_{2}}  {u}_{l,\hom} ) \varphi \, dx  \\
& = \int_{  {\Omega} \setminus K_l}    ( \stress_{\l,\hom})_{\2} \varphi \, dx  .
\end{align*}
The proof is concluded. \qed

By Proposition \ref{p.Ehom} and standard arguments it follows the convergence of the reduced energies, i.e., $ \E_n \to \E_\hom$ pointwise in $[0,L]$.  
Further, it is not difficult to prove also the following result.

\begin{corollary} \label{c.enhom} Let $l_n \in [0,L]$ such that $l_n \to l$.~Then $\nabla u_{l_n,n} \weakto \nabla {u}_{l,\hom}$ and $\C_n \nabla  {u}_{l_n,n}^T \weakto  \C_\hom \nabla  {u}_{l,\hom}^T$ in $L^2 (   {\Omega} ; \R^2)$. Moreover $ \E_n ( l_n) \to \E_\hom (l) $. 
\end{corollary}

As  by-product of the previous Corollary the energy $\E_n$ converge uniformly to $\E_\hom$ in $[0,L]$.


\end{document}